%% file: article_inputs.tex
\documentclass[a4paper, 12pt]{article}

\input{Commands}

\RequirePackage[authoryear]{natbib}

\newcommand{\SR}[2]{\textcolor{gray}{#1}\textcolor{blue}{#2}}

\date{}
\begin{document}

\title{Motif-based tests for bipartite networks}

\author{
  Sarah Ouadah$^{1,2}$, Pierre Latouche$^2$, St\'ephane Robin$^{1,3}$ \\
  \footnotesize{($^1$) UMR MIA-Paris, AgroParisTech, INRAE, Universit\'e Paris-Saclay, 75005 Paris, France} \\
  \footnotesize{($^2$) MAP5, UMR CNRS 8145, Universit\'{e} de Paris, 75006 Paris, France}  \\
  \footnotesize{($^3$) CESCO, UMR 7204, MNHN - CNRS - UPMC, Paris, France}  
  }
  
\maketitle


\abstract{
Bipartite networks are a natural representation of the interactions between entities from two different types. The organization (or topology) of such networks gives insight to understand the systems they describe as a whole. Here, we rely on motifs which provide a meso-scale description of the topology. Moreover, we consider the bipartite expected degree distribution (B-EDD) model which accounts for both the density of the network and possible imbalances between the degrees of the nodes. Under the B-EDD model,  we prove the asymptotic normality of the count of any given motif, considering sparsity conditions. We also provide close-form expressions for the mean and the variance of this count. This allows to avoid computationally prohibitive resampling procedures. Based on these results, we define a goodness-of-fit test for the B-EDD model and propose a family of tests for network comparisons. We assess the asymptotic normality of the test statistics and the power of the proposed tests on synthetic experiments and illustrate their use on ecological data sets.
}

\section{Introduction}  \label{sec:intro}

Bipartite interaction networks are used to represent a diverse range of interactions in various fields such as biology, ecology, sociology or economics. For instance, in ecology, bipartite graphs depict interactions between two groups of species such as plants and pollinators \citep[see e.g.][]{SCB19, DFT20} or host and parasites \citep[see e.g.][]{VPD08,DCBA20}, in agroethnology, they may involve interactions between farmers and  crop species \citep[see][]{TVB15} and in economics, country-product trades as signals of the 2007-2008 financial crisis \citep[see][]{SDGS16}. Formally, a bipartite interaction network can be viewed as a bipartite graph, the nodes of which being individuals pertaining to two different groups, and an edge between two nodes being present if these two individuals interact. In the sequel, the two types of nodes will be referred to as top nodes and bottom nodes, respectively. Characterizing the general organization of such a network, namely its topology, is key to understand the behavior of the system as a whole. 

\bigskip
The topology of a network can be studied at various scales. Micro-scale analyses  typically focus on the degree of each node, the betweenness of each edge or on the closeness between each pair of nodes. On the opposite, macro-scale analysis focus on global properties of the network such as its density or its modularity. The reader may refer to \cite{New03} or \cite{SCB19} for a general discussion. In this paper, we are mostly interested in the meso-scale description of the network that is provided by the frequency of motifs \citep{MSI02}. \\
A motif is defined as a given subgraph depicting the interactions between a small number of nodes; the count of a motif consists in the number of occurrences of this subgraph in the observed network. Figures \ref{fig:motifs1-5} and \ref{fig:motifs6} display the set of all bipartite motifs involving up to 6 top or bottom nodes. Counting the occurrences of a motif is a computationally challenging task \citep[see][for simple--~i.e. non-bipartite~-- networks]{MSI02,PDK08}; efficient tools have been recently proposed by \cite{SSS19,SCB19} for bipartite networks.

\bigskip
Whatever the description scale, the analysis must account for a series of characteristics of the network at hand (such as its dimension or its density) to make the results comparable. A convenient way to account for such peculiarities is to define a null model capable to fit the network characteristics. We consider here a bipartite and exchangeable version of the expected degree distribution model proposed by \cite{ChL02} for simple binary graphs. The bipartite expected degree distribution (\BEDD) model simply states that each (top or bottom) node is associated with an expected degree and that a pair of nodes is connected with a probability that is proportional to the product of their respective expected degrees. \\
The \BEDD model can obviously accommodate to the network dimension (number of top and bottom nodes), for its density but also for some existing imbalances between the degrees of the nodes. Such imbalances play an important role in many fields: in ecology they are related to the opposition between generalist insects (capable of pollinate a large number of plant species) and specialist insects (interacting with a limited number of plant species) \citep{SCB19}. \\
In addition to its interpretation, this model is attractive because we can calculate the expected frequency of motifs under \BEDD such as their variance.

\bigskip
The distribution of motif counts in simple graphs has been widely studied, especially for simple motifs like triangles \citep[see e.g.][]{NoW88,Sta01,PDK08}.
In this paper, we prove the asymptotic normality of the count of any given motif under the \BEDD model, under sparsity conditions. One important feature of the \BEDD model is that the mean and the variance of the count have close form expressions. The strategy to derive these moments is related to the one introduced by \cite{PDK08} for simple networks. \\
This property has a major practical impact as the expectation and the variance of a motif count could not be evaluated via resampling, because of the computational cost of motif counting event for networks with intermediate size. The knowledge of the asymptotic distribution of the motif counts opens a series of possible applications, including goodness-of-fit tests for the \BEDD model and a series of tests for network comparison in the \BEDD framework.

\bigskip
The paper is organized as follows. Section \ref{sec:model} is devoted to the definition and properties of motifs in the \BEDD model and Section \ref{sec:tests} to tests for bipartite networks.  More specifically, we establish the asymptotic normality of motif frequencies in Section \ref{sec:asymp} and propose a goodness-of-fit test for the \BEDD model and comparison tests for two bipartite networks in Section \ref{sec:test} and Section \ref{sec:comparison}, respectively.
The accuracy of the normal approximation for finite graphs and the power of the proposed tests are assessed via a simulation study in Section \ref{sec:simul}. Finally, proofs are given in Section  \ref{sec:proof}.


\section{Motifs in the bipartite expected degree model} \label{sec:model}

We consider a bipartite graph $\Gcal = (\Vcal, \Ecal)$ with $N$ nodes. The set of nodes is $\Vcal = (\Vcal^t, \Vcal^b)$, where $\Vcal^t=\llbracket{1,m} \rrbracket$ (resp. $\Vcal^b=\llbracket{1,n} \rrbracket$) stands for the set of top (resp. bottom) nodes, and the set of edges is $\Ecal \subset \Vcal^t \times \Vcal^b$, meaning than an edge can only connect a top node with a bottom node. The total number of nodes is therefore $N = n+m$. We denote by $G$ the corresponding $m \times n$ incidence matrix where the entry $G_{ij}$ of $G$ is 1 if $(i,j)\in \Ecal$, and 0
otherwise. 

\subsection{Bipartite expected degree model}
 
The bipartite expected degree (\BEDD) model is defined as follows: 
\begin{align} \label{eq:BEDD}
  \{U_i\}_{1 \leq i \leq m} \text{ iid}, \qquad U_1 & \sim \Ucal_{[0, 1]}, \nonumber \\
  \{V_j\}_{1 \leq j \leq n} \text{ iid}, \qquad V_1 & \sim \Ucal_{[0, 1]}, \\
  \{G_{ij}\}_{1 \leq i \leq m, 1 \leq j \leq n} \text{ indep. } | \{U_i\}_{1 \leq i \leq m}, \{V_j\}_{1 \leq j \leq n}, \qquad G_{ij} | U_i, V_j & \sim \Bcal\left(\rho g(U_i) h(V_j)\right), \nonumber 
\end{align}
where $\gt, \gb: [0, 1] \mapsto \Rbb^+$, such that $\int \gt(u) \d u = \int \gb(v) \d v = 1$ and $0 \leq \rho \leq 1$. 

The parameter $\rho$ controls the density of the graph ($\Esp G_{ij} = \rho$) whereas the function $g$ (resp. $h$) encodes the heterogeneity of the expected degrees of the top (resp. bottom) nodes. More specifically, denoting $K_i = \sum_{1\leq j \leq n} G_{ij}$ the degree of the top node $i$, we have that $\Esp(K_i \mid U_i) = n \rho g(U_i)$. The symmetric property holds for bottom nodes.

\begin{remark}
\cite{LoS06} and \cite{DiJ08} introduced a generic model for exchangeable random graphs called the $W$-graph, which is based on a {\sl graphon} function $\Phi: [0, 1]^2 \mapsto [0, 1]$. The \BEDD model is a natural extension of the $W$-graph for bipartite graphs with a product-form graphon function $\Phi(u, v) = \rho g(u) h(v)$. The \BEDD model is obviously exchangeable is the sens that the distribution of the incidence matrix $G$ is preserved under permutation of the top nodes and/or the bottom nodes.
\end{remark}

\begin{remark}
The \BEDD model can also be seen has an exchangeable bipartite version of the expected degree sequence model studied in \cite{ChL02} and of the configuration model from \cite{New03}. 
Under these two models, the degree of each node is fixed which makes them non exchangeable.
\end{remark}

\subsection{Bipartite motifs in the \BEDD model}

\paragraph{Bipartite motifs.} 
We are interested in the distribution of the count of motifs (or subgraphs) in bipartite graphs arising from the \BEDD model. 
A bipartite motif $s$ is defined by its number of top nodes $\pt_s$, its number of bottom nodes $\pb_s$ and a $\pt_s \times \pb_s$ incidence matrix $A^s$. 
Figures \ref{fig:motifs1-5} and \ref{fig:motifs6} display the 44 bipartite motifs involving between two and six nodes, from which we see that
$$
A^2 = \left(\begin{array}{cc}
            1 & 1 
            \end{array} \right), \qquad
A^5 = \left(\begin{array}{cc}
            1 & 1 \\ 0 & 1
            \end{array} \right), \qquad
A^{15} = \left(\begin{array}{ccc}
            1 & 1 & 0 \\ 1 & 1 & 1
            \end{array} \right).
$$
An important characteristic of a graph motif $s$ is its number of automorphisms $r_s$ \citep{Sta01}, that is the number of non-redundant permutations of its incidence matrix (see, e.g. section 2.4 in \cite{PDK08}):
\begin{eqnarray} \label{eq:rs}
r_s = \left|\left\{A^s_{\sigma^t,         \sigma^b}         =
  \left(A^s_{\sigma^t(u),\sigma^b(v)}\right)_{1 \leq  u \leq  \pt_s, 1
    \leq v \leq \pb_s}: 
  \sigma^t \in \sigma\left(\llbracket 1,\pt_s\rrbracket\right), 
  \sigma^b \in \sigma\left(\llbracket 1,\pb_s\rrbracket\right) \right\} \right|.
\end{eqnarray}
Note that, because pairs of permutations $(\sigma^t, \sigma^b)$ yielding the same matrix $A^s_{\sigma^t, \sigma^b}$ are not counted twice, we obviously have that $r_s \leq (\pt_s!) \times (\pb_s!)$. In many cases, $r_s$ turns out to be much smaller: in particular, $r_s = 1$ for star-motifs, which will be defined later.
We further denote by $\dt^s_u$ the degree of the top node $u$ ($1 \leq u \leq \pt_s$) within motif $s$, that is $\dt^s_u = \sum_{1 \leq v \leq \pb_s} A^s_{u, v}$. The degree of the bottom node $v$ within $s$ is defined similarly as $\db^s_v = \sum_{1 \leq u \leq \pt_s} A^s_{u, v}$. \\

\paragraph{Motif occurrence.} 
Counting the  occurrences of motif  $s$ in $\Gcal$ simply  consists in considering all possible of $\pt_s$ (resp. $\pb_s$) top (resp. bottom) nodes  among  the  $m$  (resp.   $n$)  and  check  for  each  possible automorphism of $s$ if an occurrence is observed. More   formally, let us define the  set $\Pcal_s$ of possible positions for motif   $s$  as   the  Cartesian   product   of  the   set  of   the   $\binom{m}{\pt_s} \binom{n}{\pb_s}$ possible  locations with the set of the $r_s$ (top, bottom) permutations giving  rise to each of the automorphisms of $s$. So, a {\sl position} results from the combination of a {\sl location} with a {\sl permutation}. Because the graph is bipartite, any position $\alpha$ from $\Pcal_s$ decomposes as $\alpha = (\alpha^t, \alpha
^b)$ where $\alpha^t$ stands for an ordered list of top nodes and $\alpha^b$ for an ordered list of bottom nodes. The number of positions for motif $s$ in $\Gcal$ is precisely
\begin{eqnarray} \label{eq:cs}
c_s := |\Pcal_s| = r_s \binom{m}{\pt_s} \binom{n}{\pb_s}.
\end{eqnarray}
Now, for a given position $\alpha = (\alpha^t,\alpha^b) \in \Pcal_{s}$, we define $Y_s(\alpha)$ as the indicator for motif $s$ to occur in position $\alpha$:
\begin{eqnarray} \label{eq:Yalpha}
Y_s(\alpha)=\prod_{i \in \alpha^t, j \in \alpha^b} G_{ij}^{A^s_{ij}}.
\end{eqnarray}
\begin{remark}
Note that the occurrence defined by Equation \eqref{eq:Yalpha} corresponds to an \emph{induced} occurrence, which means that we consider that a motif $s$ is observed at position $\alpha$ as soon as all the present edges that are specified by its incidence matrix $A^s$ are observed, even if additional edges are also observed. In other words, we do not check for the {\sl absent} edges specified by $A^s$. 
\end{remark}

\begin{remark}
As opposed to an induced occurrence, an \emph{exact} occurrence is observed when both the presence and the absence of edges are satisfied. The indicator variable corresponding to an exact occurrence writes $\prod_{i \in \alpha^t, j \in \alpha^b} G_{ij}^{A^s_{ij}} (1 - G_{ij})^{1-A^s_{ij}} $. Counting induced and exact occurrences in a graph is actually equivalent, as these counts are related in a deterministic manner. For example, each exact occurrence of motif 6 corresponds to two induced occurrences of motif 5.
\end{remark}

\paragraph{Motif probability.} 
The \BEDD model is an exchangeable bipartite graph model in the sense that, for any pair of permutations 
$(\sigma^t \in \sigma\left(\llbracket 1,m\rrbracket\right), \sigma^b \in \sigma\left(\llbracket 1,n\rrbracket\right))$, 
we have that 
$\Pr\{G = \{g_{ij}\}_{1 \leq i\leq m, 1\leq j\leq n}\} = \Pr\{G = \{g_{\sigma^t(i) \sigma^b(j)}\}_{1 \leq \sigma^t(i)\leq m, 1\leq \sigma^b(j)\leq n}\}$ 
\citep[see e.g.][for simple graphs]{LoS06,DiJ08}.
For any exchangeable graph model, we may define $\phi_s$ as the probability for motif $s$ to occur in position $\alpha = (\alpha^t, \alpha^b)$:
$$
\phi_s := \Pr\left(Y_s(\alpha) = 1\right).
$$
Importantly, because the model is exchangeable, this probability does not depend on $\alpha$.

\paragraph{Star motifs.}
We define a star as a bipartite motif $s$ for which either $\pb_s = 1$ or $\pt_s = 1$ (or both). More specifically, we name top stars (resp bottom stars) motifs for which $\pt_s = 1$ (resp. $\pb_s= 1$). The top stars in Figures \ref{fig:motifs1-5} and \ref{fig:motifs6} are motifs 1, 2, 7, 17 and 44, and the bottom stars are motifs 1, 3, 4, 8 and 18. Observe that $r_s = 1$ for all star motifs, that $\dt^s_v = 1$ for all $v$ in all top star motifs, and that $\db^s_u = 1$ for all $u$ in all bottom star motifs. \\
Because they will play a central role in the sequel, we adopt a specific notation for the probability of star motifs, denoting $\st_d$ the occurrence probability of the top star with degree $d$ and $\sb_d$ for the occurrence probability of the bottom star with degree $d$. As a consequence, we have that
\begin{align} \label{eq:lambdagamma}
  \st_1 & = \phi_1, & \st_2 & = \phi_2, & \st_3 & = \phi_7 & \st_4 & = \phi_{17}, & \st_5 & = \phi_{44}, \\
  \sb_1 & = \phi_1, & \sb_2 & = \phi_3, & \sb_3 & = \phi_4 & \sb_4 & = \phi_8, & \sb_5 & = \phi_{18}. \nonumber
\end{align}

\subsection{Moments of motif counts}

\paragraph{Expected count.} 
Let us now denote by $N_s$ the count, that is the number of occurrences of a motif $s$ in a graph $\Gcal$. We simply have that
$$
N_s = \sum_{\alpha \in \Pcal_s} Y_s(\alpha)
$$
As a consequence, the expected count of $s$ in $\Gcal$ is $\Esp(N_s) = c_s \phi_s$.
We also define the normalized frequency of motif $s$ as 
$$
F_s = N_s / c_s,
$$
which is an unbiased estimate of $\phi_s$.

\paragraph{Illustration.}
As an illustration, we consider two of the networks studied by \cite{SSS19}, which include both plant-pollinator and seed dispersal networks extracted from the Web of Life database ({\tt www.web-of-life.es}). More specifically, we consider the two largest networks of each type, which were first published by \cite{Rob29} and \cite{Sil02}, respectively. The plant-pollinator network involves 546 plant species and 1044 insects and the seed dispersal network 207 plant species and 110 seed dispersers (birds or insects). Our purpose is not to provide a thorough ecological analysis of these networks, but to examplify the proposed methodology. Table~\ref{tab:stars} gives the counts and the frequency of the star motifs with up to four branch. For the sake of clarity, we will limit ourselves to motifs up to five nodes\SR{ in the illustrations}{}. Observe that both the counts $N_s$ and the number of possible positions $c_s$ range over huge order of magnitudes.
\begin{table}[ht]
  \begin{center}
    \begin{tabular}{cl|lll|lll}
      \multicolumn{8}{c}{plant-pollinator: $m=546, n=1044$ \citep{Rob29}} \\ 
      \hline
      & \multicolumn{1}{c}{edge}  & \multicolumn{3}{c}{top stars} & \multicolumn{3}{c}{bottom stars} \\
      $s$ & \quad 1 & \quad 2 & \quad 7 & \quad 17 & \quad 3 & \quad 4 & \quad 8 \\
      \hline
      $c_s$ & 
      4.76\;$10^{5}$ & 2.48\;$10^{8}$ & 8.62\;$10^{10}$ & 2.24\;$10^{13}$ & 1.08\;$10^{8}$ & 1.64\;$10^{10}$ & 1.86\;$10^{12}$ \\
      $N_s$ & 
      1.53\;$10^{4}$ & 2.61\;$10^{5}$ & 3.04\;$10^{6}$ & 2.72\;$10^{7}$ & 3.07\;$10^{5}$ & 6.82\;$10^{6}$ & 1.48\;$10^{8}$ \\
      $F_s$ & 
      3.20\;$10^{-2}$ & 1.05\;$10^{-3}$ & 3.52\;$10^{-5}$ & 1.21\;$10^{-6}$ & 2.84\;$10^{-3}$ & 4.16\;$10^{-4}$ & 7.99\;$10^{-5}$ \\
      \\ 
      \multicolumn{8}{c}{seed dispersal: $m=207, n=110$ \citep{Sil02}} \\ 
      \hline
      & \multicolumn{1}{c}{edge}  & \multicolumn{3}{c}{top stars} & \multicolumn{3}{c}{bottom stars} \\
      $s$ & \quad 1 & \quad 2 & \quad 7 & \quad 17 & \quad 3 & \quad 4 & \quad 8 \\
      \hline
      $c_s$ & 
      2.28\;$10^{4}$ & 1.24\;$10^{6}$ & 4.47\;$10^{7}$ & 1.20\;$10^{9}$ & 2.35\;$10^{6}$ & 1.60\;$10^{8}$ & 8.17\;$10^{9}$ \\
      $N_s$ & 
      1.12\;$10^{3}$ & 6.50\;$10^{3}$ & 4.07\;$10^{4}$ & 2.32\;$10^{5}$ & 1.24\;$10^{4}$ & 1.31\;$10^{5}$ & 1.23\;$10^{6}$ \\
      $F_s$ & 
      4.92\;$10^{-2}$ & 5.23\;$10^{-3}$ & 9.11\;$10^{-4}$ & 1.94\;$10^{-4}$ & 5.28\;$10^{-3}$ & 8.16\;$10^{-4}$ & 1.50\;$10^{-4}$ \\
    \end{tabular}
    \caption{Coefficients $c_s$, counts $N_s$ and frequency $F_s$ of all star motifs.
    Top: plant-pollinator network, bottom: seed dispersal network. The motif number $s$ refers to Figure~\ref{fig:motifs1-5}.} \label{tab:stars}
  \end{center}
\end{table}

\paragraph{Main property of motif probabilities under \BEDD.} 
The tests we propose rely on the comparison between the observed count (or normalized frequency) of a motif, with its theoretical counterpart under a \BEDD model. More specifically, the motif probabilities have a close form expression under the \BEDD model. 
\begin{proposition} \label{prop:phiBEDD}
  Under the \BEDD model \eqref{eq:BEDD}, we have that
  \begin{equation} \label{eq:phi}
    \phibar_s
    = {\prod_{u = 1}^{\pt_s}  \st_{\dt^s_u} \prod_{v = 1}^{\pb_s}\sb_{\db^s_v}} \left/ {(\phi_1)^{d_+^s}} \right. .
  \end{equation}
  where $\dt_+^s := \sum_u \dt_u^s =  \sum_v \db_v^s$ stands for the total number of edges in $s$. 
\end{proposition}

\proofbegin
This follows from the fact that, under \BEDD, the edges are independent conditionally on the latent coordinates $U_i$ and $V_j$ defined in \eqref{eq:BEDD}, which are all independent with respect to one other. Consider an arbitrary position $\alpha = (\alpha^t, \alpha^b)$; for the sake of clarity, we identify the elements of $\alpha^t$ with $\llbracket 1, \pt_s \rrbracket$ and the elements of $\alpha^b$ with $\llbracket 1, \pb_s \rrbracket$. We have 
\begin{align*}
\phibar_s
  & = \Esp_{(U_i)_{1 \leq i \leq \pt_s}, (V_j)_{1 \leq j \leq \pb_s}}
  \left( \left.
  \Pr\left\{\prod_{1\leq i \leq \pt_s, 1 \leq v \leq \pb_s} G_{ij}^{A^s_{ij}} = 1
   \right| (U_i)_{1 \leq i \leq \pt_s}, (V_j)_{1 \leq j \leq \pb_s} \right\} \right) \\
  & = \Esp_{(U_i)_{1 \leq i \leq \pt_s}, (V_j)_{1 \leq j \leq \pb_s}}
  \left( 
    \prod_{1 \leq i \leq \pt_s, 1 \leq j \leq \pb_s: A_{ij}^s = 1} \rho \gt(U_i) \gb(V_j) \right) \\
  & = \Esp_{(U_i)_{1 \leq i \leq \pt_s}, (V_j)_{1 \leq j \leq \pb_s}}
  \left( 
    \rho^{\dt^s_+} 
    \prod_{1 \leq i \leq \pt_s} \gt(U_i)^{\dt^s_i}
    \prod_{1 \leq j \leq \pb_s} \gb(V_j)^{\db^s_j} \right) \\
  & = 
  \rho^{\dt^s_+} 
  \prod_{1 \leq i \leq \pt_s} \left(\int \gt(u)^{\dt^s_i} \d u\right)
  \prod_{1 \leq j \leq \pb_s} \left(\int \gb(v)^{\db^s_j} \d v\right).
\end{align*}
The result then results from the fact that
\begin{equation} \label{eq:int_g_h}
\st_d = \rho^d \int g(u)^d \d u, \qquad
\sb_d = \rho^d \int h(v)^d \d v, \qquad
\rho = \phi_1.
\end{equation}
\proofend

\begin{figure}[ht]
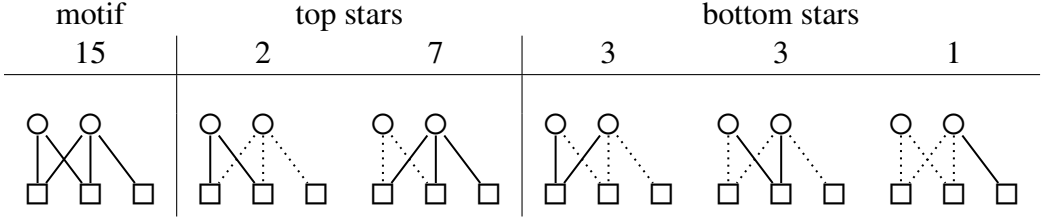

  \begin{center}
    $\begin{array}{c|cc|ccc}
      \multicolumn{1}{c}{\text{motif}} & \multicolumn{2}{c}{\text{top stars}} & \multicolumn{3}{c}{\text{bottom stars}} \\
      15 & 2 & 7 & 3 & 3 & 1 \\
      \hline 
      & & & & & \\
      ~ \input{./Figs/15-bipartite-Xlambda} ~ &
      ~ \input{./Figs/15-bipartite-decomposition-2} ~ &
      ~ \input{./Figs/15-bipartite-decomposition-7} ~ &
      ~ \input{./Figs/15-bipartite-decomposition-3a} ~ &
      ~ \input{./Figs/15-bipartite-decomposition-3b} ~ &
      ~ \input{./Figs/15-bipartite-decomposition-1} ~ 
    \end{array}$
  \end{center}
  \caption{Decomposition of motif 15 as an overlap of 2 top stars (motifs 2 and 7) and 3 bottom stars (motifs 3, 3 and 1). Because each edge is accounted for twice, we get $\phibar_{15} = \phi_2 \phi_7 \phi_3 \phi_3 \phi_1 / \phi_1^5 = \phi_2 \phi_7 \phi_3^2 / \phi_1^4$. \label{fig:decompMotif}}
\end{figure}

An important consequence of Proposition \ref{prop:phiBEDD} is that, under \BEDD, the motif probability of any motif can be expressed in terms of probabilities of star motifs. Figure~\ref{fig:decompMotif} provides an intuition of this: a motif can be decomposed in terms of top and bottom stars arising from each of its nodes. 

In the sequel, to distinguish the motif probability $\phi_s$ under an arbitrary exchangeable model from the probability under the \BEDD model, we will denote by $\phibar_s$ the probability of motif $s$ under \BEDD. Figure~\ref{fig:motifs1-5} provides the list of all $\phibar_s$ expressions.

\paragraph{Probability estimate under \BEDD.} 
Proposition \ref{prop:phiBEDD} suggests a natural plug-in estimator for the \BEDD motif probability $\phibar_s$: 
\begin{equation} \label{eq:F.bar}
 \Fbar_s 
= \frac{\prod_{u = 1}^{\pt_s}  \St_{\dt^s_u} \prod_{v = 1}^{\pb_s} \Sb_{\db^s_v}}{F_1^{d_+^s}},
\end{equation}
where $\St_d$ (resp $\Sb_d$) denotes the normalized frequency of the top (resp. bottom) star motif with degree $d$. Obviously, $\St_d$ (resp $\Sb_d$) is an unbiased estimated of $\st_d$ (resp. $\sb_d$).

\begin{figure}[ht]
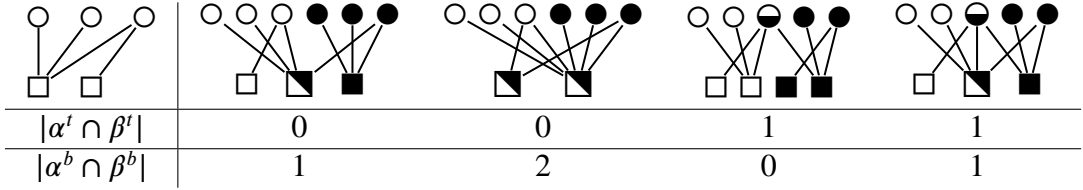

  \begin{center}
    $\begin{array}{c|ccccc}
      ~ \input{./Figs/9-bipartite-star3upN} ~ &
      ~ \input{./Figs/15-bipartite-supermotif-9-0-1} ~ &
      ~ \input{./Figs/15-bipartite-supermotif-9-0-2} ~ &
      ~ \input{./Figs/15-bipartite-supermotif-9-1-0} ~ &
      ~ \input{./Figs/15-bipartite-supermotif-9-1-1} ~ \\
      \hline
      |\alpha^t \cap \beta^t| & 0 & 0 & 1 & 1 \\
      \hline
      |\alpha^b \cap \beta^b| & 1 & 2 & 0 & 1 \\
    \end{array}$
  \end{center}
  \caption{Some super-motifs from $\Scal_2(s)$ for motif $s = 9$ (top left) with $\pt_s = 3$ top nodes and $\pb_s = 2$ bottom nodes. $|\alpha^t \cap \beta^t|$ (resp. $|\alpha^b \cap \beta^b|$): number of top (resp. bottom) nodes shared by the overlapping positions $\alpha$ and $\beta$. 
  Black: nodes from $\alpha$, white: nodes from $\beta$, black/white: nodes from $\alpha \cap \beta$.
  There are actually $|\Scal_2(9)| = 396$ such super-motifs of motif 9.  \label{fig:superMotif}}
\end{figure}

\paragraph{Variance of the count.} 
We now consider the variance of the count, that is 
\begin{align} \label{eq:varNs}
  \Var(N_s) & = \Esp(N^2_s) - \Esp(N_s)^2, \nonumber \\
  \text{where} \qquad 
  N^2_s 
  & = \sum_{\alpha, \beta \in \Pcal_s} Y_s(\alpha) Y_s(\beta) \\
  & 
  =   \sum_{\alpha \in \Pcal_s} Y_s(\alpha) + \sum_{\alpha, \beta \in \Pcal_s: |\alpha \cap \beta| = 0} Y_s(\alpha) Y_s(\beta)  + \sum_{\alpha, \beta \in \Pcal_s: \alpha \neq \beta, |\alpha \cap \beta| > 0} Y_s(\alpha) Y_s(\beta). \nonumber
\end{align}
When positions $\alpha$ and $\beta$ are equal, the product $Y_s(\alpha) Y_s(\beta)$ is simply given by $Y_s(\alpha)$, the indicator of the presence of $s$ at position $\alpha$.  
Then, when positions $\alpha$ and $\beta$ do not overlap ($|\alpha \cap \beta| = 0$), the product $Y_s(\alpha) Y_s(\beta)$ simply indicates that two occurrences of motif $s$ occur in position $\alpha$ and $\beta$, which are independent under the \BEDD model. When positions $\alpha$ and $\beta$ are different and do overlap ($|\alpha \cap \beta| > 0$), the product $Y_s(\alpha) Y_s(\beta)$ becomes the indicator of a super-motif, that is a motif made of two overlapping automorphisms of $s$. We denote by $\Scal_2(s)$ the set of super-motifs generated by the overlaps of two occurrences of the motif $s$; Figure~\ref{fig:superMotif} provides some examples of super-motifs. 

An expression similar to \eqref{eq:varNs} can be derived for the covariance between two counts:
\begin{align} \label{eq:covNsNt}
  \Cov(N_s, N_t) & = \Esp(N_sN_t) - \Esp(N_s)\Esp(N_t), \nonumber \\
  \text{where} \qquad 
  N_sN_t 
  & = \sum_{\alpha \in \Pcal_s, \beta \in \Pcal_t} Y_s(\alpha) Y_t(\beta) \\
  & 
  = \sum_{\alpha \in \Pcal_s, \beta \in \Pcal_t: |\alpha \cap \beta| = 0} Y_s(\alpha) Y_t(\beta)  + \sum_{\alpha \in \Pcal_s, \beta \in \Pcal_t: \alpha \neq \beta, |\alpha \cap \beta| > 0} Y_s(\alpha) Y_t(\beta). \nonumber
\end{align}
Again, the last term corresponds to occurrences of super-motifs resulting from an overlap between an occurrence of motif $s$ and an occurrence of motif $t$. We denote by $S_2(s, t)$ the set of these super-motifs. We use the strategy described in \cite{PDK08} to determine the sets of super-motifs $\Scal_2(s)$ and $\Scal(s, s')$. Observe that these sets do not depend on the observed networks, so, to alleviate the computational burden, they can be determined and stored once for all.

Eq.~\eqref{eq:varNs} shows that $\Esp(N^2_s)$ only depends on $\Esp(Y_s(\alpha)Y_s(\beta))$, which is $\phi_s^2$ when positions $\alpha$ and $\beta$ do not overlap and the probability of the corresponding super-motif when they overlap. As a consequence, we have that
\begin{align} \label{eq:espNs2}
\Esp(N^2_s) 
= \kappa_{m, n, s} \phi_s + \kappa'_{m, n, s} \phi_s^2 
+ \sum_{S \in \Scal_2(s)} \kappa''_{m, n, s, S} \phi_S,
\end{align}
where the $\kappa_{m, n, s}$, $\kappa'_{m, n, s}$, $\kappa''_{m, n, s, S}$ are constants, which depend on the dimensions of the graph, on the motif $s$ and on the super-motif $S$. The order of magnitude of $\kappa_{m, n, s, S}$ for large $m$ and $n$ will be studied in Section \ref{Sec:tech:lem}.

Because super-motifs are actually motifs, their respective occurrence probability $\phibar_S$ under \BEDD are given by Proposition \ref{prop:phiBEDD} as well, so the expectation and the variance of $N_s$ under \BEDD can be expressed as functions of the $\phibar_s$ and $\{\phibar_S\}_{S \in \Scal_2(s)}$. An estimate $\Fbar_S$ of each $\phibar_S$ can be obtained using Eq.~\eqref{eq:F.bar} in the same way.

\begin{remark}
 The estimate defined in \eqref{eq:F.bar} is only based on  empirical quantities (the counts of stars motifs) and does not depend on any parameter estimation. Especially, the functions $g$ and $h$ do not need to be estimated as the frequency of star motifs provides all necessary information about the degree distributions. As a consequence, we may define plug-in estimates of the occurrence probability, the expected count and the variance of the count of any motif under \BEDD.
\end{remark}

\paragraph{Illustration.}
Table~\ref{tab:freqs} compares the empirical frequencies $F_s$ of a selection of motifs with their respective estimated probability $\Fbar_s$. The probability estimates are computed according to Equation~\ref{eq:F.bar}, using the star motifs frequencies $\St_{\dt}$ and $\Sb_{\db}$ given in Table~\ref{tab:stars}. Observe that the difference between the observed frequency $F_s$ and their estimated expectation under the \BEDD model $\Fbar_s$ are of the same order of magnitude, if not smaller, than their estimated standard deviations. 

\begin{table}[ht]
  \begin{center}
    \begin{tabular}{clllll}
      \multicolumn{6}{c}{plant-pollinator} \\ 
      \hline
      $s$ & \quad 5 & \quad 6 & \quad 10 & \quad 15 & \quad 16 \\
      \hline
      $F_s$ & 
      9.21\;$10^{-5}$ & 1.00\;$10^{-5}$ & 8.12\;$10^{-6}$ & 3.32\;$10^{-7}$ & 4.47\;$10^{-8}$ \\
      $\Fbar_s$ & 
      9.29\;$10^{-5}$ & 8.41\;$10^{-6}$ & 8.23\;$10^{-6}$ & 2.82\;$10^{-7}$ & 2.62\;$10^{-8}$ \\
      $\sqrt{\widehat{\Var}(F_s)}$ & 
      1.26\;$10^{-5}$ & 1.61\;$10^{-6}$ & 1.58\;$10^{-6}$ & 6.54\;$10^{-8}$ & 7.60\;$10^{-9}$ \\
      \\ 
      \multicolumn{6}{c}{seed dispersal} \\ 
      \hline
      $s$ & \quad 5 & \quad 6 & \quad 10 & \quad 15 & \quad 16 \\
      \hline
      $F_s$ & 
      5.13\;$10^{-4}$ & 1.15\;$10^{-4}$ & 5.07\;$10^{-5}$ & 1.79\;$10^{-5}$ & 5.96\;$10^{-6}$ \\
      $\Fbar_s$ & 
      5.61\;$10^{-4}$ & 1.30\;$10^{-4}$ & 6.02\;$10^{-5}$ & 2.26\;$10^{-5}$ & 8.59\;$10^{-6}$ \\
      $\sqrt{\widehat{\Var}(F_s)}$ & 
      2.25\;$10^{-4}$ & 7.24\;$10^{-5}$ & 3.25\;$10^{-5}$ & 1.59\;$10^{-5}$ & 7.38\;$10^{-6}$ \\
    \end{tabular}
    \caption{Empirical frequency $F_s$, estimated probability $\Fbar_s$ and estimated standard-deviation of the frequency according to the \BEDD model for a selection of motifs. All estimates are derived from the star motifs frequencies given in Table~\ref{tab:stars}.} \label{tab:freqs}
  \end{center}
\end{table}


\section{Tests for bipartite networks} \label{sec:tests}
\paragraph{Asymptotic framework.}
We consider a sequence of \BEDD random graphs defined as follows.

$\{\Gcal_N\}_{N \geq 2}$ is a sequence of independent graphs, where $\Gcal_N$ is a \BEDD random graph with $m = \lfloor \lambda N \rfloor$ top nodes with $\lambda\in(0,1)$, $n = N -m$ bottom nodes and parameters $\rho_N$, $\gb$ and $\gt$, where the sequence $\{\rho_N\}_{N \geq 2}$ satisfies $\rho_N = \Theta(m^{-a} n^{-b})$ with $a, b > 0$.
All quantities computed on $\Gcal_N$ should be indexed by $N$ as well but for the sake of clarity, we will drop that index in the rest of the paper.

\subsection{Asymptotic normality of motif frequencies}  \label{sec:asymp}

This section is devoted to the asymptotic normality of motif frequencies under the \BEDD model. More precisely,
 our first main result states the asymptotic normality of the following statistic $W_s$ relying on $F_s$ the empirical frequency of a given motif $s$ in $\Gcal$:
\begin{equation}\label{eq::Ws}
    W_s = \frac{F_s - \Fbar_s}{\sqrt{\widehat{\Var}(F_s)}},
\end{equation}
where $\Fbar_s$ denotes the estimator of $\phibar_s$ defined in \eqref{eq:F.bar} and  $\widehat{\Var}(F_s)$ the one of $\Var(F_s)$ obtained by the plug-in of $\Fbar_S$ ($S$ being any super-motif generated by two occurrences of $s$) in the expressions of $\Var(N_s)$ given in \eqref{eq:varNs}-\eqref{eq:espNs2}.

\begin{theorem}\label{Thm:AsympNorm}
{If $a+b<2/d_+^{s}$, then for all non-star motif $s$ and under the \BEDD model,} the statistic $W_s$ is asymptotically normal as $m\sim n\to\infty$:
$$W_s\overset{D}{\longrightarrow} \Ncal(0, 1).$$
\end{theorem}

The proof is based on three results given hereafter in Proposition  \ref{Prop:Ls},  Lemma \ref{Lemma:Cs} and Lemma \ref{Lemma:Var}.

\bigskip
\noindent{\sl Sketch of proof.} 
Let first consider the following decomposition of the numerator of $W_s$ : 
\begin{align*} 
 F_s - \Fbar_s  & := L_s + C_s &
 \text{where  } L_s & = F_s - \phibar_s 
 \text{ and } C_s  = \phibar_s - \Fbar_s.
\end{align*}
Under the null \BEDD model, we show that, $(i)$ $L_s/\sqrt{\Var(F_s)}$ is asymptotically normal in Proposition \ref{Prop:Ls}, it is the leader term, $(ii)$ $C_s/\sqrt{\Var(F_s)}$ is negligible in Lemma \ref{Lemma:Cs}, it is the reminder term. Then, we conclude using Slutsky Theorem Lemma \ref{Lemma:Var} which states that $\widehat{\Var}(F_s)/\Var(F_s)\rightarrow 1$ in probability.
\proofend 

\begin{remark} 
Like $\Fbar_s$, $W_s$ is only based on empirical quantities, that is $i$) the empirical frequency of motif $s$ and $ii$) the empirical frequencies of the stars motifs forming $s$. The expected frequencies of the supermotifs of $s$ involved in $\widehat{\Var}(F_s)$ also depend only on empirical star frequencies.
\end{remark}

\begin{remark} 
\cite{GaL17b} proved a similar result as Theorem \ref{Thm:AsympNorm} in the EDD model, for a test statistic which is a linear combination of edges, vees and triangles empirical frequencies in the case of simple graphs, and under a specific condition on the graph density. Though their result is not comparable to ours since triangles can not occur in bipartite graphs and we do not account for stars motifs.
Although they seem similar, a fair comparison between Theorem \ref{Thm:AsympNorm} and the result from \cite{GaL17b} is not easy ($i$) because the model is not the same (we consider bipartite graphs whereas they consider simple graphs) and ($ii$) because they only consider vees (which are star-motifs) and triangles (which do not occur in bipartite graphs).
\end{remark}

In the following proposition, the asymptotic normality of the statistic ruling the law of $W_s$ is stated under the null. This statistic involves the empirical frequency of {a given non star} motif $s$ and its theoretical expectation and variance. The proof of its asymptotic normality mostly relies on tools of martingale theory. We show that we can exhibit conditional martingale difference sequences relative to a specific filtration. This filtration is generated by the sequence of graphs $\Gcal_N$ (see a proper definition of the filtration in Section \ref{Sec:Notations}). So, we could apply the central limit theorem of \cite{HaHy14}.

\begin{proposition}\label{Prop:Ls}
{If $a+b<2/d_+^{s}$, then for all star motif $s$ and under the \BEDD model,} we have, as $m\sim n\to\infty$,
 $$
 \frac{F_s - \phibar_s}{\sqrt{\Var(F_s)}}\overset{D}{\longrightarrow} \Ncal(0, 1).
 $$
\end{proposition}

The complete proof is given in Section \ref{sec:prop}, it relies especially on Lemma \ref{Lemma:Rs} and Lemma \ref{Lemma:Ms}.

\bigskip 
\noindent{\sl Sketch of proof.} 
We first consider the decomposition $L_s=F_s -\phibar_s=M_s+R_s$ with $M_s$ being the difference between $F_s$ and its expectation conditionally to the considered filtration and $U,V$, and $R_s$ the difference between the latter conditional expectation and $\phibar_s$; the proper definitions are given in Section \ref{subsec:Ls}. Lemma \ref{Lemma:Rs} shows that, under the null \BEDD model, the reminder term $R_s/\sqrt{\Var(F_s)}|U,V\to 0$ a.s. as $m\sim n\to\infty$. Lemma \ref{Lemma:Ms} shows that, under {the \BEDD model},  the leader term $M_s/\sqrt{\Var(F_s)}| U,V$ is asymptotically normal with variance $\Var(N_s | U,V)/\Var(N_s)$. Slutsky theorem implies the asymptotic normality of $L_s/\sqrt{\Var(F_s)}$ conditional on $(U, V)$. 
Then, Lemma \ref{Lemma:Var:deconditionning} shows that $\Var(N_s | U,V)/\Var(N_s)$ tends to 1 in probability for all $(U, V)$, which allows deconditionning.
\proofend 

The two following lemmas combined with Proposition \ref{Prop:Ls} permit to conclude to Theorem \ref{Thm:AsympNorm}. Their proofs are given in sections \ref{sec:lemma:cs} and \ref{sec:lemma:var}  respectively. 
\begin{lemma}\label{Lemma:Cs}
{If $a+b<2/d_+^{s}$, then for all non-star motif $s$ and under {the \BEDD model},}  we have, as $m\sim n\to\infty$,
$$\frac{\Fbar_s - \phibar_s}{\sqrt{\Var(F_s)}}\rightarrow 0\mbox{ a.s.} $$
\end{lemma}

\begin{lemma}\label{Lemma:Var}
{If $a+b<2/d_+^{s}$, then for all star motif $s$ and under {the \BEDD model},}  we have, as $m\sim n\to\infty$,
$$\widehat{\Var}(F_s)/\Var(F_s)\rightarrow 1\mbox{ a.s.}$$
\end{lemma}

\subsection{Goodness-of-fit tests for the \BEDD model}  \label{sec:test}

{
We consider a bipartite network $\Gcal$ and we want to test if it arises from the \BEDD model:
$$
\left\{\begin{array}{lcl}
H_0 : \Gcal\mbox{ follows a \BEDD model}, \\
H_1 : \Gcal\mbox{ does not follow a \BEDD model}.
\end{array}\right.
$$
}

{To this aim, we consider the test statistic $W_s = (F_s - \Fbar_s)/\sqrt{\widehat{\Var}(F_s)}$ defined in \eqref{eq::Ws}. The idea is thus to compare the frequency of a motif observed in the network with its expected value under the \BEDD model.}

\begin{remark}
We can consider more specific hypothesis. Suppose we want to test the top node heterogeneity under \BEDD, more specifically $H_0 : \Gcal\mbox{ follows a \BEDD model}$ and $\gt$ is constant. Then, according to \eqref{eq:int_g_h}, we have that $\st_d = \rho^d$ under $H_0$, so a similar statistic to $W_s$ can be designed by considering 
$\Fbar_s 
= \prod_{u = 1}^{\pt_s}  F_1^{\dt^s_u} \prod_{v = 1}^{\pb_s} \Sb_{\db^s_v}/F_1^{d_+^s}$.
In the same manner, a statistic can be designed to test the bottom node heterogeneity.
\end{remark}

\paragraph{Illustration.}
Table~\ref{tab:tests} gives the test statistics $W_s$ for goodness of fit to the \BEDD model for the same motifs as in  Table~\ref{tab:freqs}. According to Theorem~\ref{Thm:AsympNorm}, these statistics should be compared with the quantiles of standard normal distribution $\Ncal(0, 1)$. Almost no motif frequency displays a significant deviation from its expectation under the \BEDD model. Only motif 16 in the plant-pollinator network displays a higher frequency than expected under \BEDD (with $p$-value 7.5\;$10^{-3}$).

\begin{table}[ht]
  \begin{center}
    \begin{tabular}{clllll}
      \multicolumn{6}{c}{plant-pollinator} \\ 
      \hline
      $s$ & \quad 5 & \quad 6 & \quad 10 & \quad 15 & \quad 16 \\
      \hline
      $W_s$ & 
      -6.45\;$10^{-2}$ & 9.96\;$10^{-1}$ & -6.63\;$10^{-2}$ & 7.52\;$10^{-1}$ & 2.43 \\ 
      \\
      \multicolumn{6}{c}{seed dispersal} \\
      \hline
      $s$ & \quad 5 & \quad 6 & \quad 10 & \quad 15 & \quad 16 \\
      \hline
      $W_s$ & 
      -2.14\;$10^{-1}$ & -2.14\;$10^{-1}$ & -2.93\;$10^{-1}$ & -2.95\;$10^{-1}$ & -3.56\;$10^{-1}$
    \end{tabular}
    \caption{Test statistics $W_s$ for the goodness-of-fit of \BEDD for the same motifs as in Table~\ref{tab:freqs}. } \label{tab:tests}
  \end{center}
\end{table}

\subsection{Tests for the comparison of two bipartite networks} \label{sec:comparison}

This section is devoted to network comparison test. More specifically, considering two networks assumed to arise from two \BEDD models, we want to test if they arise from the same \BEDD model, or for, instance, from two different \BEDD model with same function $\gt$. The rational behind the tests we propose is to compare the frequency of a motif observed in one network with its expected value according to the parameters of the other network.
To this aim, we need to introduce specific notations.

\paragraph{Notations.} The \BEDD model is parametrized with the $(m, n, \rho, g, h)$ but all moments depend on $(m, n, \rho, \st, \sb)$, where $\st$ (resp. $\sb$) stands for the sequence of occurrence probability of all the top (resp. bottom) star motifs. In the sequel we denote by $E_s$ the expected frequency of motif $s$:
$$
E_s(m, n, \rho, \st, \sb) := \phi_s,
$$
so its plug-in estimate is $E_s(m, n, F_1, \St, \Sb) = \Fbar_s$.
Similarly, we denote the variance of the frequency by $V_s(m, n, \rho, \st, \sb) := \Var(F_s)$ and its plug-in estimate $V_s(m, n, F_1, \St, \Sb) := \widehat{\Var}_s(F_s)$.

\paragraph{A global test.} We consider two bipartite networks $\Gcal^A$ and $\Gcal^B$ supposed to arise from \BEDD models with respective dimensions and parameters $(m^A, n^A, \rho^A, \st^A, \sb^A)$ and $(m^B, n^B, \rho^B, \st^B, \sb^B)$. We want to test
$$
\left\{\begin{array}{lcl}
H_0 : \left\{(\rho^A, \gt^A, \gb^A) = (\rho^B, \gt^B, \gb^B)\right\}, \\
H_1 : \left\{\rho^A \neq \rho^B \text{ or } \gt^A \neq \gt^B \text{ or } \gb^A \neq \gb^B \right\}.
\end{array}\right.
$$
This is to test that, although the two networks may have different dimensions ($m, n$), they have the same density ($\rho$), the same top node heterogeneity ($\gt$) and the same bottom node heterogeneity ($\gb$).

\paragraph{Test statistics.} The test statistic is based on $F_s^A$ and $F_s^B$ the empirical frequencies of motif $s$ in $\Gcal^A$ and $\Gcal^B$ respectively. The superscript $A$ (resp. $B$) is added to all quantities observed in $\Gcal^A$ (resp. $\Gcal^B$).
\begin{equation} \label{eq:statComparison}
W_s
= \frac{
\left(F^A_s - E_s(m^A, n^A, F_1^B, \St^B, \Sb^B)\right) - 
\left(F^B_s - E_s(m^B, n^B, F_1^A, \St^A, \Sb^A)\right)}{\sqrt{
V_s(m^A, n^A, F_1^B, \St^B, \Sb^B) 
+ V_s(m^B, n^B, F_1^A, \St^A, \Sb^A)}}.
\end{equation}

\begin{theorem}
If both $m^A/m^B$ and $n^A/n^B$ tends to constants, 
if $a+b<2/d_+^{s}$, then for all non-star motif $s$ and under $H_0$, the statistic $W_s$ is asymptotically normal as $m^A\sim n^A\sim m^B\sim n^B\to\infty$:
$$
W_s\overset{D}{\longrightarrow} \Ncal(0, 1).
$$
\end{theorem}

\proofbegin 
We decompose 
\begin{align*}
F^A_s - E_s(m^A, n^A, F_1^B, \St^B, \Sb^B) & = L_s^A + C_s^A \\
\text{where } \qquad L_s^A & = F^A_s - E_s(m^A, n^A, \phi_1^B, \st^B, \sb^B) \\
\text{and } \qquad C_s^A & = E_s(m^A, n^A, \phi_1^B, \st^B, \sb^B) - E_s(m^A, n^A, F_1^B, \St^B, \Sb^B). 
\end{align*}Because $(m^A, n^A)$ go to infinity at the same speed as $(m^B, n^B)$, under $H_0$, $L_s^A / V_s(m^A, n^A, \phi_1^B, \st^B, \sb^B)$ is asymptotically normal according to Proposition \ref{Prop:Ls}, whereas $C_s^A / V_s(m^A, n^A, \phi_1^B, \st^B, \sb^B)$ tends to zero according to Lemma \ref{Lemma:Cs}. Using the same arguments for the symmetric term, we get that and the negligible one $\left(C_s^A / V_s(m^A, n^A, \phi_1^B, \st^B, \sb^B),C_s^B / V_s(m^A, n^A, \phi_1^A, \st^A, \sb^A)\right)$, replacing $V_s(m^A, n^A, \phi_1^B, \st^B, \sb^B)$ and $V_s(m^B, n^B, \phi_1^A, \st^A, \sb^A)$ with their plug-in estimate $V_s(m^A, n^A, F_1^B, \St^B, \Sb^B)$ 
and $V_s(m^B, n^B, F_1^A, \St^A, \Sb^A)$. We conclude using Lemma \ref{Lemma:Var} and Slutsky Theorem.
\proofend 

\paragraph{Testing equal top nodes heterogeneity.} 
Suppose we want to test that, although the two networks may have different dimensions, different densities, and different bottom node heterogeneity, they have the same top node heterogeneity, that is
$$
\left\{\begin{array}{lcl}
H_0 : \left\{\gt^A = \gt^B \right\}, \\
H_1 : \left\{\gt^A \neq \gt^B  \right\}.
\end{array}\right.
$$

{
Since we allow the two networks to have different densities, one might normalize the probabilities of star motifs given in \eqref{eq:lambdagamma} as follows:
\begin{align*}
 \widetilde{\st}_1 & = 1, &  \widetilde{\st}_2 & = \phi_2/ \phi_1^2 &  \widetilde{\st}_3 & = \phi_7/ \phi_1^3  &  \widetilde{\st}_4 & = \phi_{17}/ \phi_1^4 , &  \widetilde{\st}_5 & = \phi_{44}/ \phi_1^4, \\
   \widetilde{\sb}_1 & = 1, &  \widetilde{\sb}_2 & = \phi_3/ \phi_1^2 , &  \widetilde{\sb}_3 & = \phi_4/ \phi_1^3  &  \widetilde{\sb}_4 & = \phi_8/ \phi_1^4, &  \widetilde{\sb}_5 & = \phi_{18}/ \phi_1^4. \nonumber
\end{align*} 
This allows to see that we can rewrite $E_s(m, n, \rho, \st, \sb)=\phi_s$ as an expression of $\gt$ on which relies the test we consider. According to \eqref{eq:phi} and to the definition of $\phi_s $ under the \BEDD model, we get:
  \begin{equation*}
    E_s(m, n, \rho, \gt, \gb) 
     =  \rho^{d_+^s}  {\prod_{u = 1}^{\pt_s}  \tilde{\st}_{\dt^s_u} \prod_{v = 1}^{\pb_s} \tilde{\sb}_{\db^s_v}} \\
     = \rho^{d_+^s} \prod_{u = 1}^{\pt_s} \prod_{v = 1}^{\pb_s} \gt_{\dt^s_u} \gb_{\db^s_v},
  \end{equation*}
where $\gt_\dt = \int \gt(u)^\dt \d u$ and $\gb_\db = \int \gb(v)^\db \d v$. 
We may consider the following test statistic:
$$
W^\gt_s
= \frac{
\left(F^A_s - E_s(m^A, n^A, F_1^A, \widetilde{\St}^B, \widetilde{\Sb}^A)\right) - 
\left(F^B_s - E_s(m^B, n^B, F_1^B, \widetilde{\St}^A, \widetilde{\Sb}^B)\right)}{\sqrt{
V_s(m^A, n^A, F_1^A, \widetilde{\St}^B, \widetilde{\Sb}^A) 
+ V_s(m^B, n^B, F_1^B, \widetilde{\St}^A, \widetilde{\Sb}^B)}},
$$ 
where $\widetilde{\St}$ and $\widetilde{\Sb}$ are the plug-in estimates of $\widetilde{\st}$ and $\widetilde{\sb}$ respectively.
Similar statistics can be designed to test $\rho^A = \rho^B$, $\gb^A = \gb^B$ or any combination.
}

\paragraph{Illustration.} Both the plant-pollinator and the seed dispersal networks involve plants species. Although these species are not the same, one may be interested in comparing if the level of heterogeneity across plants (encoded in the function $\gt$) is the same in both networks. From an ecological point of view, this amounts to test if there is the same the degree of imbalance between specialists and generalists among plants regarding pollination and seed dispersion, that are two of the main reproduction means. \\
Table \ref{tab:testsComp} provides the results of the network comparison test presented above. No significant difference is observed, suggesting that, although generalist and specialist plants may exist for both types of interactions (no assumption is made about the shape of $g$), the degree of imbalance between them is comparable ($g^A \simeq g^B$).

\begin{table}[ht]
  \begin{center}
    \begin{tabular}{clllll}
        $s$ & \quad 5 & \quad 6 & \quad 10 & \quad 15 & \quad 16 \\
        \hline
        $F_s^A$ &  9.21\;$10^{-5}$ &  1.00\;$10^{-5}$ &  8.12\;$10^{-6}$ &  3.32\;$10^{-7}$ &  4.47\;$10^{-8}$ \\
        $\displaystyle{\widehat{\Esp}_0 F_s^A}$ &  1.96\;$10^{-4}$ &  3.75\;$10^{-5}$ &  1.74\;$10^{-5}$ &  4.25\;$10^{-6}$ &  1.33\;$10^{-6}$ \\
        \hline
        $F_s^B$ &  5.13\;$10^{-4}$ &  1.15\;$10^{-4}$ &  5.07\;$10^{-5}$ &  1.79\;$10^{-5}$ &  5.96\;$10^{-6}$ \\
        $\displaystyle{\widehat{\Esp}_0 F_s^B}$ &  2.66\;$10^{-4}$ &  2.92\;$10^{-5}$ &  2.85\;$10^{-5}$ &  1.50\;$10^{-6}$ &  1.69\;$10^{-7}$ \\
        \hline
        $F_s^B - F_s^A$ & -4.21\;$10^{-4}$ & -1.05\;$10^{-4}$ & -4.26\;$10^{-5}$ & -1.76\;$10^{-5}$ & -5.91\;$10^{-6}$ \\
        $\displaystyle{\widehat{\Esp}_0(F_s^B - F_s^A)}$ & -6.96\;$10^{-5}$ & 8.37\;$10^{-6}$ & -1.11\;$10^{-5}$ & 2.75\;$10^{-6}$ & 1.16\;$10^{-6}$ \\
        $\displaystyle{\sqrt{\widehat{\Var}_0(F_s^A) + \widehat{\Var}_0(F_s^B)}}$ & 2.25\;$10^{-4}$ & 7.24\;$10^{-5}$ & 3.26\;$10^{-5}$ & 1.59\;$10^{-5}$ & 7.38\;$10^{-6}$ \\
        \hline
        $W_s$ & -1.56 & -1.56 & -0.97 & -1.28 & -0.96
    \end{tabular}
    \caption{Network comparison test for $H_0 = \{\gt^A = \gt^B\}$ as defined in Section \ref{sec:comparison} for the same motifs as in Table~\ref{tab:freqs}. Networks: $A=$ plant-pollinator, $B=$ seed dispersal. $\widehat{\Esp}_0(\cdot)$ is a shorthand for the notation $E_s(\cdots)$ (idem for $\widehat{\Var}_0(\cdot)$ and $V_s(\cdots)$).} \label{tab:testsComp}
  \end{center}
\end{table}

\section{Simulation study}  \label{sec:simul}

{We designed a simulation study to illustrate Theorem \ref{Thm:AsympNorm} and to assess the performance of the goodness-of-fit test and the comparison test described in Section \ref{sec:test} and Section \ref{sec:comparison} respectively. More specifically, our purpose is to illustrate the asymptotic normality of the test statistics and evaluate the power of the tests for various graph sizes, densities and sparsity regimes. }

\subsection{Asymptotic normality} \label{sec:asympNormality}

\paragraph{Simulation design.}
{
We simulated series of networks with parameters ($m, n, \rho, \mut, \mub$) varying according to the following design:
\begin{description}
\item[Network dimension:] We simulated networks with equal dimensions $m=n$, with values in $\{50, 100, 200, 500, 1000, 2000\}$; 
\item[Sparsity regime:] We considered equal parameters $a=b$ in $\{1/3, 1/4, 1/5, 1/6\}$; 
\item[Network density: ] The resulting density is $\rho = \rho_0 m^{-a} n^{-b}$, $\rho_0$ being fixed so that $\rho = .01$ when $m = n = 100$; 
\item[Degree imbalance:] We considered the functions $\gt(u) = \mut u^{\mut-1}$ and $\gb(v) = \mub v^{\mub-1}$;  observe that $\mut = 1$ means that $\gt$ is constant so no imbalance does exist top nodes (resp. for $\mub$, $\gb$ and bottom nodes). We set $\mut = 2$, $\mub = 3$.
\end{description}
}
For each configuration, $S=100$ networks were sampled and the test applied.

\begin{figure}[ht]
    \centering
    \begin{tabular}{cc|cc}
    motif 6 & motif 15 & motif 5 & motif 10 \\
    \includegraphics[width=.25\textwidth, trim=10 10 10 10]{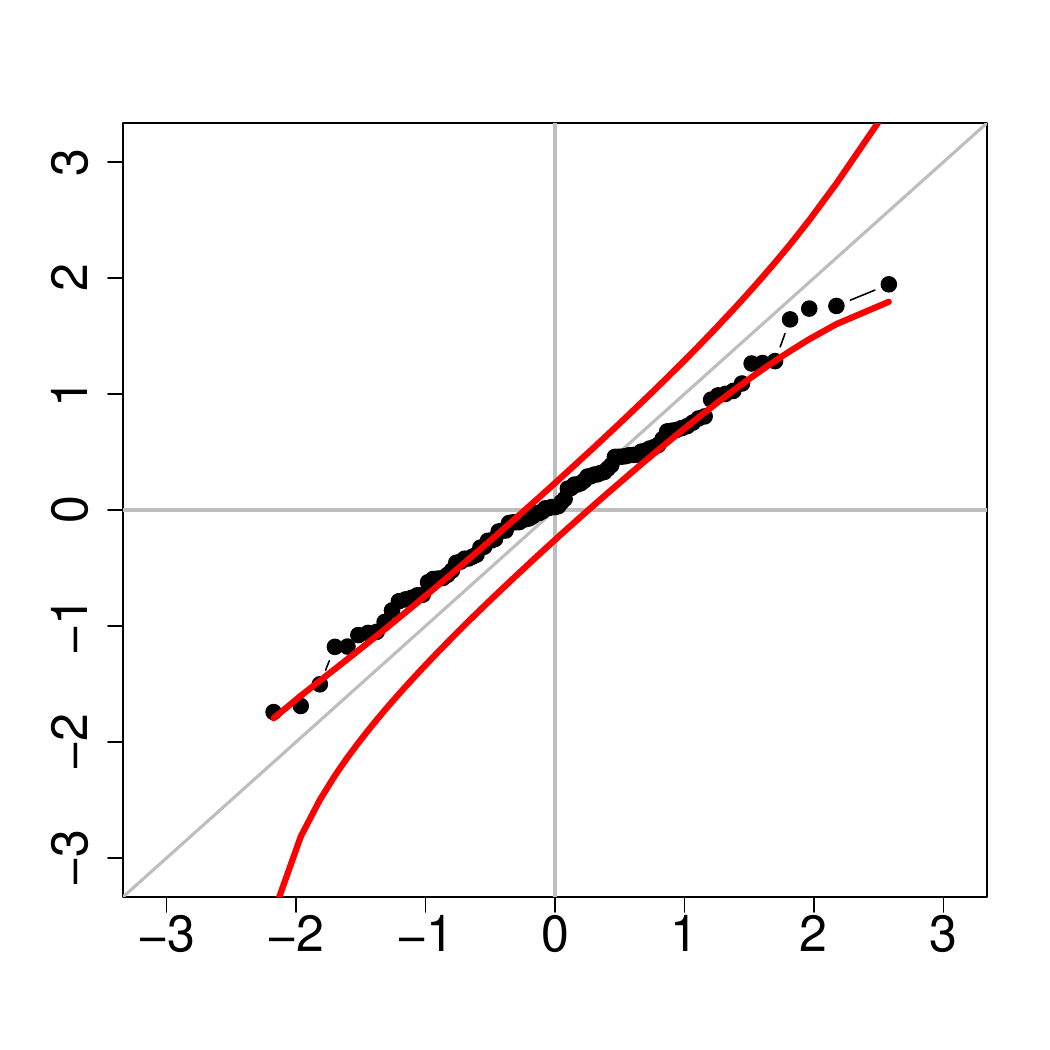} & 
    \includegraphics[width=.25\textwidth, trim=10 10 10 10]{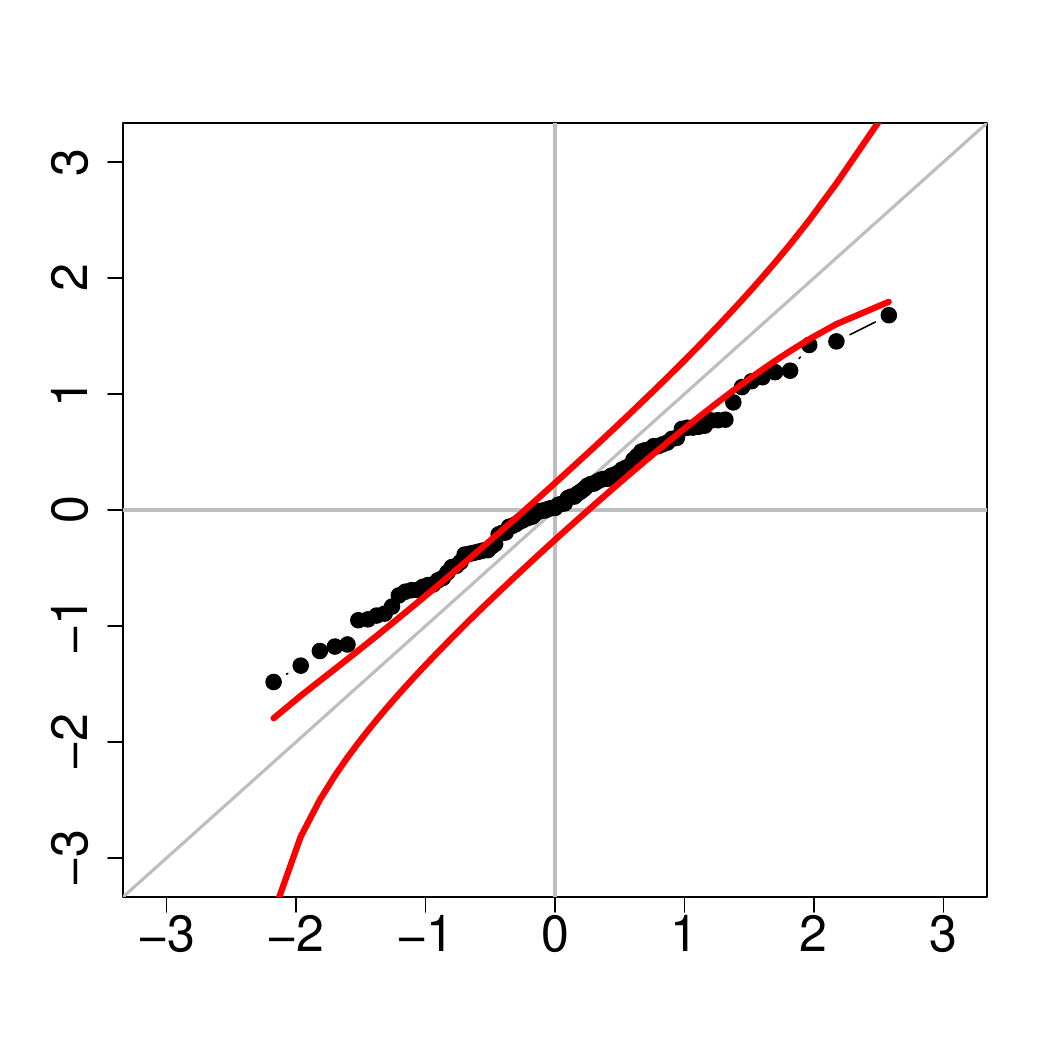} &
    \includegraphics[width=.25\textwidth, trim=10 10 10 10]{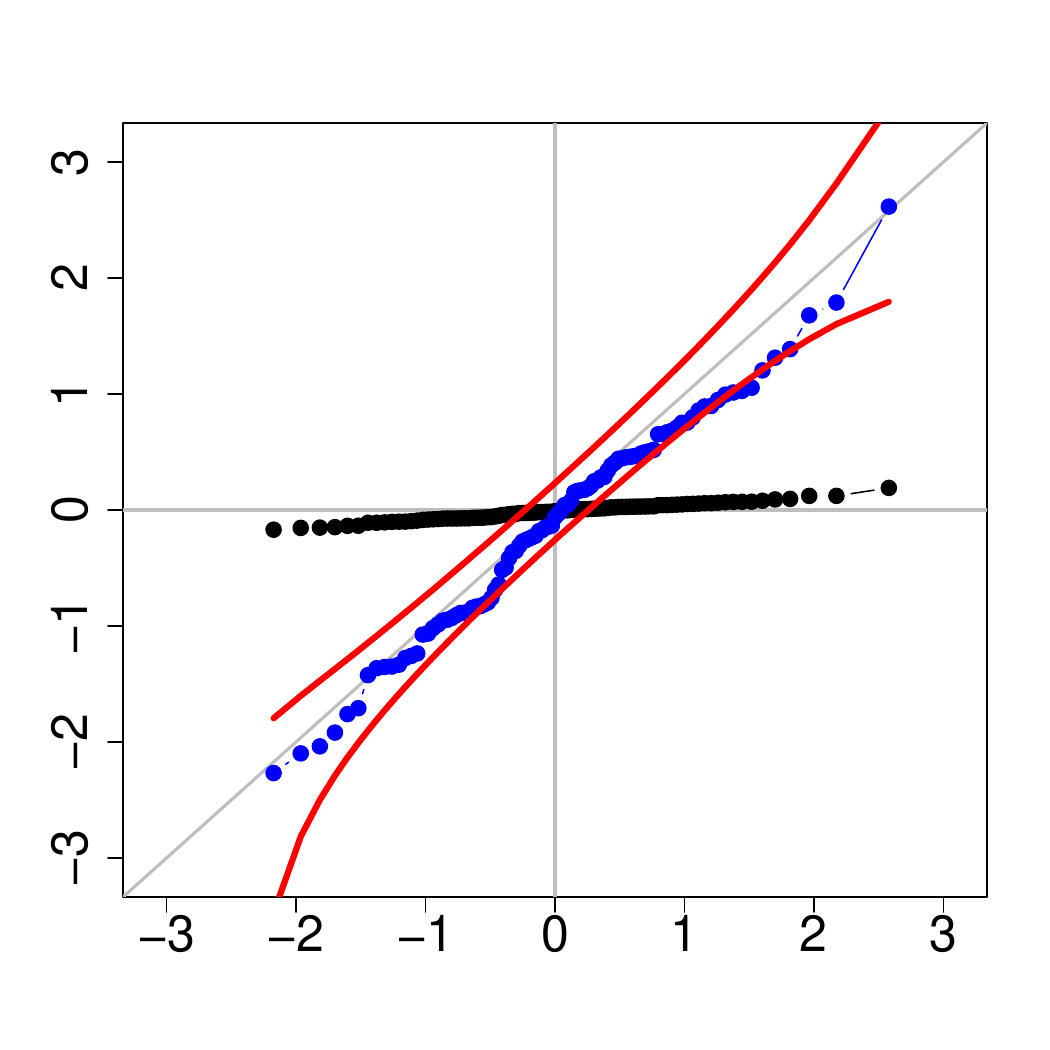} & 
    \includegraphics[width=.25\textwidth, trim=10 10 10 10]{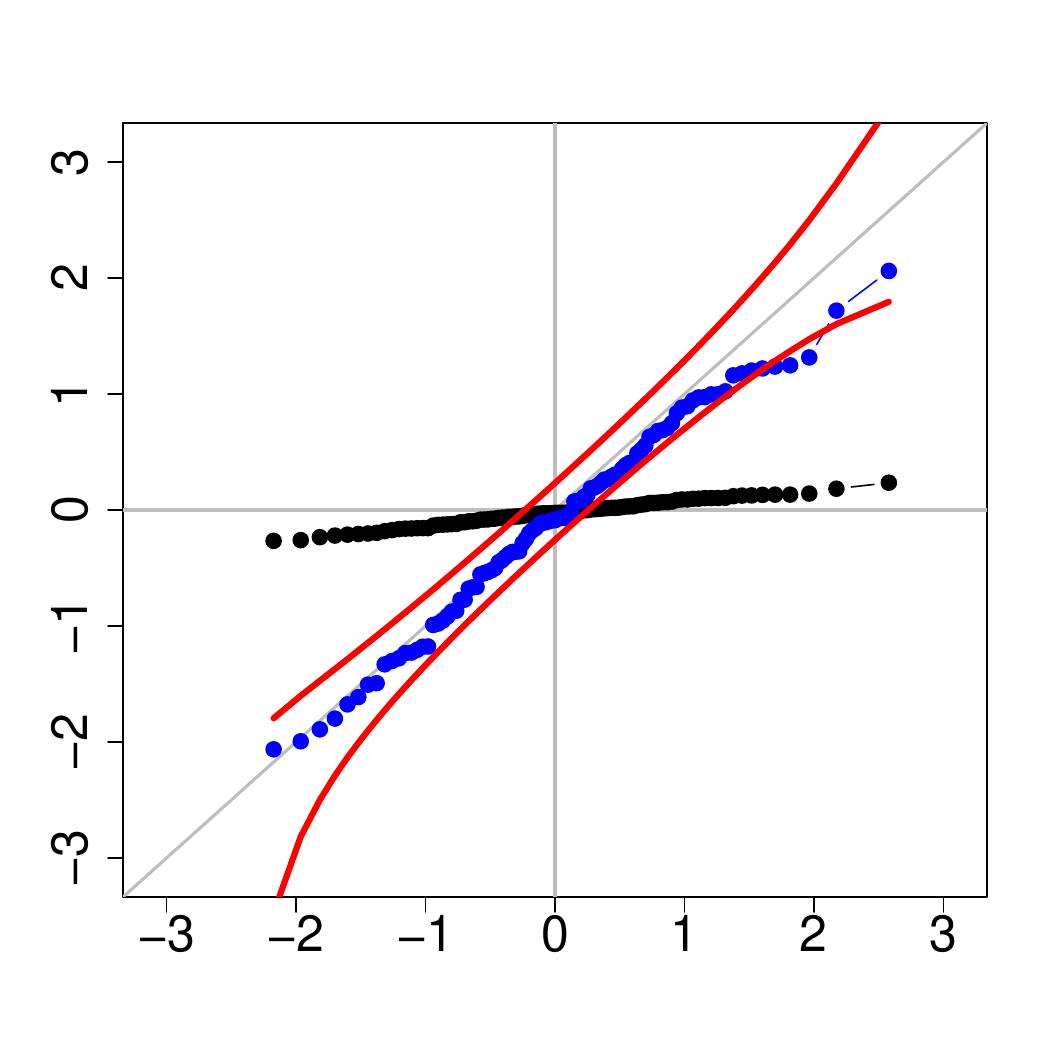}
    \end{tabular}
    \caption{Qq-plots of the test statistics ${W}_s$ for 4 motifs in a network with dimension $m = n = 2000$ and sparsity regime $a=1/3$ (black dots). Blue dots: qq-plot for the corrected statistic $\widetilde{W}_s$ defined in Equation \eqref{eq:corrected:stat}. Red line: 95\% confidence interval for a qq-plot with sample size $S=100$.}
    \label{fig:rawAsymptotic}
\end{figure}

\begin{figure}[ht]
    \centering
    \begin{tabular}{cccc}
    $\Esp(N_s) < 1$ & $1 \leq \Esp(N_s) < 10$ &
    $10 \leq \Esp(N_s) < 100$  & $\Esp(N_s) \geq 100$ \\
    \includegraphics[width=.25\textwidth, trim=10 10 10 10]{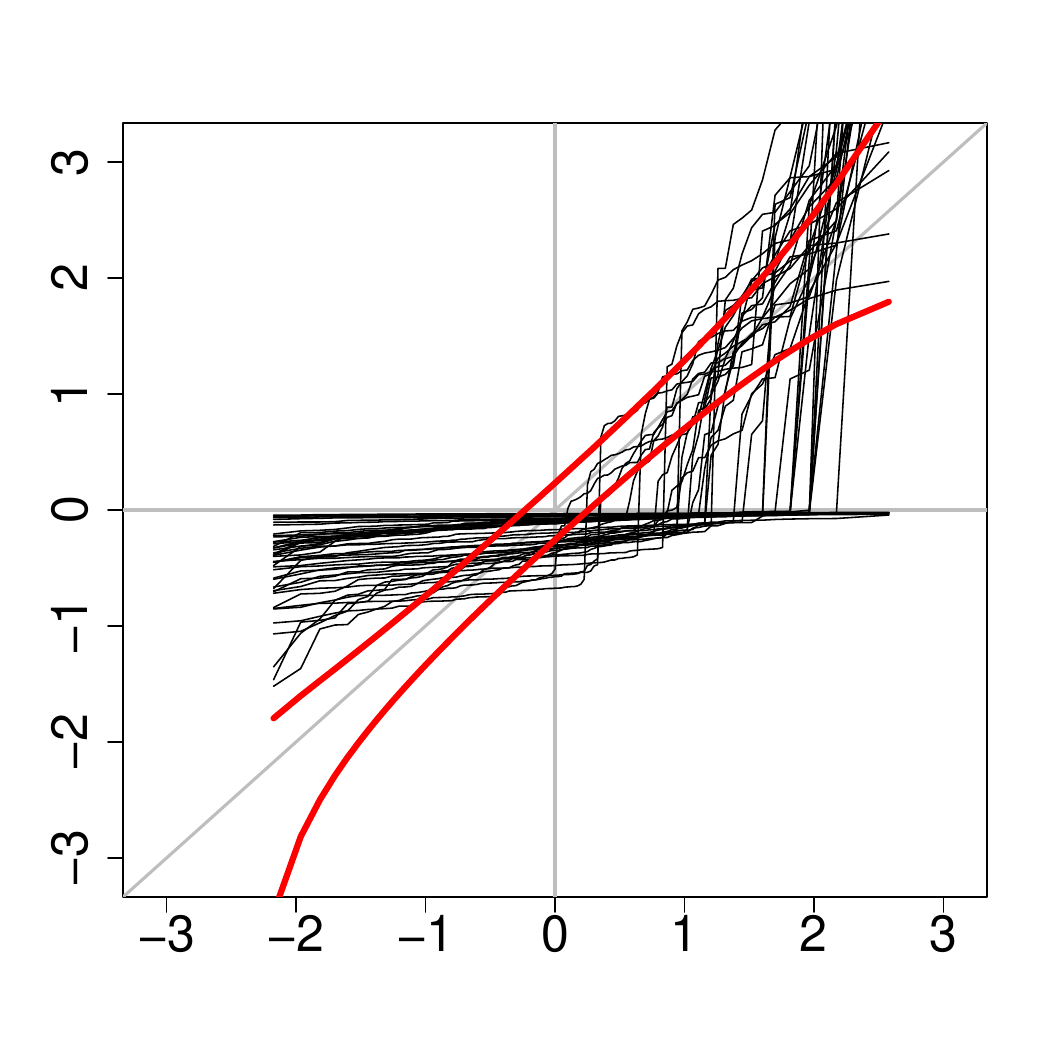} & 
    \includegraphics[width=.25\textwidth, trim=10 10 10 10]{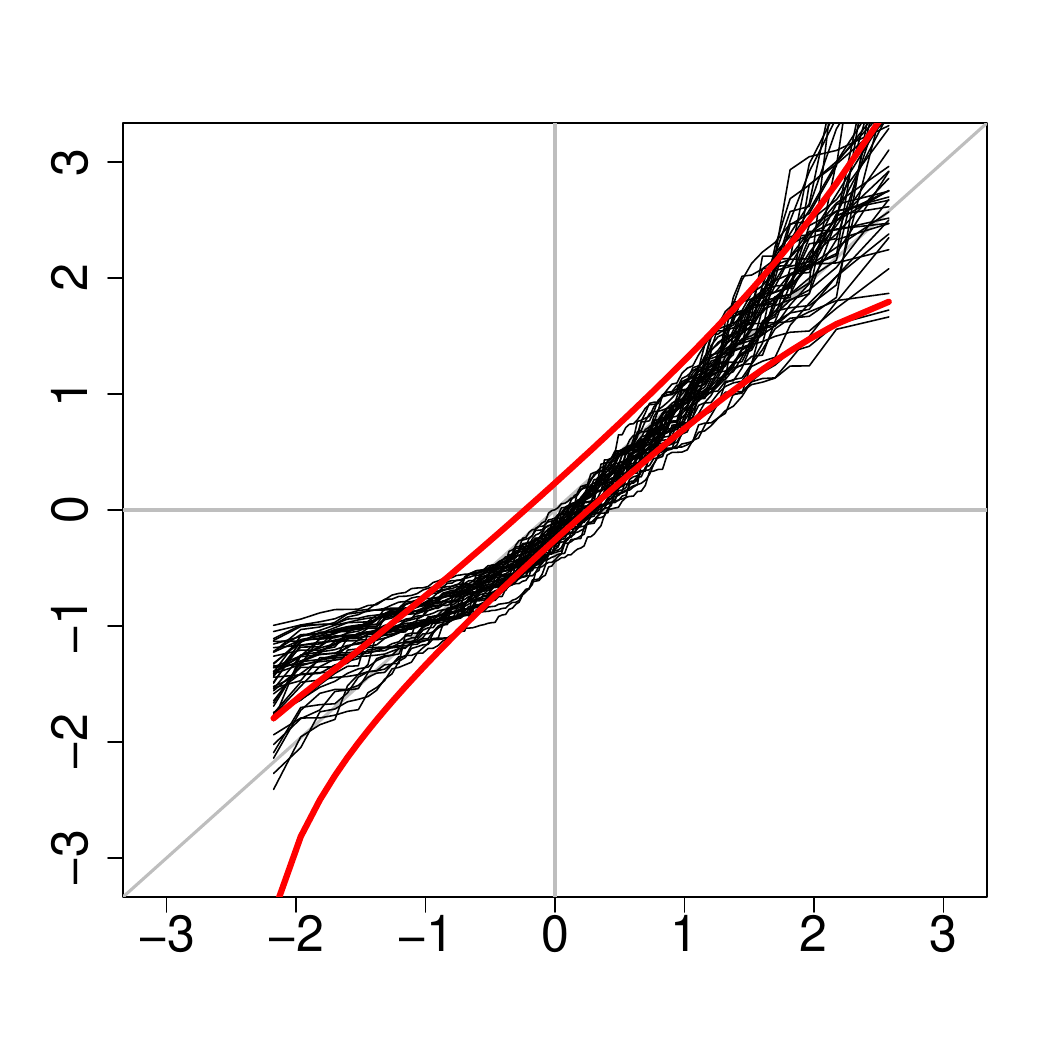} & 
    \includegraphics[width=.25\textwidth, trim=10 10 10 10]{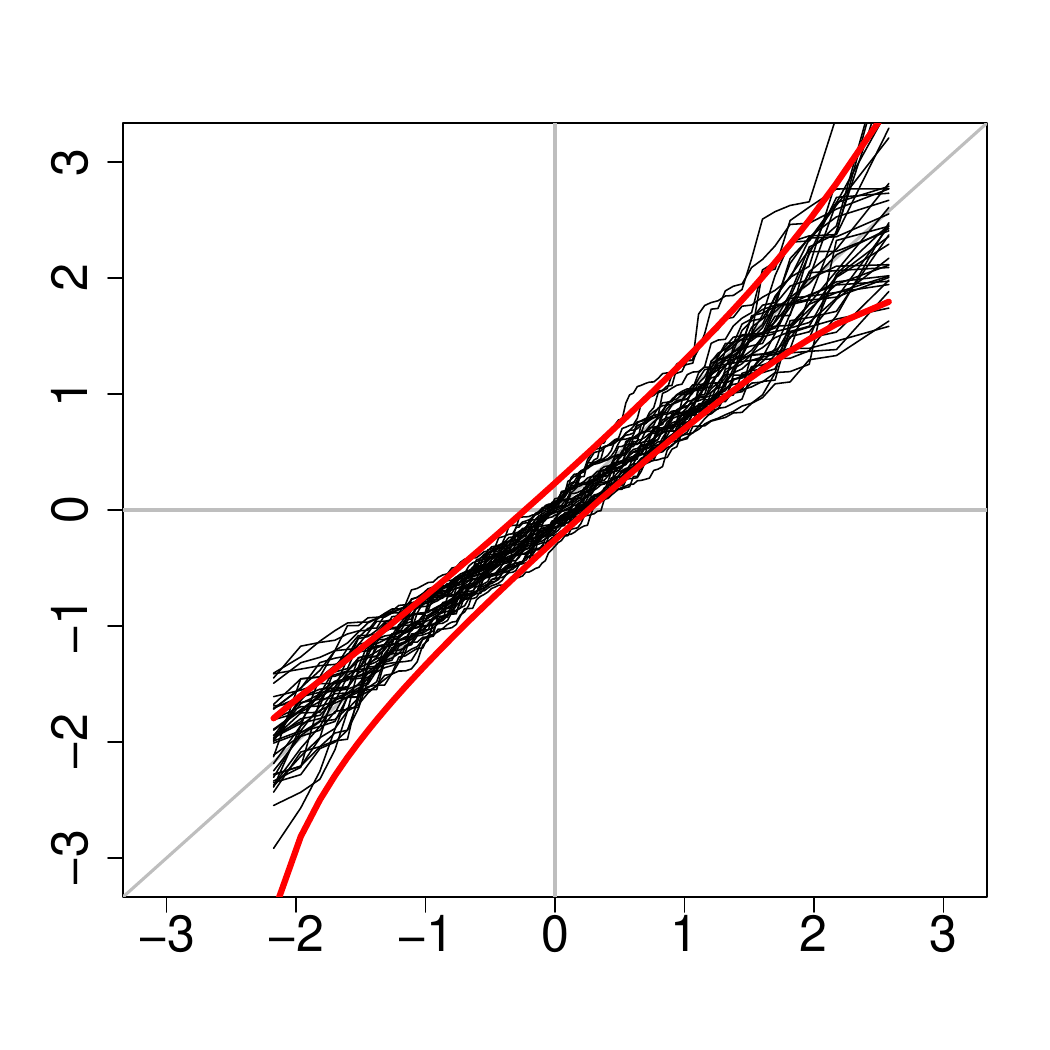} & 
    \includegraphics[width=.25\textwidth, trim=10 10 10 10]{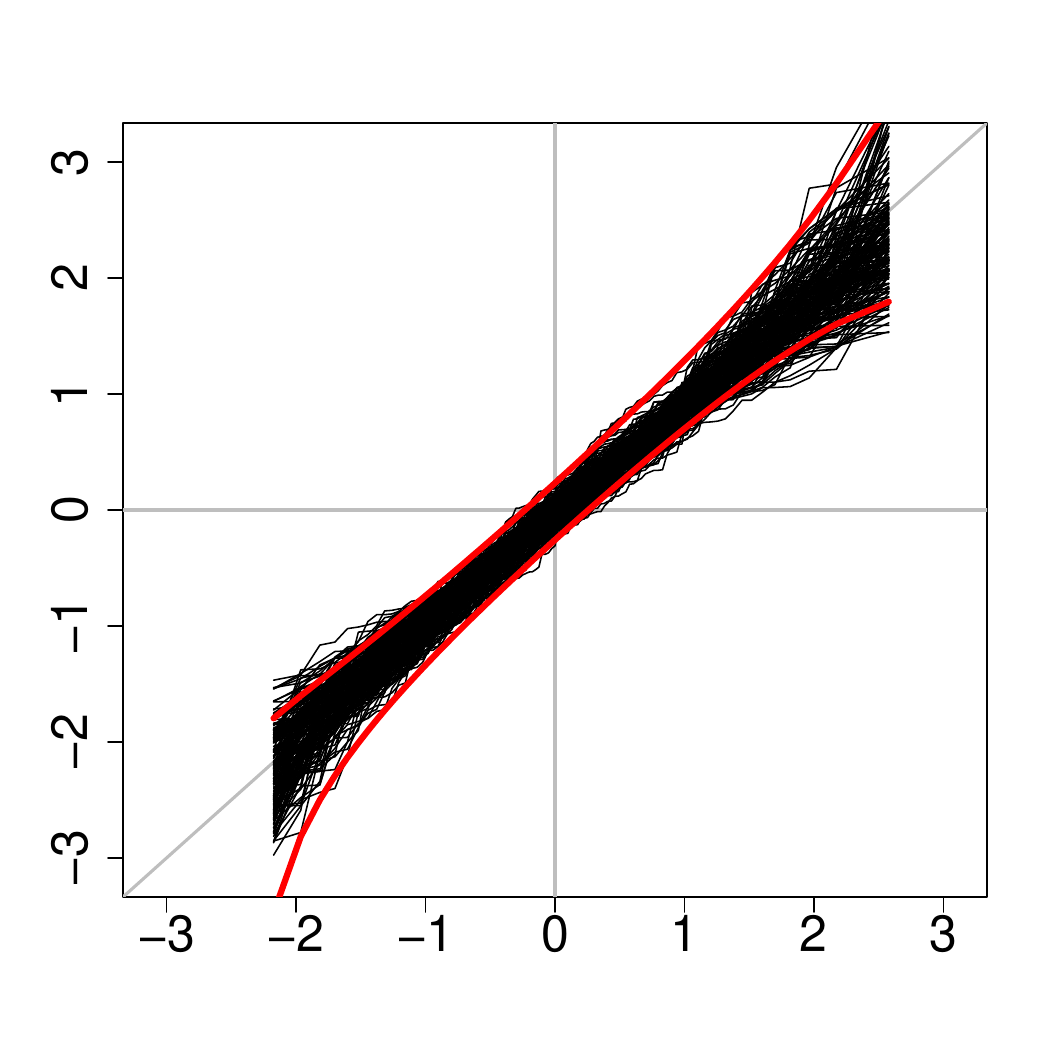} \\
    \end{tabular}
    \caption{Qq-plots of the corrected test statistics $\widetilde{W}_s$. The plot displays the results of the simulation design (i.e for all network size $n$, sparsity regime $a$ and non-star motifs $s$. The qq-plots are gathered accorded to the order of magnitude of the expected count $\Esp(N_s)$, from the smallest (top left) to the largest (bottom right). Red line: same legend as Figure \ref{fig:rawAsymptotic}.}
    \label{fig:unbiasedAsymptotic}
\end{figure}

\paragraph{Results.}
{The results are displayed in Figure \ref{fig:rawAsymptotic} and Figure \ref{fig:unbiasedAsymptotic}. In Figure \ref{fig:rawAsymptotic}, the QQ-plots of the $W_s$ statistic (black dots) defined in \eqref{eq::Ws} and the $\widetilde{W}_s$ statistic (blue dots) defined in \eqref{eq:corrected:stat} hereafter, are given for four motifs in a network with dimension $m = n = 1000$ and sparsity regime $a=b=1/3$ . Remember that the larger the power $a$, the sparser the graph. We observe that normality of $W_s$ holds for motifs 6 and 15, but not for motifs 5 and 10. \\
Actually, the latter case is due to the fluctuations of $\Fbar_s$. More specifically, for non-star motifs, $\Fbar_s$ is not an unbiased estimate of $\phibar_s$ and it is not independent from $F_s$. As a consequence, for finite dimensions $m$ and $n$, we both have that $\Esp(\Fbar_s) \neq \phi_s = \Esp(F_s)$ and $\Var(F_s - \Fbar_s) \neq \Var(F_s)$. Both the bias of $\Fbar$: $\Bias(\Fbar_s) = \Esp(\Fbar) - \phi_s$ and the variance of the numerator of $W_s$: $\Var(F_s - \Fbar_s)$ can be estimated using the delta method, which requires the covariance given in Equation~\eqref{eq:covNsNt}. This enables us to define a corrected version $\widetilde W_s$ of the test statistic $W_s$:
\begin{equation}\label{eq:corrected:stat}
    \widetilde W_s := \widehat\Var\left(F_s - \Fbar_s\right)^{-1/2}\left(F_s - \Fbar_s + \widehat\Bias(\Fbar_s) \right),
\end{equation}
where the bias $\widehat\Bias(\Fbar_s)$ and $\widehat\Var\left(F_s - \Fbar_s\right)$ are both plug-in estimates.
}

\paragraph{Illustration.} We provide in Table \ref{tab:testsCor} the values of corrected corrected statistics $\widetilde{W}_s$ for the plant-pollinator and the seed dispersal networks, to be compared with Table \ref{tab:tests}. Observe that the correction does not yield in different conclusions, in terms of fit to the \BEDD model for both networks.

\begin{table}[ht]
  \begin{center}
    \begin{tabular}{clllll}
      $s$ & \quad 5 & \quad 6 & \quad 10 & \quad 15 & \quad 16 \\
      \hline
      plant-pollinator & -0.05 &  1.03  & -0.03 &  0.79 &  2.49 \\
      seed dispersal & -0.17 & -0.14 & -0.20 & -0.19 & -0.22
    \end{tabular}
    \caption{{Corrected test statistics $\widetilde{W}_s$ for the goodness-of-fit of \BEDD for the same motifs as in Table~\ref{tab:freqs}.} } \label{tab:testsCor}
  \end{center}
\end{table}

{Figure \ref{fig:unbiasedAsymptotic} displays the QQ-plots of the corrected test statistics $\widetilde{W}_s$ gathered according to the order of magnitude of the expected motif frequencies. All network sizes, sparsity regimes and non-star motifs are thus considered here together. As expected, the normality becomes more accurate when the motifs frequency increases.}

\subsection{Power of the goodness-of-fit test} \label{sec:power}

\paragraph{Simulation design.} In order to illustrate the power of the goodness-of-fit test, we simulated a series of networks from a mixture of a \BEDD model and a latent block model (LBM) \citep{govaert2008block}, characterizing the presence of clusters of rows and columns in incidence matrices. Thus, a mixing weight $\alpha$ varying from $0$ to $1$ was considered so that $\alpha=0$ corresponds to a \BEDD that is $H_0$. In details, the following simulation setup was investigated:
\begin{description}
\item[Network dimension and density:] We considered dimensions similar to the pollination and seed dispersal binary networks studied in \cite{SCB19}, that is
$m = n \in \{10^1, \dots, 10^3\}$.
To mimic the sparsity of the same networks, we fitted the density via a linear regression and obtained $\log_{10}(\rho) = 0.3457 -0.3958 \log_{10}(mn)$;
\item[\BEDD model:] {We used the same functions $\gt$ and $\gb$ as in Section \ref{sec:asympNormality}, with $\mut = 2$, $\mub = 3$;}
\item[LBM model:] {
We considered $2$ groups in rows and $2$ groups in columns, all groups with proportion $1/2$ and all connection probabilities $\gamma_{k\ell} = C\gamma_{\min}$ for all $1 \leq k, \ell \leq 2$, except $\gamma_{22} = C\gamma_{\max}$, with $C$ set such that  $C(\gamma_{\max} + 3 \gamma_{\min})/4 = 1$. Two regimes were considered: $\gamma_{\max} = 0.95$ (scenario I: easy) and  $\gamma_{\max} = 0.5$ (scenario II: hard);} 
\item[Connection probability:] {We sampled the $\{U_i\}_{1 \leq i \leq m}$ and $\{V_j\}_{1 \leq j \leq n}$ all independently and uniformly over $[0, 1]$, and set the $\{Z_i\}_{1 \leq i \leq m}$ and $\{W_j\}_{1 \leq j \leq n}$ as $Z_i = \Ibb\{U_i > .5\}+1$  and $W_j= \Ibb\{V_j > .5\}+1$. Finally, the edges were sampled with probability 
$$
\Pr\{G_{ij}=1 \mid U_i, V_j\} = \rho \left((1-\alpha) \gt(U_i) \gb(V_j) + \alpha \gamma_{Z_i W_j}\right).
$$}
\end{description}
For each configuration, $S=500$ networks were sampled and the test applied. Again the test corrected  statistic $\widetilde{W}_s$ was used.

\begin{figure}[ht]
    \centering
    \begin{tabular}{cccc}
    $m=n=50$ & $m=n=100$ & $m=n=200$ & $m=n=500$ \\
    \hline
    \includegraphics[width=.25\textwidth, trim=10 10 10 10]{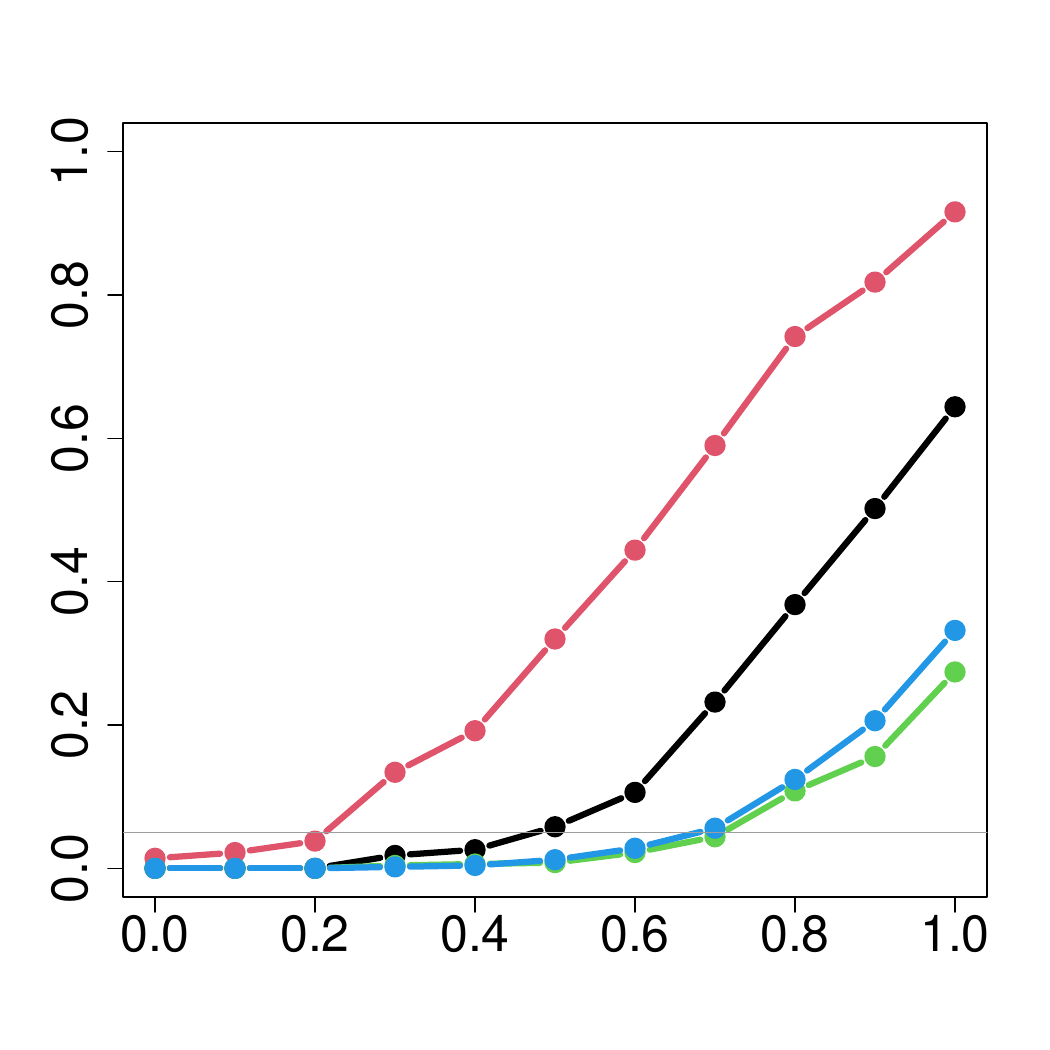} &
    \includegraphics[width=.25\textwidth, trim=10 10 10 10]{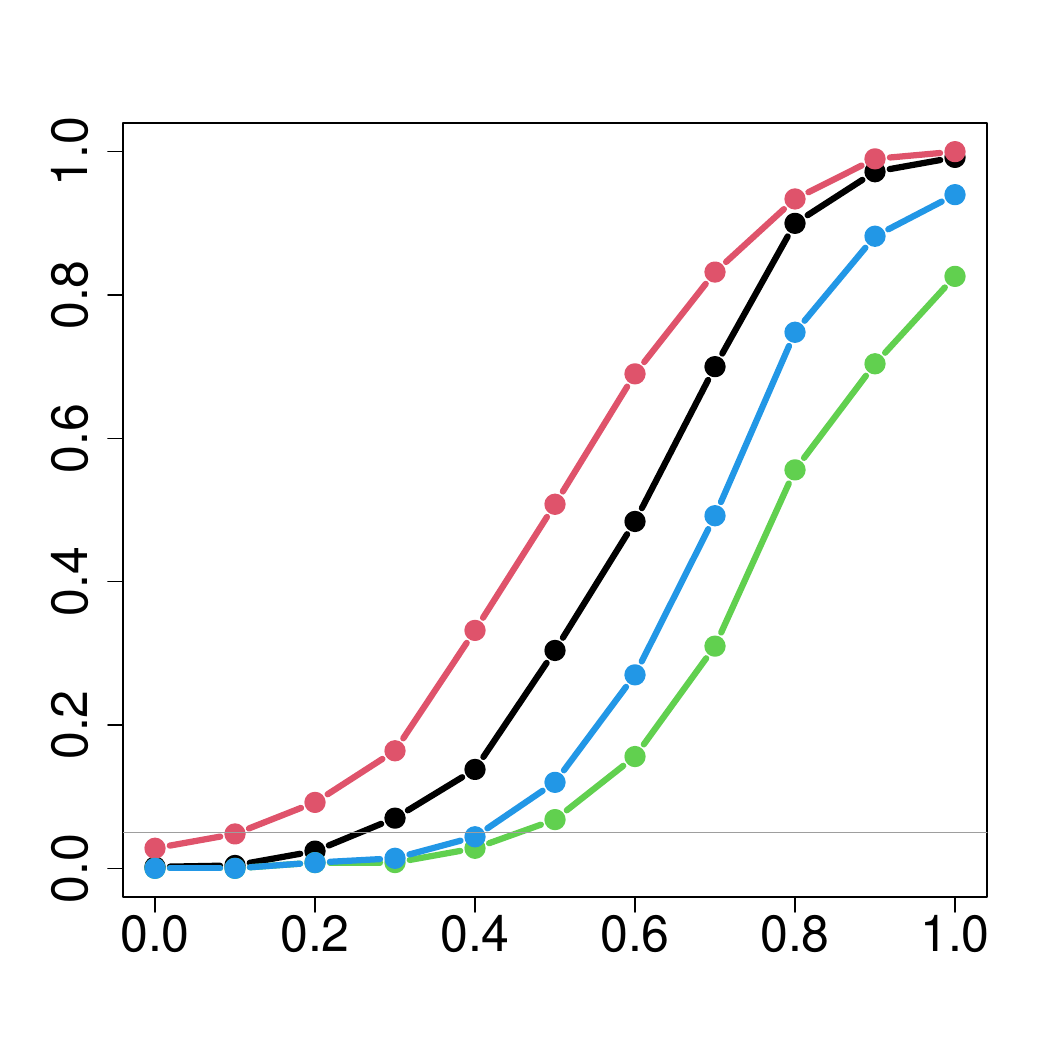} &
    \includegraphics[width=.25\textwidth, trim=10 10 10 10]{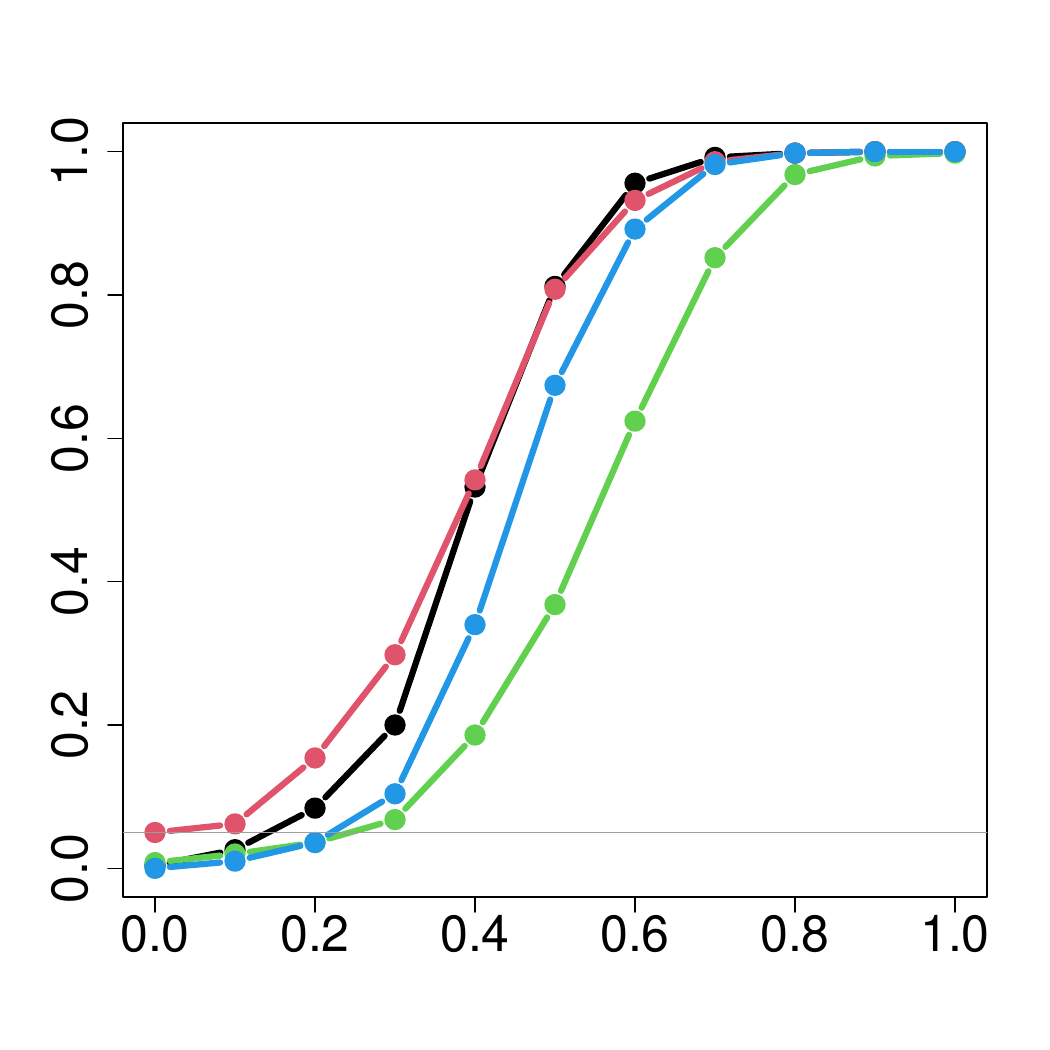} &
    \includegraphics[width=.25\textwidth, trim=10 10 10 10]{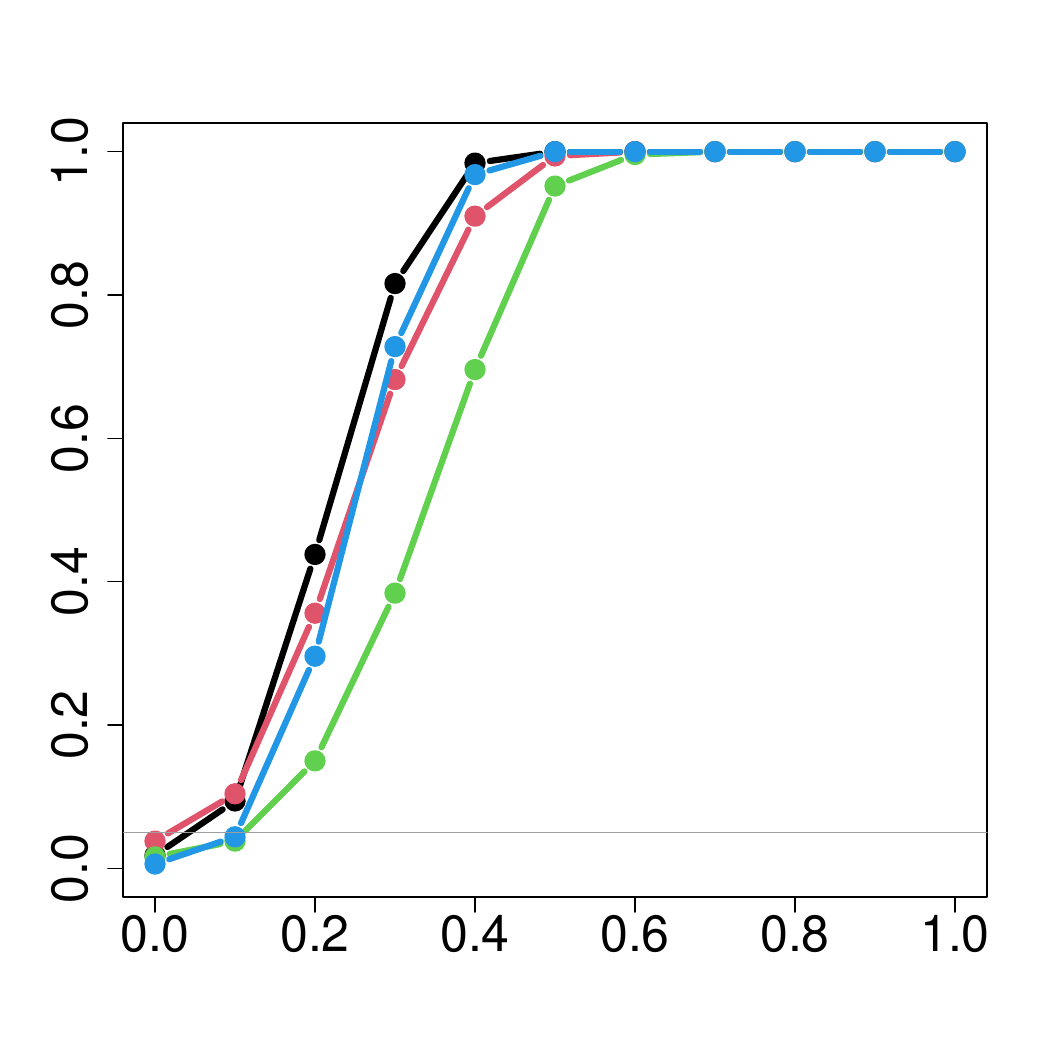} \\
    \includegraphics[width=.25\textwidth, trim=10 10 10 10]{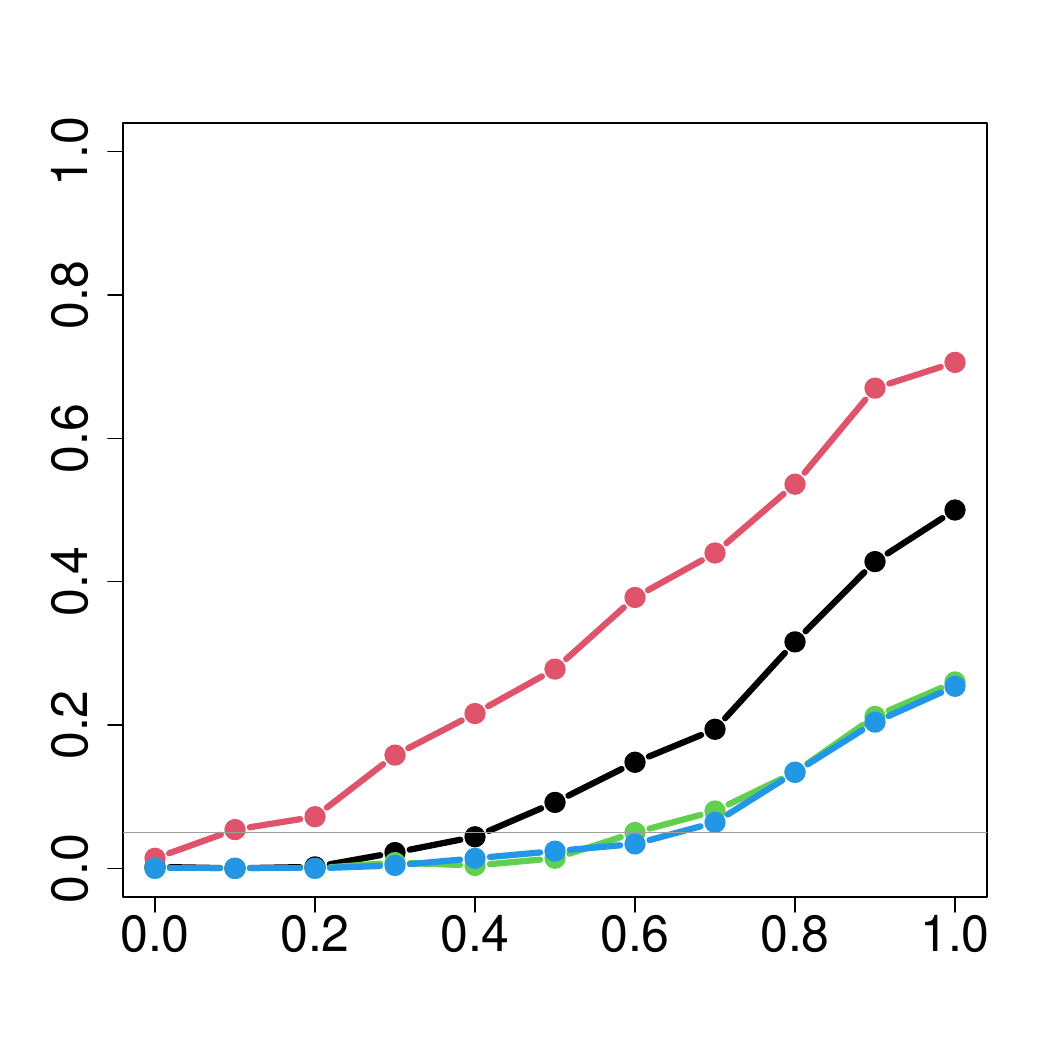} &
    \includegraphics[width=.25\textwidth, trim=10 10 10 10]{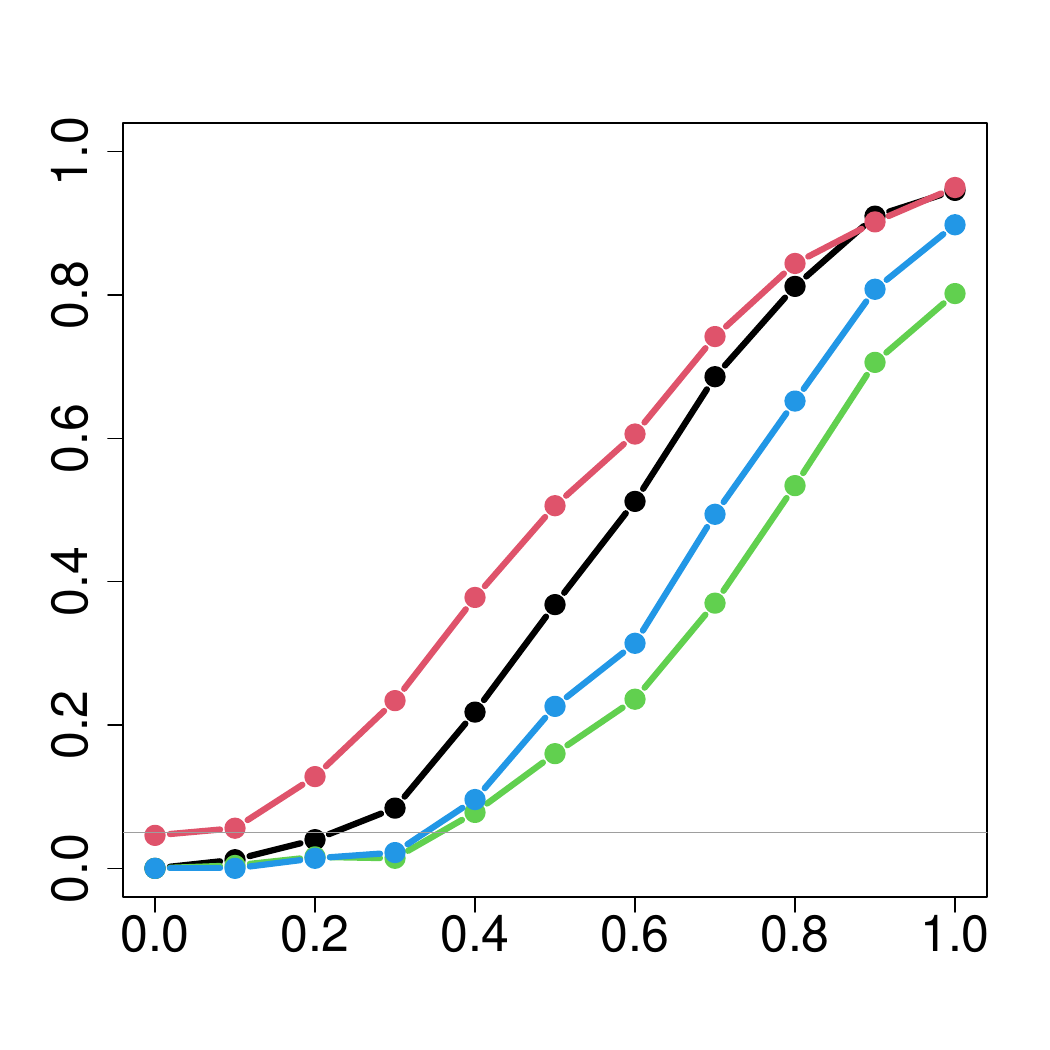} &
    \includegraphics[width=.25\textwidth, trim=10 10 10 10]{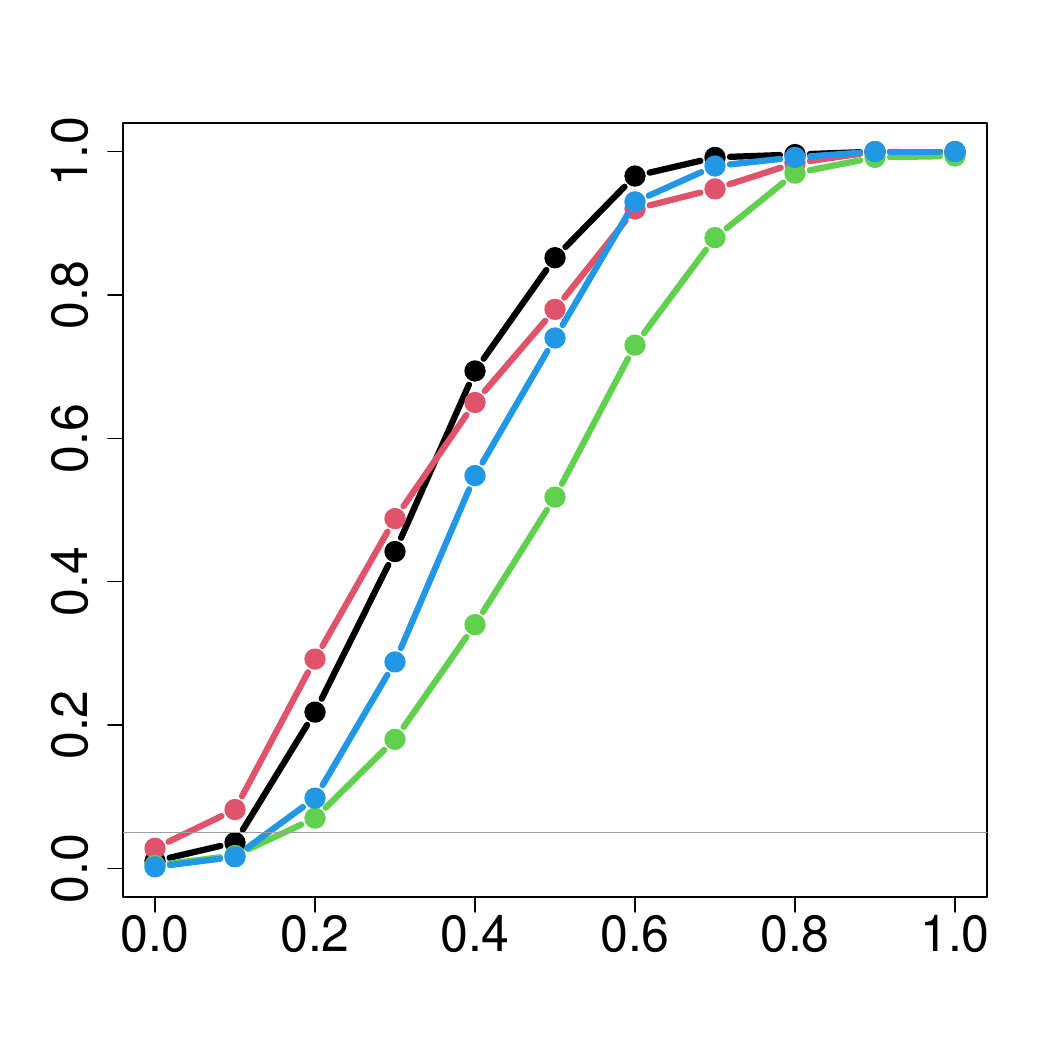} &
    \includegraphics[width=.25\textwidth, trim=10 10 10 10]{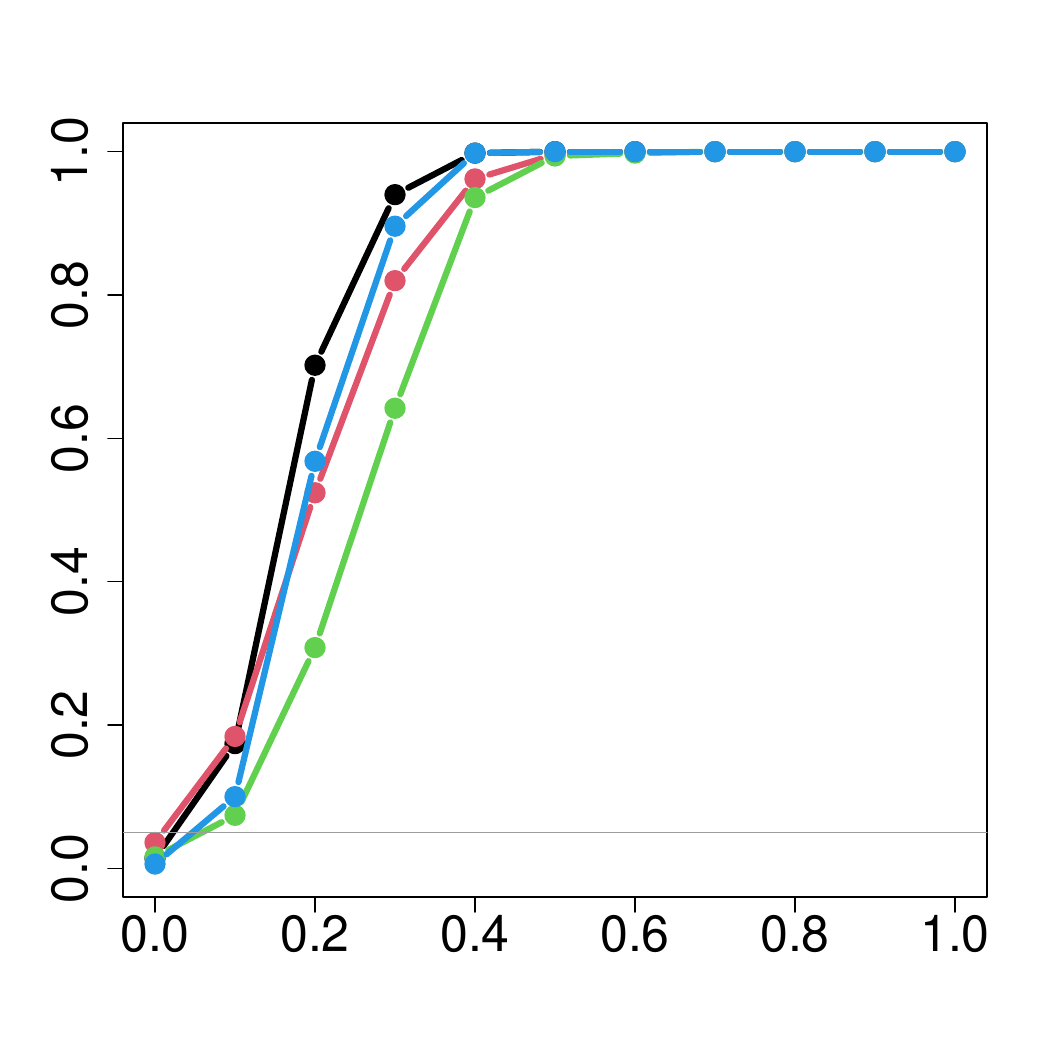}
    \end{tabular}
    \caption{Empirical power of the goodness-of-fit tests, averaged over $S=500$ simulations. Top: scenario I (easy: $\gamma_{\max} = 0.95$); bottom: scenario II (hard: $\gamma_{\max} = 0.5$). From left to right: $m = n = 50, 100, 200, 500$. Color = motif: black=5, red=\textcolor{red}{6}, green=\textcolor{green}{10}, blue=\textcolor{blue}{15}.}
    \label{fig:powerGOF}
\end{figure}

\paragraph{Results.}
The results are given in Figure \ref{fig:powerGOF}. For illustration purposes, we only present the results we obtained for $m=n$ ranging from $50$ to $500$. Moreover, for the sake of clarity, we only consider motifs 5, 6, 10, and 15 which constitute a representative panel of the set of motifs with size $4$ and $5$.\\
As the network dimensions increase, we can clearly observe that the tests become more powerful. For small networks with $m=n=50$ and $m=n=100$, the LBM regime with $\gamma_{\max} = 0.95$ is easier and leads to tests associated with motifs $5$ and $6$ with higher power. These differences vanish for larger values of $n$ and $m$. Overall, we found that motifs $5$ and $6$ lead to more powerful tests. These results illustrate that the methodology proposed is relevant and that the goodness-of-fit tests for different motifs can be used to detect the departure from a \BEDD model. 

\subsection{Power of the network comparison test}

\paragraph{Simulation design.}
We also studied the power of the test for network comparison introduced in Section \ref{sec:comparison}. To this aim, we simulated series of networks $A$ with parameters ($m_A, n_A, \rho^A, \mut^A, \mub^A$) varying according to the same design as in Section \ref{sec:asympNormality}, where $\mut^A$ was set to $2$. \\
We focused on the test of $H_0 = \{\gt^A = \gt^B\}$ so, for each network $A$, we simulated a sequence of networks $B$ with same dimensions ($m_B = m_A$, $n_B = n_A$), but a with a different parameter $\mut^B$. More specifically, setting $\mut^* = 1$ (absence of degree imbalance between top nodes), we sampled networks $B$ with $\mut^B = (1-\alpha) \mu^A + \alpha \mut^*$, with $\alpha = 0, 0.1, 0.2, \dots 1$, so that $\alpha = 0$ corresponds to $H_0$. \\
Regarding the two remaining parameters $\rho^B$ and $\mub^B$, we considered two scenarios:
\begin{description}
\item[I (easy):] $\rho^B = \rho^B$, $\mub^B = \mub^A$, so that the two networks only differ with respect to $\mut$;
\item[II (hard):] $\rho^B = \rho^A/2$, $\mub^B = 2$, so that the two network differ in all parameters, but only the difference in $\mut$ is tested.
\end{description}
The 'hard' scenario is designed to assess the ability of the proposed test statistic to accommodate to differences in density and bottom node imbalance between the two networks, when testing the equality of their top node imbalance. For each configuration, $S=500$ pairs of networks ($A$, $B$) were sampled and compared.\\
Following the simulation results presented in Section \ref{sec:asympNormality}, we used the delta-method to derive a corrected version $\widetilde{W}_s$ of the test statistic $W_s$ defined in Equation \eqref{eq:statComparison}. Similarly to Section \ref{sec:asympNormality}, the performances of the uncorrected test statistic $W_s$ become similar to these of the corrected version $\widetilde{W}_s$ for large networks (results not shown).

\paragraph{Illustration.} Again, to illustrate the effect of the proposed correction, we provide in Table \ref{tab:testsCompCor} the values of corrected  statistics $\widetilde{W}_s$ testing $H_0 = \{\gt^A = \gt^B\}$, network $A$ being plant-pollinator and network $B$ being seed dispersal. These results can be compared with Table \ref{tab:testsComp}: The correction yields in (moderately) higher absolute values, suggesting a gain of power.

\begin{table}[ht]
  \begin{center}
    \begin{tabular}{clllll}
      $s$ & \quad 5 & \quad 6 & \quad 10 & \quad 15 & \quad 16 \\
      \hline
      $\widetilde{W}_s$ & -2.71 & -1.90 & -1.76 & -1.34 & -0.96 
    \end{tabular}
    \caption{{Corrected test statistics $\widetilde{W}_s$ for $H_0 = \{\gt^A = \gt^B\}$ for the same motifs as in Table~\ref{tab:freqs} and same networks as in Table \ref{tab:testsComp}.}} \label{tab:testsCompCor}
  \end{center}
\end{table}

\begin{figure}[ht]
    \centering
    \begin{tabular}{cccc}
    $m=n=50$ & $m=n=100$ & $m=n=200$ & $m=n=500$ \\
    \hline
    \includegraphics[width=.25\textwidth, trim=10 10 10 10]{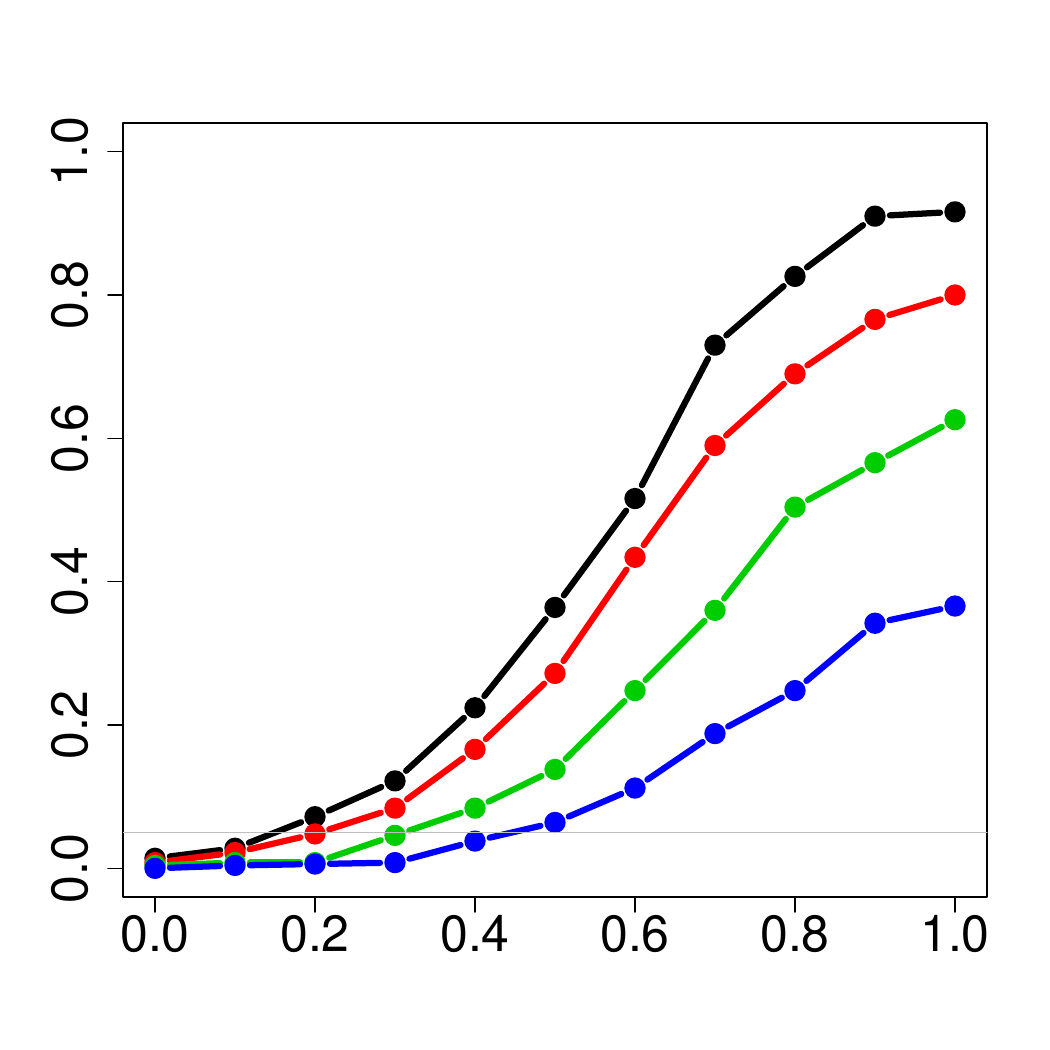} &
    \includegraphics[width=.25\textwidth, trim=10 10 10 10]{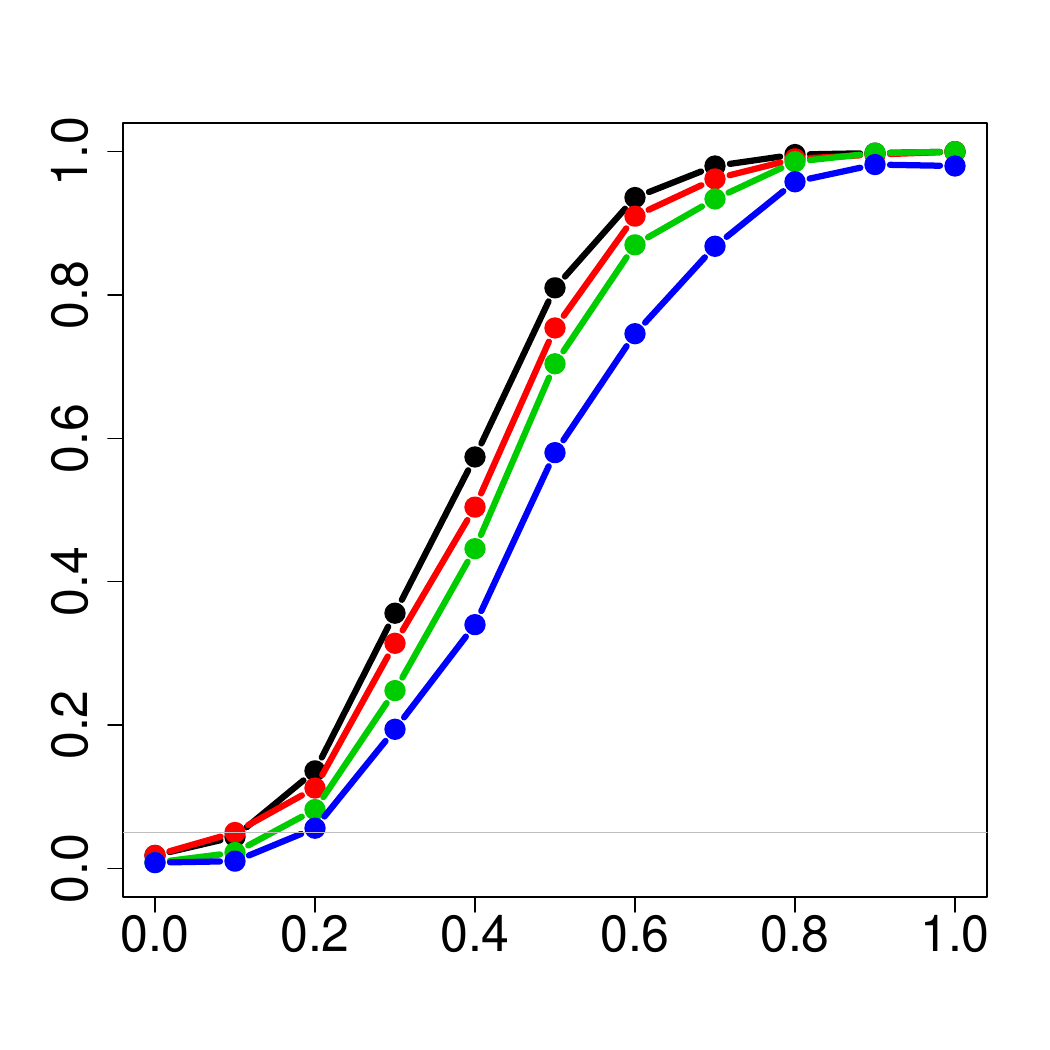} &
    \includegraphics[width=.25\textwidth, trim=10 10 10 10]{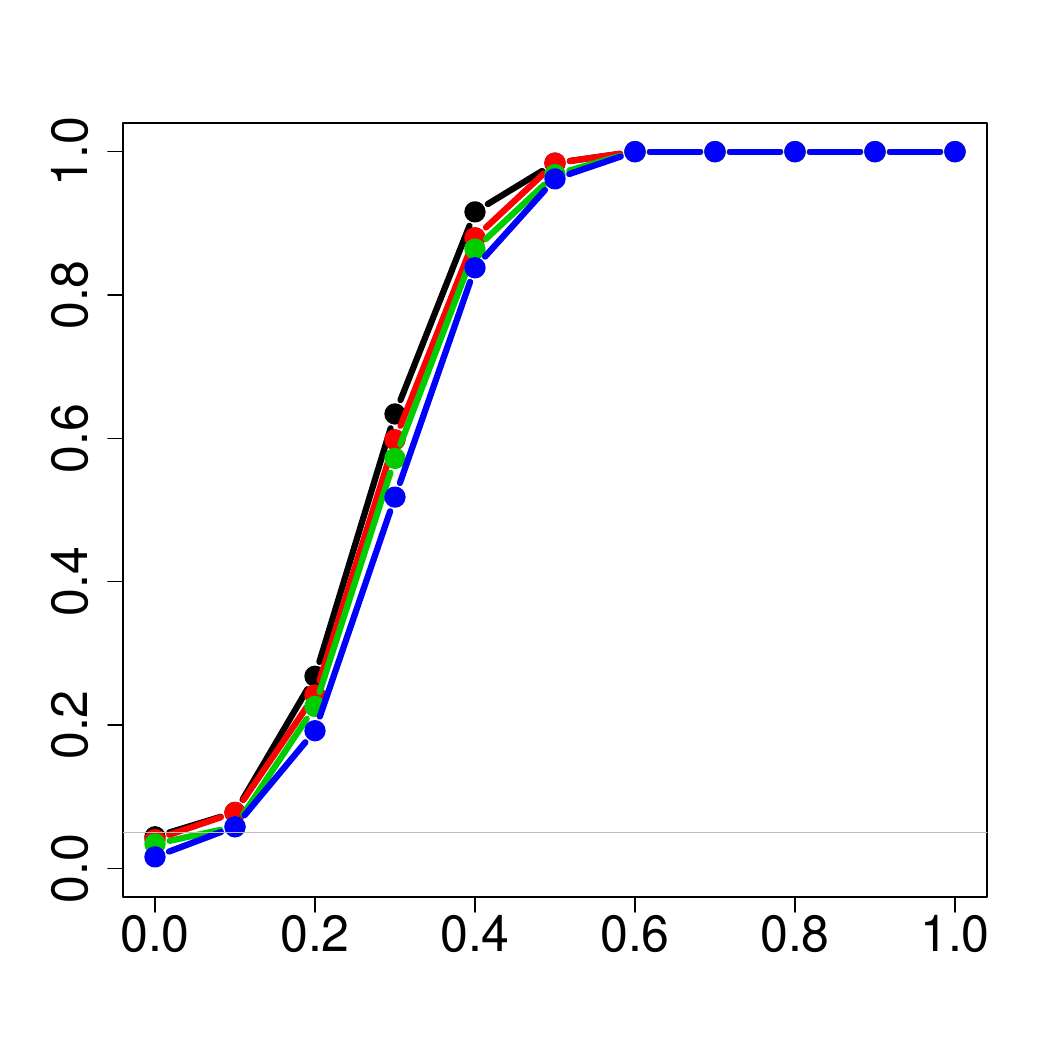} &
    \includegraphics[width=.25\textwidth, trim=10 10 10 10]{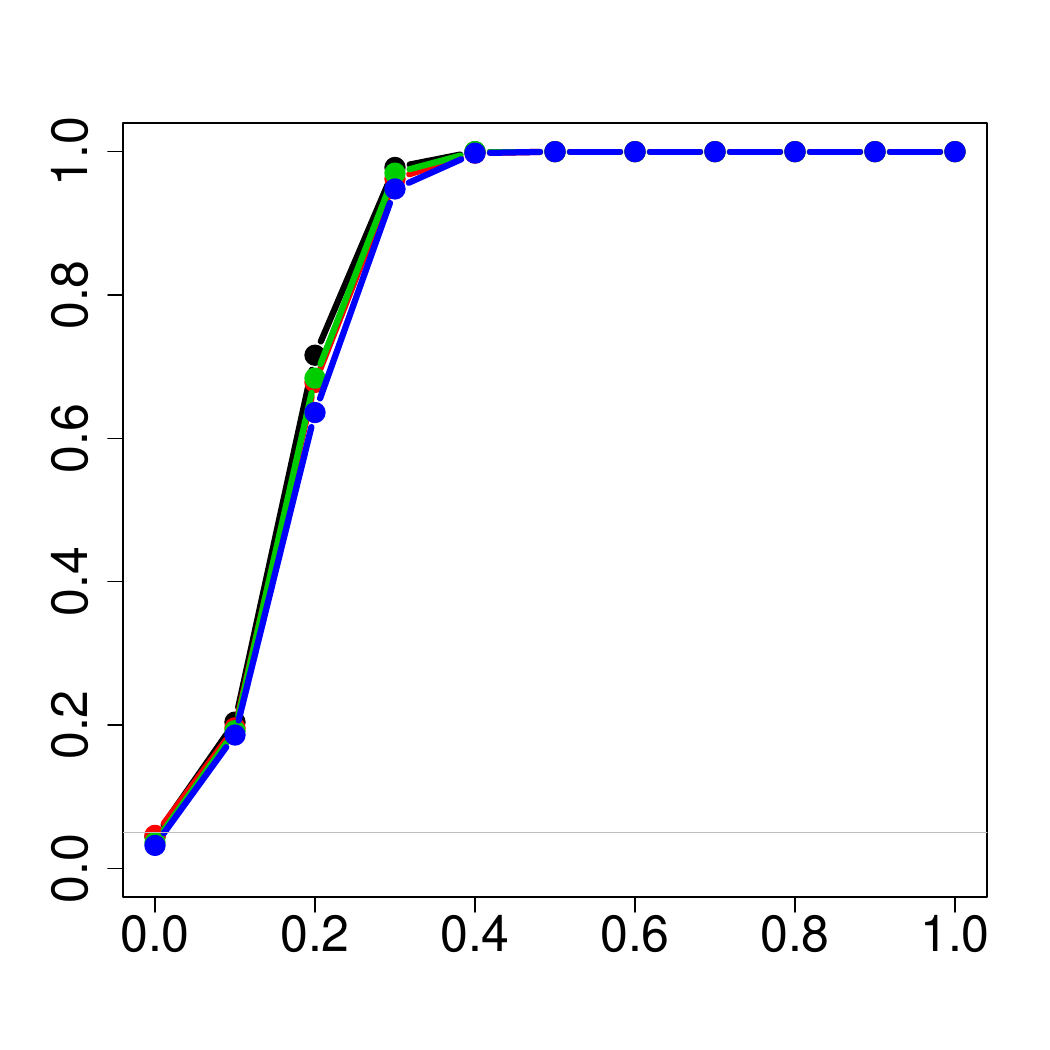} \\
    \includegraphics[width=.25\textwidth, trim=10 10 10 10]{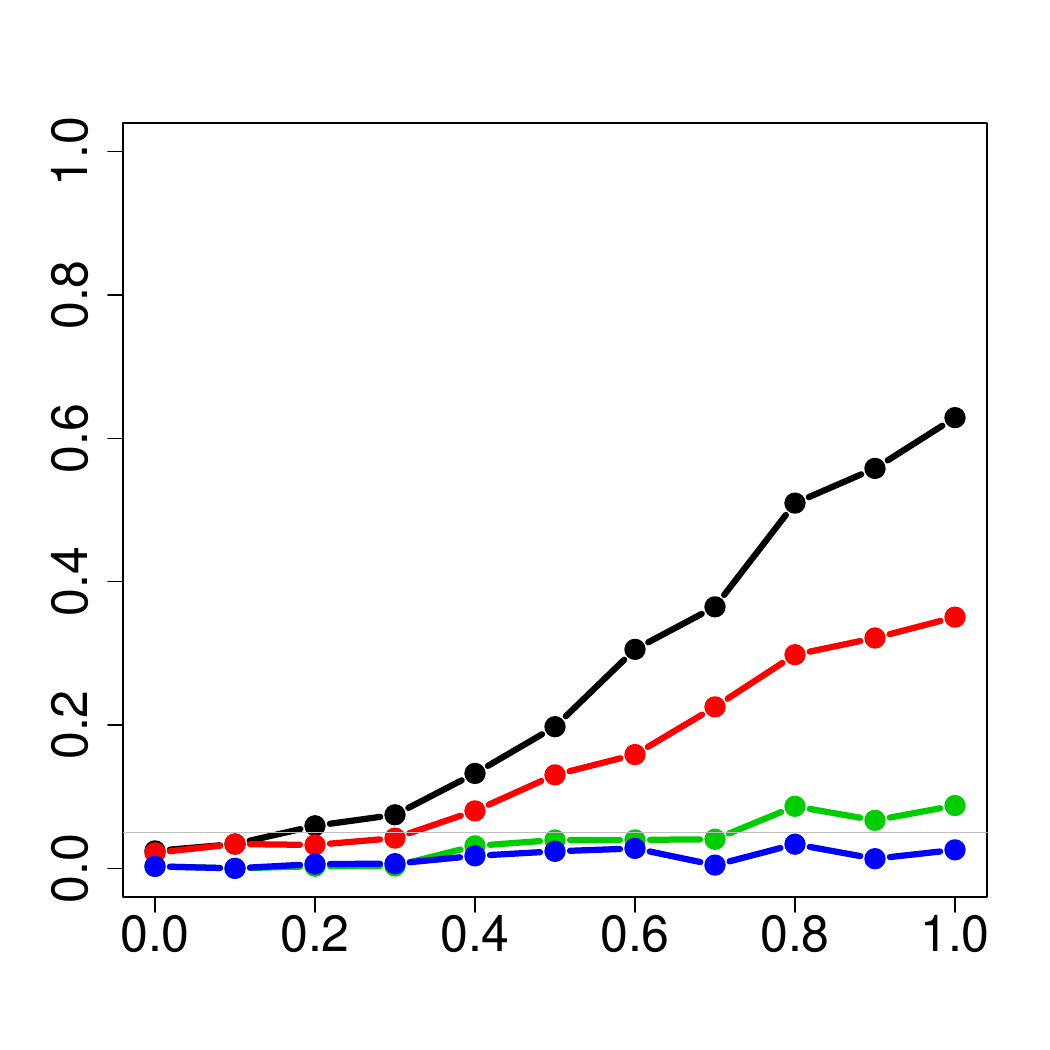} &
    \includegraphics[width=.25\textwidth, trim=10 10 10 10]{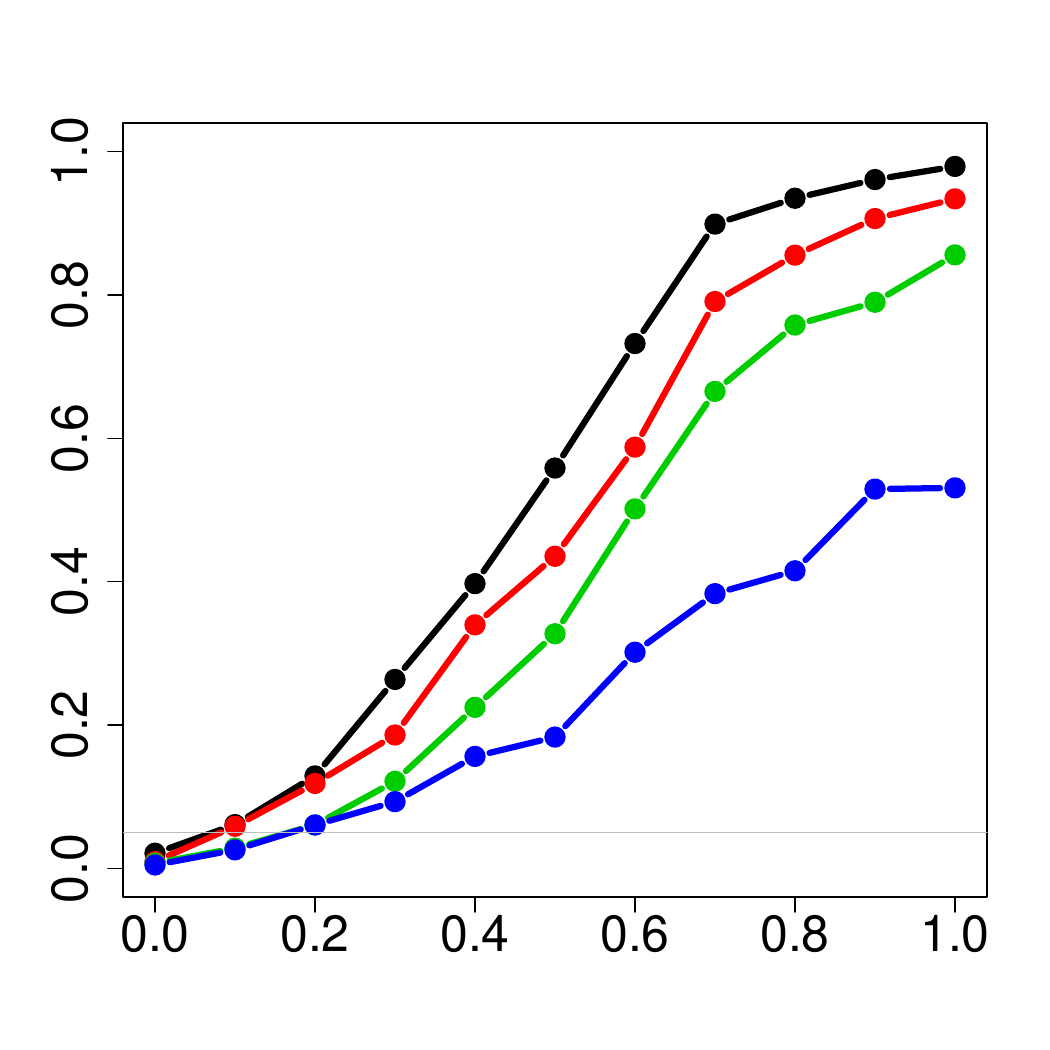} &
    \includegraphics[width=.25\textwidth, trim=10 10 10 10]{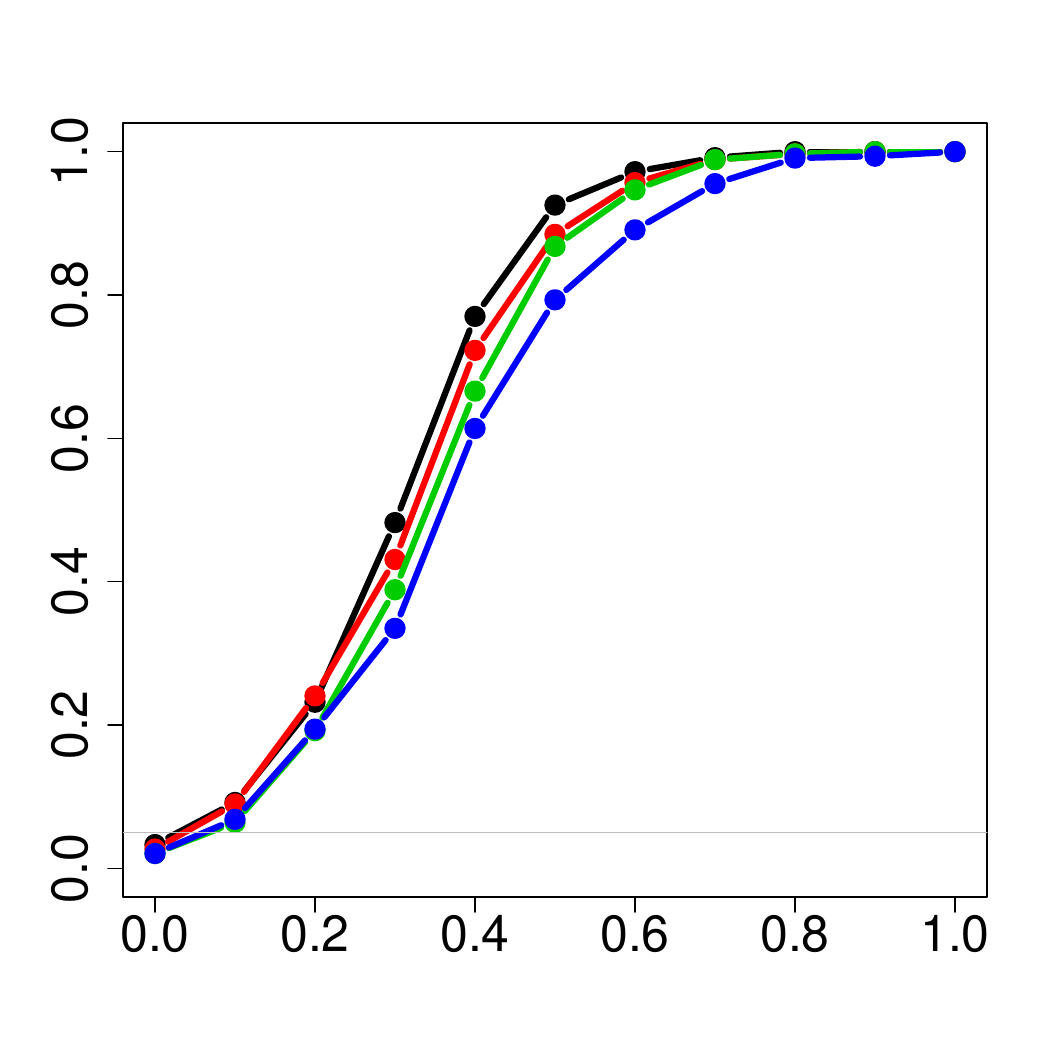} &
    \includegraphics[width=.25\textwidth, trim=10 10 10 10]{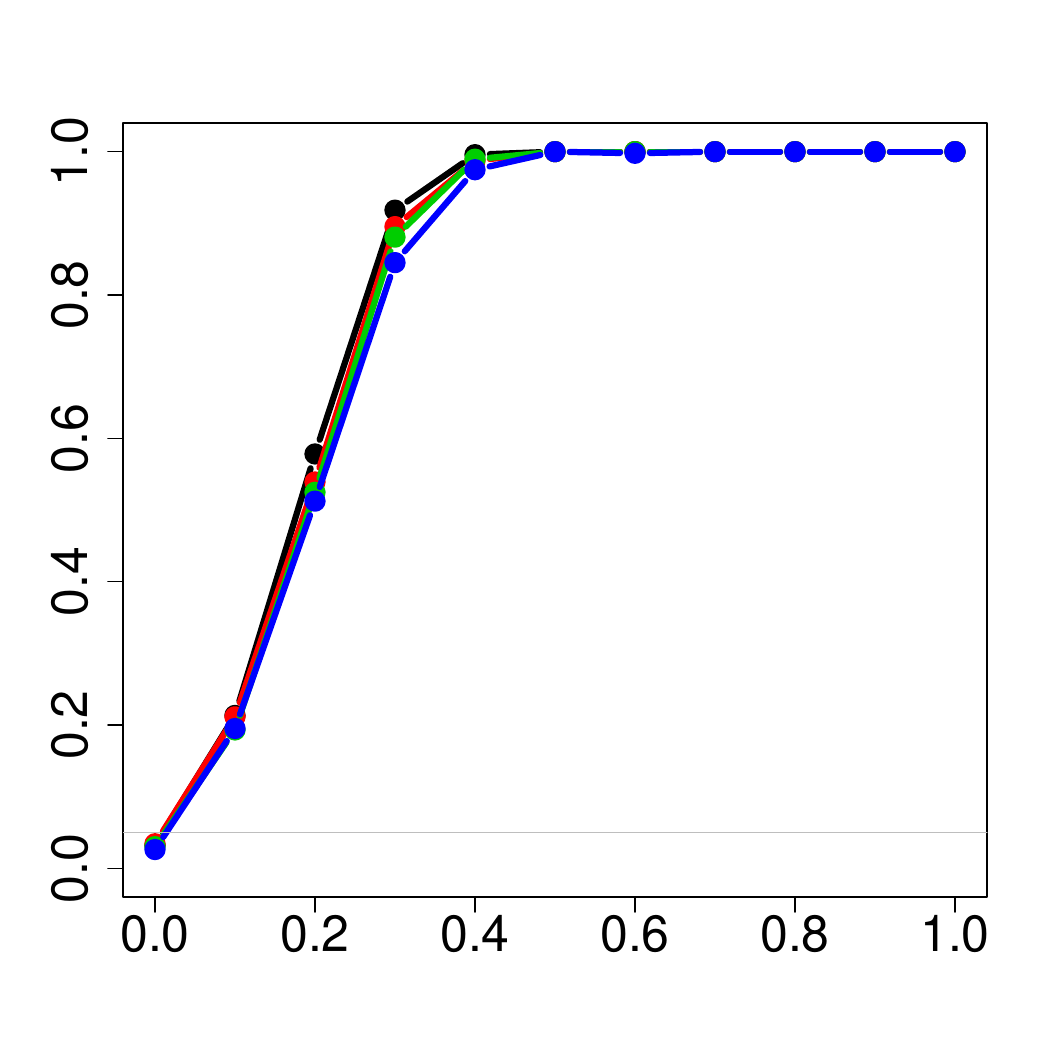}
    \end{tabular}
    \caption{Empirical power of the network comparison test for $H_0 = \{\gt^A = \gt^B\}$, averaged over $S=500$ simulations. Top: scenario I (easy); bottom: scenario II (hard). From left to right: $m = n = 50, 100, 200, 500$. Color = motif: same legend as Figure~\ref{fig:powerGOF}.}
    \label{fig:powerComparison}
\end{figure}

\paragraph{Results.}
The results are displayed in Figure \ref{fig:powerComparison}. We only present the results for $m_A = n_A = m_B = n_B$ ranging for 50 to 500. Moreover, as in the previous section, we only consider motifs 5, 6, 10 and 15. \\
As expected, the test becomes more powerful when the networks dimensions increase. More interestingly, for small networks, the smaller motifs (5 and 6, with size 4) turn out to yield a higher power. The difference vanishes when the dimensions increase. \\
These conclusions hold under the two scenarios, which shows that the proposed test statistic does accommodate for departures that may exist between two networks, not being the departure under study (scenario II 'hard'). Still, the power is always better under scenario I: obviously, the test performs better when focusing on the only difference that actually exists (scenario I 'easy').

\section{Proofs}  \label{sec:proof}

\subsection{Definitions and technical lemmas}

In this section, we introduce notations and useful technical lemmas for establishing proofs of Proposition \ref{Prop:Ls} in Section \ref{sec:prop}, Lemma \ref{Lemma:Cs} in Section  \ref{sec:lemma:cs} and Lemma \ref{Lemma:Var} in Section \ref{sec:lemma:var}.

\subsubsection{Definitions} \label{Sec:Notations}
Let remind that we consider a bipartite graph $\Gcal = (\Vcal, \Ecal)$ with $N$ nodes. The set of nodes is $\Vcal = (\Vcal^t, \Vcal^b)$, where $\Vcal^t=\llbracket{1,m} \rrbracket$ (resp. $\Vcal^b\llbracket{1,n} \rrbracket$) stands for the set of top (resp. bottom) nodes, and the set of edges is $\Ecal \subset \Vcal^t \times \Vcal^b$, meaning than an edge can only connect a top node with a bottom node. The total number of nodes is therefore $N = n+m$. We denote by $G$ the corresponding $m \times n$ incidence matrix where the entry $G_{ij}$ of $G$ is 1 if $(i,j)\in \Ecal$, and 0
otherwise. 

Let consider now a collection of bipartite graphs $(\Gcal_{\ell})_{\ell\in\llbracket 1,N \rrbracket} =(\Vcal_{\ell}, \Ecal_{\ell})$ with $\ell$ nodes. In the following, we introduce notations for subsets of interest and a filtration we will use to construct differences of martingales involving motif counts.

\paragraph{Subsets definitions.}
Let introduce the following subsets definitions:
\begin{itemize}
\item $\Vcal_{\ell}=\{(k_1,\ldots,k_{\ell})\subset\Vcal^t\cup\Vcal^b\mbox{ with at least one top node and one bottom node}\}$, $\ell\in\llbracket{2, N} \rrbracket$, it is the set of nodes of  $\Gcal_{\ell}$ meaning the $\ell$ selected nodes among $\Gcal$, and $k_{\ell}$ denotes the $\ell$-th and last selected one; we will use $k_{\ell}$ several times hereafter;
\item $V_{\ell}^t=\Vcal_{\ell}\cap\Vcal^t$ and $V_{\ell}^b=\Vcal_{\ell}\cap\Vcal^b$, these are the sets of top and bottom nodes in $\Vcal_{\ell}$;
\item $\Pcal_{s,\ell}=\left\{(i_1,\ldots,i_{\pt_s})\subset V_{\ell}^{t} \right\}\times\left\{(j_1,\ldots,j_{\pb_s})\subset V_{\ell}^{b} \right\}$, $\ell\in\llbracket{\pt_s+\pb_s, N} \rrbracket$, it is the positions set of motif $s$ in $\Gcal_{\ell}$; 
\item $T_{\ell}=\{k_{\ell}\in \Vcal^t
\}$ is an event;
\item $\Qcal_{s,\ell}=\left\{\begin{array}{lcl}
\{\Pcal_{s,\ell -1}\setminus i_{\pt_s}\}\cup\{i_{\pt_s}=k_{\ell}\} & & \mbox{ if } T_{\ell},\\
\{\Pcal_{s,\ell -1}\setminus j_{\pb_s}\}\cup\{j_{\pb_s}=k_{\ell}\}& &\mbox{ otherwise},\end{array}\right.$ \label{Qsl}

it is the positions set of motif $s$ in $\Gcal_{\ell}$ 
with the particularity that $k_{\ell}$ the last node added to $\Vcal_{\ell}$ is part of motif $s$.
\end{itemize}

\paragraph{Filtration.}
The filtration $(\Fcal_{\ell})_{\ell\in\llbracket{2,N} \rrbracket}$ is defined by the $\sigma$-algebra $\Fcal_{\ell}=\sigma\left(\Gcal_{\ell}\right)$.

\subsubsection{Technical lemmas} \label{Sec:tech:lem}

We present here three lemmas which are key arguments in the proofs of Proposition \ref{Prop:Ls}, 
 Lemma \ref{Lemma:Cs} and Lemma \ref{Lemma:Var}.

The following lemma gives the order of magnitude of the variance of a count. Before, its statement let give the order of magnitude of the expected count of a motif $s$ with $\pt_s$ top nodes and $\pb_s$ bottom nodes. It writes $
\Esp(N_s) = c_s \phi_s $, with

\begin{eqnarray} 
c_s &=& \Theta(m^{\pt_s} n^{\pb_s})\quad \mbox{(normalizing coefficient specific to }s)\label{order:cs} \\
\rho &=&\Theta(m^{-a} n^{-b}), \quad\mathrm{ with }\quad a,b>0\quad\mbox{(graph density)}\label{order:rho}\\
\phi_s &=& \Theta(\rho^{d_+^s} )=\Theta(m^{-ad_+^s} n^{-bd_+^s})\quad \mbox{(expected frequency of }s)\label{order:phis},
\end{eqnarray}
where $\dt_+^s$ stands for the total number of edges in $s$ and $c_s$ being defined in \eqref{eq:cs}.

\begin{lemma}\label{Lemma:var:count}
We have,
\begin{align*}
 \Var(N_s) = \Theta\left(\max(m^{2\pt_s-2ad_+^{s}-1} n^{2\pb_s-2bd_+^{s}}, m^{2\pt_s-2ad_+^{s}} n^{2\pb_s-2bd_+^{s}-1}, m^{2\pt_s-ad_+^{s}-1} n^{2\pb_s-bd_+^{s}-1})\right).
\end{align*}
\end{lemma}
 
\proofbegin
Let observe that, for $\alpha,\beta\in\Pcal_{s,N}$,
$$
N_s^2 = \sum_{\alpha} Y_s(\alpha) 
+ \sum_{\alpha \cap \beta \neq \emptyset} Y_s(\alpha) Y_s(\beta)
+ \sum_{\alpha \cap \beta = \emptyset} Y_s(\alpha) Y_s(\beta).
$$
Thus, a general form for the variance is the following:
\begin{equation}\label{eq:Var:Ns}
\Var(N_s) = \Esp(N_s) + \sum_{t \in \Scal_2(s)} \Esp(N_t) + \left(|\Ocal_{s}| - c_s^2\right) \phi_s^2,
\end{equation}
where $\Ocal_{s}=\left\{\alpha,\beta\in\Pcal_{s,N}: \alpha\cap\beta=\emptyset \right\}$ and $\Scal_2(s)$ denotes the set of supermotifs of $s$ which are formed by two overlapping occurrences of $s$.

Let evaluate the orders of the three added terms of assertion \eqref{eq:Var:Ns}.
Considering that $\rho =\Theta(m^{-a} n^{-b})$, the first term of \eqref{eq:Var:Ns} is $\Theta(m^{\pt_s-ad_+^{s}} n^{\pb_s-bd_+^{s}})$.
Then denoting $(a)_b = a (a-1) \dots (a-b)$, we see that
\begin{align}\label{order:calculus}
 |\Ocal_{s}| - c_s^2
 & = \frac{(m_{2\pt_s-1})}{(\pt_s!)^2} \frac{(n)_{2\pb_s-1}}{(\pb_s!)^2} - \frac{(m_{\pt_s-1})^2}{(\pt_s!)^2} \frac{(n)^2_{\pb_s-1}}{(\pb_s!)^2} \\
 & =\Theta\left( \frac{(-1)^{2\pt_s-1} \pt_s^2 m^{2\pt_s-1} n^{2\pb_s} + (-1)^{2\pb_s-1} \pb_s^2 m^{2\pt_s} n^{2\pb_s-1}}{(\pt_s!)^2(\pb_s!)^2} \right) \nonumber\\
 & =\Theta\left( \max(m^{2\pt_s-1} n^{2\pb_s}, m^{2\pt_s} n^{2\pb_s-1}) \right).\nonumber
\end{align}
Thus the third term is $\Theta\left( \max(m^{2\pt_s-2ad_+^{s}-1} n^{2\pb_s-2bd_+^{s}}, m^{2\pt_s-2ad_+^{s}} n^{2\pb_s-2bd_+^{s}-1}) \right).$

Let focus now on the second term. When $t \in \Scal_k(s)$, it can result of an overlap of ($i$) only top nodes, ($ii$) only bottom nodes,   or ($iii$) both. For each case we have
\begin{itemize}
 \item[($i$)] $\pt_t < 2 \pt_s, \pb_t = 2\pb_s, d_+^{t} = 2 d_+^{s}$ so $\Esp N_t = O(m^{2\pt_s-1} n^{2\pb_s} \rho^{2d_+^{s}})  = O(m^{2\pt_s-2ad_+^{s}-1} n^{2\pb_s-2bd_+^{s}})$;
 \item[($ii$)] $\pt_t = 2 \pt_s, \pb_t < 2\pb_s, d_+^{t} = 2 d_+^{s}$ so $\Esp N_t = O(m^{2\pt_s} n^{2\pb_s-1} \rho^{2d_+^{s}})  = O(m^{2\pt_s-2ad_+^{s}} n^{2\pb_s-2bd_+^{s}-1})$;
 \item[($iii$)] $\pt_t < 2 \pt_s, \pb_t < 2\pb_s, d_+^{s} < d_+^{t} < 2 d_+^{s}$ so $\Esp N_t = O(m^{2\pt_s-ad_+^{s}-1} n^{2\pb_s-bd_+^{s}-1})$.
\end{itemize}

Combining the orders of the three terms of assertion \eqref{eq:Var:Ns}, we get that the order of magnitude of the variance of a count is
\begin{align*}
 \Var(N_s) = \Theta\left(\max(m^{2\pt_s-2ad_+^{s}-1} n^{2\pb_s-2bd_+^{s}}, m^{2\pt_s-2ad_+^{s}} n^{2\pb_s-2bd_+^{s}-1}, m^{2\pt_s-ad_+^{s}-1} n^{2\pb_s-bd_+^{s}-1})\right).
\end{align*}
\proofend 

The last argument of proof of Proposition \ref{Prop:Ls}, Lemma \ref{Lemma:Ms} and Lemma \ref{Lemma:Cs} relies on the following result.
\begin{lemma}\label{Lemma:Var:deconditionning}
We have, as $m\sim n\to\infty$,
$$\Var(N_s| U,V)/\Var(N_s)\rightarrow 1\mbox{ in probability}.$$
\end{lemma}

\proofbegin 
First let us write that  
\begin{align*}
&\Esp(N_s|U,V)= \sum_{\alpha\in\Pcal_s}\mathbb{P}\left(Y_s(\alpha)=1 |  U_{\alpha^t},V_{\alpha^b}\right)\\
& \Esp(N_s^2|U,V)=\sum_{\alpha,\beta \in\Pcal_s}\mathbb{P}\left(Y_s(\alpha)Y_s(\beta)=1  |  U_{\alpha^t},V_{\alpha^b},U_{\beta^t},V_{\beta^b}\right).
\end{align*}
The proof relies on showing the convergence in probability of the two above expectations towards $\sum_{\alpha\in\Pcal_s}\mathbb{P}\left(Y_s(\alpha)=1\right)$ and $\sum_{\alpha,\beta \in\Pcal_s}\mathbb{P}\left(Y_s(\alpha)Y_s(\beta)=1\right)$, respectively.
Let us now introduce the equivalence relation $\textgoth{R}_s$ and the set $R_s$ defined as follows:
$$\textgoth{R}_s:(\sigma_t,\sigma_b)\sim (\tilde{\sigma_t},\tilde{\sigma_b}) \Leftrightarrow A_{\sigma_t,\sigma_b}^s=A_{\tilde{\sigma_t},\tilde{\sigma_b}}^s\ \mbox{ and }R_s=\left( \sigma\left(\llbracket 1,\pt_s\rrbracket\right)\otimes\sigma\left(\llbracket 1,\pb_s\rrbracket\right) \right)/\textgoth{R}_s.$$
Then, we can exhibit the two following quantities which are two-samples U-Statistics (see Section 12.2, p.165 in \cite{Vaart00}):
\begin{eqnarray*}
\frac{r_s}{c_s}\sum_{\alpha\in\mathcal{L}_s} k_1\left(U_{\alpha^t},V_{\alpha^b}\right)\quad \mbox{ and }\quad \frac{r_s}{c_s}\sum_{\alpha\in\mathcal{L}_s} k_2\left(U_{\alpha^t},V_{\alpha^b},U_{\beta^t},V_{\beta^b}\right),
\end{eqnarray*}
with $r_S$ and $c_s$ being defined in \eqref{eq:rs}, \eqref{eq:cs}, respectively,  $\mathcal{L}_s$ denoting the location relative to a given position for motif $s$ and
where 
\begin{align*}
& k_1\left(U_{\alpha^t},V_{\alpha^b}\right)=\sum_{\sigma\in R_s} 
\mathbb{P}\left(
Y_s(\sigma(\alpha))=1 |  U_{\sigma_t(\alpha^t)},V_{\sigma_b(\alpha^b)}
\right)\\
& k_2\left(U_{\alpha^t},V_{\alpha^b},U_{\beta^t},V_{\beta^b}\right)=
\sum_{(\sigma_{\alpha},\sigma_{\beta}) \in R_s} 
\mathbb{P}\left(
Y_s(\sigma_{\alpha}(\alpha))Y_s(\sigma_{\beta}(\beta))=1  |  U_{\sigma_{\alpha}^t(\alpha^t)},V_{\sigma_{\alpha}^b(\alpha^b)},U_{\sigma_{\beta}^t(\beta^t)},V_{\sigma_{\beta}^b(\beta^b)}
\right),
\end{align*}
with $k_1(\cdotp)$ and $k_2(\cdotp)$ being permutation symmetric kernels in $(U_i)_i$ and $(V_j)_j$ separetely.
We conclude by applying the central limit theorem for two-sample U-Statistics (see Theorem 12.6 in \cite{Vaart00}) which holds under the assumption that the kernel of the U-statistic has a finite moment of order two. Here, as it concerns probabilities this assumption is obviously fulfilled.
\proofend


In proofs of Lemma \ref{Lemma:Var}, Lemma \ref{Lemma:Rs} and Lemma \ref{Lemma:Ms}, we need to know the cardinal order of the sets $\Qcal_{s,\ell}^{\otimes k}\setminus\Ocal_{s,\ell}^{(k)}, k=2,4$ which contains only dependent $k$-uplets of positions of motif $s$ on the event $T_{\ell}$ for which the last node added to $\Vcal_{\ell}$ is a top node. Recall that
$\Qcal_{s,\ell}$ is the positions set of motif $s$ in the subgraph of $\Gcal$ with nodes in $\Vcal_{\ell}$ and the particularity that $ k_{\ell}$ the last node added to $\Vcal_{\ell}$ is part of motif $s$. The definition of the other set of interest is the following: 

$$\Ocal_{s,\ell}^{(k)}=\left\{\alpha_1,\ldots\alpha_k\in\Qcal_{s,\ell}:(\alpha_1^t\setminus k_{\ell})\times\alpha_1^b\cap\ldots\cap(\alpha_k^b\setminus k_{\ell})\times\alpha_k^b=\emptyset \right\}.$$

%

\begin{lemma}\label{Lemma:Cardinal}
We have, on $T_{\ell}$,
$$|\Qcal_{s,\ell}|^k-|\Ocal_{s,\ell}^{(k)}|=\Theta\left(\ell_t^{k(\pt_s-1)-1 }\ell_b^{k\pb_s}\right),
$$
with $\ell_t$ and $\ell_b$ denoting respectively top and bottom nodes in $\Vcal_{\ell}$.
\end{lemma}

\proofbegin 
Let observe that
\begin{align*}
& |\Qcal_{s,\ell}|=\binom{\ell_t-1}{\pt_s-1}\binom{\ell_b}{\pb_s}, \\
& |\Ocal_{s,\ell}^{(k)}|=\binom{\ell_t-1}{\pt_s-1}^k\binom{\ell_b}{\pb_s\ldots \pb_s \ \ell_b-k\pb_s}+\binom{\ell_t}{\pt_s\ldots \pt_s \ \ell_t-k\pt_s}\binom{\ell_b}{\pb_s}^k\\
&\qquad\qquad +\binom{\ell_t}{\pt_s\ldots \pt_s \ \ell_t-k\pt_s}\binom{\ell_b}{\pb_s\ldots \pb_s \ \ell_b-k\pb_s}.
\end{align*}
The leader term of order $\Theta\left(\ell_t^{k(\pt_s-1) }\ell_b^{k\pb_s}\right)$ obviously vanishes and imply the lost of one order (the calculation omitted here are simply based on the same arguments as in \eqref{order:calculus}).
\proofend

\subsection{Proof of Proposition \ref{Prop:Ls}} \label{sec:prop}

For establishing the proof of  Proposition \ref{Prop:Ls}, we first consider a decomposition of $L_s=F_s - \phibar_s$ in Section \ref{subsec:Ls} , then we focus on the reminder term of this decomposition in Lemma \ref{Lemma:Rs} and finally show the asymptotic normality of the leading term in Lemma \ref{Lemma:Ms}.

\subsubsection{Decomposition of $L_s$}\label{subsec:Ls}
Let use the sets introduced in Section \ref{Sec:Notations} to express $L_s$ as follows:
\begin{eqnarray*}
L_s(U,V)=F_s - \phi_s(U,V)&=& \frac{1}{c_s}\sum_{\alpha=(\alpha^t,\alpha^b)\in\Pcal_{s,N}}\{Y_s(\alpha) - \phibar_s(U_{\alpha^t}, V_{\alpha^b})\}\\
&=&\frac{1}{c_s}\sum_{\ell=1}^N \sum_{\alpha\in \Qcal_{s,\ell}}\{Y_s(\alpha) - \phibar_s(U_{\alpha^t}, V_{\alpha^b})\},
\end{eqnarray*}
with the random variables $U,V$ of the \BEDD model \eqref{eq:BEDD}.
Then let decompose 
$L_s$ as the sum of two expressions, the first one corresponding to a martingale difference sequence relative to the filtration $(\Fcal_{\ell})_{\ell\in\llbracket{2,N} \rrbracket}$, the second one being a term of rest:
\begin{eqnarray*}
L_s(U,V)&:=&M_s(U,V)+R_s(U,V),
\end{eqnarray*}
where
\begin{align*}
& M_s(U,V)=\frac{1}{c_s}\sum_{\ell=1}^N\sum_{\alpha\in \Qcal_{s,\ell}}\{Y_s(\alpha) - \Esp(Y_s(\alpha)|\Fcal_{\ell-1}; U,V)\}\\
& R_s(U,V)= \frac{1}{c_s}\sum_{\ell=1}^N\sum_{\alpha\in \Qcal_{s,\ell}}\{\Esp(Y_s(\alpha)|\Fcal_{\ell-1}; U,V) - \phi_s(U_{\alpha^t}, V_{\alpha^b})\}. 
\end{align*}
Observe that by construction, $M_{s,\ell}=\sum_{\alpha\in \Qcal_{s,\ell}}\{Y_s(\alpha) - \Esp(Y_s(\alpha)|\Fcal_{\ell-1}; U,V)\}$ is a conditional martingale difference with respect to $(\Fcal_{\ell})_{\ell\in\llbracket{2,N} \rrbracket}$: $$\Esp(M_{s,\ell}(U,V)|\Fcal_{\ell-1}; U,V)=0.$$

\subsubsection{Study of $R_s$}
\begin{lemma}\label{Lemma:Rs}
Under the \BEDD model and condition $a+b<2/d_+^{s}$, 
$$R_s(U,V)/\sqrt{\Var(F_s)}| U,V\to 0\mbox{ a.s. as }  m\sim n\to\infty,$$
where $R_s(U,V)=\frac{1}{c_s}\sum_{\ell=1}^N\sum_{\alpha\in \Qcal_{s,\ell}}\{\Esp(Y_s(\alpha)|\Fcal_{\ell-1}; U,V) - \phibar_s(U_{\alpha^t}, V_{\alpha^b})\} $.
\end{lemma}

\proofbegin
The proof consists in showing the two following assertions:
\begin{itemize}
\item[(A1)] $\Esp\left(R_s(U,V)/\sqrt{\Var(F_s)}| U,V\right)=0$;
\item[(A2)] $\Var\left(c_sR_s(U,V)/\sqrt{\Var(N_s)}| U,V\right)\rightarrow 0$ almost surely as $n$ tends to infinity  under condition $a+b<2/d_+^{s}$.
\end{itemize}
Let show assertion (A1): 
\begin{eqnarray*}
\Esp\left(R_s(U,V)| U,V\right)
&=&\Esp\left(\frac{1}{c_s}\sum_{\ell=1}^N\sum_{\alpha\in \Qcal_{s,\ell}}\{\Esp(Y_s(\alpha)|\Fcal_{\ell-1}; U,V) - \phibar_s(U_{\alpha^t}, V_{\alpha^b})\} |U,V\right)\\
&=&\Esp\left(\frac{1}{c_s}\sum_{\ell=1}^N\sum_{\alpha\in \Qcal_{s,\ell}}\Esp(Y_s(\alpha)|\Fcal_{\ell-1}; U,V)\right)
-\frac{1}{c_s}\sum_{\ell=1}^N\sum_{\alpha\in \Qcal_{s,\ell}}\phibar_s(U_{\alpha^t}, V_{\alpha^b})\\
&=&\frac{1}{c_s}\sum_{\alpha=(\alpha^t,\alpha^b)\in\Pcal_{s,N}}\Esp(Y_s(\alpha)|U,V)
-\frac{1}{c_s}\sum_{\alpha=(\alpha^t,\alpha^b)\in\Pcal_{s,N}}\phibar_s(U_{\alpha^t}, V_{\alpha^b})=0.
\end{eqnarray*}

Let focus now on assertion (A2).
Let first observe that,
\begin{align*}
&\Var\left(R_s(U,V)| U,V\right)\\
&=\Var\left(\frac{1}{c_s}\sum_{\ell=1}^N\sum_{\alpha\in \Qcal_{s,\ell}}\{\Esp(Y_s(\alpha)|\Fcal_{\ell-1}; U,V) - \phibar_s(U_{\alpha^t}, V_{\alpha^b})\}|U,V\right)\\
&=\sum_{\ell=1}^N\Var\left(\frac{1}{c_s}\sum_{\alpha\in \Qcal_{s,\ell}}\Esp(Y_s(\alpha)|\Fcal_{\ell-1}; U,V)|U,V\right),
\end{align*}
by independance of successive  choices of $\mathcal{G}_{\ell}$. Using definition (\ref{eq:Yalpha}) of the indicator motif, we see that
\begin{align*}
&\sum_{\ell=1}^N\Var\left(\frac{1}{c_s}\sum_{\alpha\in \Qcal_{s,\ell}}\Esp(Y_s(\alpha)|\Fcal_{\ell-1}; U,V)|U,V\right)\\
&=\sum_{\ell=1}^N\Var\left(\frac{1}{c_s}\sum_{\alpha\in \Qcal_{s,\ell}}\Esp\left(\prod_{i\in \alpha^t, j\in\alpha^b} G_{ij}^{A^s_{ij}} |\Fcal_{\ell-1}; U,V\right)|U,V\right).
\end{align*}
Then according to measurability with respect to $\Fcal_{\ell-1}$ and the position (top or bottom) of $k_{\ell}$ the last selected node, we get
\begin{align*}
&\sum_{\ell=1}^N\Var\left(\frac{1}{c_s}\sum_{\alpha\in \Qcal_{s,\ell}}\Esp\left(\prod_{i\in \alpha^t, j\in\alpha^b} G_{ij}^{A^s_{ij}} |\Fcal_{\ell-1}; U,V\right)|U,V\right)\\
&=\Pr(T_{\ell})\sum_{\ell=1}^N\Var\left(\frac{1}{c_s}\sum_{\alpha\in \Qcal_{s,\ell}}\left(\prod_{j\in\alpha^b}\Esp(G_{ k_{\ell} j}|U,V)^{A^s( k_{\ell} j)}\right)\left(\prod_{i\in \alpha^t\setminus k_{\ell}, j\in\alpha^b} G_{ij}^{A^s_{ij}} \right)|U,V\right)\\
&+(1-\Pr(T_{\ell}))\sum_{\ell=1}^N\Var\left(\frac{1}{c_s}\sum_{\alpha\in \Qcal_{s,\ell}}\left(\prod_{i\in\alpha^t}\Esp(G_{i k_{\ell}}|U,V)^{A^s(i k_{\ell})}\right)\left(\prod_{i\in\alpha^t, j\in \alpha^b\setminus k_{\ell}} G_{ij}^{A^s_{ij}} \right)|U,V\right),
\end{align*}
and using the usual notation of the conditional expectation of $G_{ij}$'s, we have
\begin{align*}
&\Var\left(R_s(U,V)| U,V\right)\\
&=\Pr(T_{\ell})\sum_{\ell=1}^N\Var\left(\frac{1}{c_s}\sum_{\alpha\in \Qcal_{s,\ell}}\left(\prod_{j\in\alpha^b}\phi_1(U_{k_{\ell}}, V_j)^{A^s( k_{\ell} j)}\right)\left(\prod_{i\in \alpha^t\setminus k_{\ell}, j\in\alpha^b} G_{ij}^{A^s_{ij}} \right)|U,V\right)\\
&+(1-\Pr(T_{\ell}))\sum_{\ell=1}^N\Var\left(\frac{1}{c_s}\sum_{\alpha\in \Qcal_{s,\ell}}\left(\prod_{i\in\alpha^t}\phi_1(U_i, V_{k_{\ell}})^{A^s( i k_{\ell} )}\right)\left(\prod_{i\in\alpha^t, j\in \alpha^b\setminus k_{\ell}} G_{ij}^{A^s_{ij}} \right)|U,V\right).
\end{align*}

Then, considering the fact that $\Var(\sum_ia_iX_i)\leq\left(\sum_i a_i\sqrt{\Var(X_i)}\right)^2$, we get 
\begin{align*}
&\Var\left(R_s(U,V)| U,V\right)\\
&\leq 2\sum_{\ell=1}^N\left(\frac{1}{c_s}\sum_{\alpha\in \Qcal_{s,\ell}}\left(\prod_{j\in\alpha^b}\phi_1(U_{k_{\ell}}, V_j)^{A^s( k_{\ell} j)}\right)\sqrt{\Var\left(\prod_{i\in \alpha^t\setminus k_{\ell}, j\in\alpha^b} G_{ij}^{A^s_{ij}} |U,V\right)}\right)^2.
\end{align*}
From now, we will work on the set 
$\Qcal_{s,\ell}^{\otimes 2}\setminus\Ocal_{s,\ell}^{(2)}$ which contains only dependent pairs of positions. 
It follows from the Bernoulli conditional distribution of $G_{ij}$ combined with the fact that $a\sqrt{b}<\sqrt{ab}$ when $a<1$, that
\begin{align*}
&\Var\left(R_s(U,V)| U,V\right)\\
&\leq 2
\sum_{\ell=1}^N\left(\frac{1}{c_s}\sum_{\alpha\in \Qcal_{s,\ell}^{\otimes 2}\setminus\Ocal_{s,\ell}^{(2)}}\sqrt{\prod_{j\in\alpha^b}\phi_1(U_{k_{\ell}}, V_j)^{A_s( k_{\ell} j)}\prod_{i\in \alpha^t\setminus k_{\ell}, j\in\alpha^b}\phi_1(U_i, V_j)^{A_s(ij)}}\right)^2\\
&\leq  2\sum_{\ell=1}^N\frac{1}{c_s^2}\left(\sum_{\alpha\in \Qcal_{s,\ell}^{\otimes 2}\setminus\Ocal_{s,\ell}^{(2)}}\sqrt{\phibar_s(U_{\alpha^t}, V_{\alpha^b})}\right)^2\\
&\leq\frac{2}{c_s^2}\sum_{\ell=1}^N \left(|\Qcal_{s,\ell}|^2-|\Ocal_{s,\ell}^{(2)}|\right)\max_{\alpha\in \Qcal_{s,\ell}^{\otimes 2}\setminus\Ocal_{s,\ell}^{(2)}}\phibar_s(U_{\alpha^t}, V_{\alpha^b}).
\end{align*}

In order to evaluate the right-hand side term of the above inequality, recall that $c_s=\Theta(m^{\pt_s} n^{\pb_s})$ by \eqref{order:cs}, $\phi_s=\Theta\left(m^{-ad_+^{s}}n^{-bd_+^{s}}\right)$ by \eqref{order:phis} and $|\Qcal_{s,\ell}|^2-|\Ocal_{s,\ell}^{(2)}|=\Theta\left(\ell_t^{2\pt_s-3}\ell_b^{2\pb_s}\right)$ by Lemma \ref{Lemma:Cardinal}, $\ell_t$ and $\ell_b$ denoting respectively top and bottom nodes in $\Vcal_{\ell}$. Thus, we get 

\begin{align*}
&\frac{2}{c_s^2}\sum_{\ell=1}^N \left(|\Qcal_{s,\ell}|^2-|\Ocal_{s,\ell}|\right)\max_{\alpha\in\Qcal_{s,\ell}^{\otimes 2}\setminus\Ocal_{s,\ell}^{(2)}
}\phibar_s(U_{\alpha^t}, V_{\alpha^b})\\
& =\Theta\left(m^{-2\pt_s} n^{-2\pb_s}\right)\sum_{\ell=\ell_t+\ell_b=1}^N \Theta\left(\ell_t^{2\pt_s-3} \ell_b^{2\pb_s}\ell_t^{-ad_+^{s}} \ell_b^{-bd_+^{s}}\right)\\
&= \Theta\left(m^{-2\pt_s} n^{-2\pb_s}\right)\sum_{\ell=\ell_t+\ell_b=1}^N \Theta\left(\ell^{2\pt_s+2\pb_s-ad_+^{s}-bd_+^{s}-3}\right)\\
&=\Theta\left(N^{-ad_+^{s}-bd_+^{s}-3}\right).
\end{align*}
By taking the normalization $\sqrt{\Var(F_s)}=\sqrt{\Var(N_s)}/c_s$ which order is $$\Theta\left(\max\left(N^{-2ad_+^{s}-2bd_+^{s}-1},N^{-ad_+^{s}-bd_+^{s}-2}\right)\right)$$ by Lemma \ref{Lemma:Var} and \eqref{order:cs}, we conclude to $\Var\left(\frac{R_s(U,V)}{\sqrt{\Var(F_s)}}| U,V\right)\rightarrow 0$ almost surely as $n$ tends to infinity under condition 
$a+b<2/d_+^{s}$.
\proofend 


\subsubsection{Study of $M_s$}\label{subsec:Ms}
\begin{lemma}\label{Lemma:Ms}
Under the \BEDD model and condition $a+b<2/d_+^{s}$, 
$$M_s(U,V)/\sqrt{\Var(F_s)}| U,V\overset{D}{\longrightarrow} \Ncal\left(0, \frac{\Var(N_s | U,V)}{\Var(N_s)})\right), \mbox{ as } m\sim n\to\infty,$$
where $M_s(U,V)=\frac{1}{c_s}\sum_{\ell=1}^N\sum_{\alpha\in \Qcal_{s,\ell}}\{Y_s(\alpha) - \Esp(Y_s(\alpha)|\Fcal_{\ell-1}; U,V)\}.$
\end{lemma}

\proofbegin
We will apply the following martingale central limit theorem to the conditional martingale difference sequence $M_{s,\ell}(U,V)=\sum_{\alpha\in \Qcal_{s,\ell}}\{Y_s(\alpha) - \Esp(Y_s(\alpha)|\Fcal_{\ell-1}; U,V)\}$ with respect to $(\Fcal_{\ell})_{\ell\in\llbracket{2,N} \rrbracket}$.
\begin{theorem}[\citep{HaHy14}]\label{tcl:hh}
Suppose that for every $n\in\Nbb$ and $k_n\to\infty$ the random variables $X_{n,1},\ldots,X_{n,k_n}$ are a martingale difference sequence relative to an arbitrary filtration $\Fcal_{n,1}\subset\Fcal_{n,1}\subset\ldots\subset\Fcal_{n,k_n}$. If
\begin{itemize}
\item[1.] $\sum_{i=1}^{k_n}\Esp(X_{n,i}^2|\Fcal_{n,i-1})\rightarrow 1$ in probability,
\item[2.] $\sum_{i=1}^{k_n}\Esp(X_{n,i}^2\Ibb\{|X_{n,i}|>\epsilon\}|\Fcal_{n,i-1})\rightarrow 0$ in probability for every $\epsilon>0$,
\end{itemize}
then $\sum_{i=1}^{k_n}X_{n,i}\overset{\Dcal}{\longrightarrow} \Ncal(0, 1)$.
\end{theorem}
Here $X_{n,i}$ and $\Fcal_{n,i}$ would be $M_{s,\ell}(U,V)/(c_s\sqrt{\Var(F_s)})$  and $\Fcal_{\ell}$ respectively, and we have to verify the two following conditions:
\begin{itemize}
\item[(C1)] $\frac{1}{\Var(N_s)}\sum_{\ell=1}^{N}\Esp(M_{s,\ell}^2(U,V)|\Fcal_{\ell-1}; U,V)\rightarrow\frac{\Var(N_s | U,V)}{\Var(N_s)}$ in probability,
\item[(C2)] $\frac{1}{\Var(N_s)}\sum_{\ell=1}^{N}\Esp\left(M_{s,\ell}^2(U,V)\Ibb\left\{\frac{|M_{s,\ell}(U,V)|}{\sqrt{\Var(N_s)}}>\epsilon\right\}|\Fcal_{\ell-1}; U,V\right)\rightarrow 0$ in probability for every $\epsilon>0$. 
\end{itemize}
Let verify condition (C1). First observe that it follows from properties of martingale differences, meaning variance decomposition, null conditional expectation and conditional orthogonality of differences, that
\begin{align*}
&\Var\left(\sum_{\ell=1}^{N}M_{s,\ell}(U,V)| U,V\right)\\
&=\Esp\left[\Var\left(\sum_{\ell=1}^{N}M_{s,\ell}(U,V)|\Fcal_{\ell-1}; U,V\right)\right]+\Var\left[\Esp\left(\sum_{\ell=1}^{N}M_{s,\ell}(U,V)|\Fcal_{\ell-1};U,V\right)\right]\\
&=\Esp\left[\Var\left(\sum_{\ell=1}^{N}M_{s,\ell}(U,V)|\Fcal_{\ell-1}; U,V\right)\right]\\
&=\Esp\left[\Esp\left(\left(\sum_{\ell=1}^{N}M_{s,\ell}(U,V)\right)^2|\Fcal_{\ell-1}; U,V\right)\right]\\
&=\Esp\left(\sum_{\ell=1}^{N}\Esp\left(M_{s,\ell}^2(U,V)|\Fcal_{\ell-1}; U,V\right)\right),
\end{align*}
and further notice that $\Var\left(\sum_{\ell=1}^{N}M_{s,\ell}(U,V)| U,V\right)=\Var\left(c_sM_s(U,V)| U,V\right)$. 
Since $M_s=L_s-R_s$ (see Section \ref{subsec:Ls}), 
\begin{align*}
&\Esp\left(\frac{1}{\Var(N_s)}\sum_{\ell=1}^{N}\Esp(M_{s,\ell}^2(U,V)|\Fcal_{\ell-1}; U,V)\right)\\
&=\Var\left(c_s\frac{L_s(U,V)-R_s(U,V)}{\sqrt{\Var(N_s)}}| U,V\right)\rightarrow \Var(N_s | U,V)/\Var(N_s),\mbox{ as }n\to\infty,
\end{align*}
in probability and under condition $a+b<2/d_+^{s}$,
because
$\Var\left(c_s\frac{L_s(U,V)}{\sqrt{\Var(N_s)}}| U,V\right)=\frac{\Var(N_s | U,V)}{\Var(N_s)}$ and $\Var\left(c_s\frac{R_s(U,V)}{\sqrt{\Var(N_s)}}| U,V\right)\rightarrow 0$ a.s. under condition $a+b<2/d_+^{s}$ by Lemma \ref{Lemma:Rs}.

Now, let verify condition (C2). 
First, by applying the Cauchy-Schwartz inequality, we get 
\begin{eqnarray*}
&&\frac{1}{\Var(N_s)}\sum_{\ell=1}^{N}\Esp\left(M_{s,\ell}^2(U,V)\Ibb\left\{\frac{|M_{s,\ell}(U,V)|}{\sqrt{\Var(N_s)}}>\epsilon\right\}|\Fcal_{\ell-1}; U,V\right)\\
&\leq &\sum_{\ell=1}^{N}\Esp\left(\frac{M_{s,\ell}^4(U,V)}{\Var(N_s)^2}|\Fcal_{\ell-1}; U,V\right)^{1/2}\times \sum_{\ell=1}^{N}\Pr\left(\frac{|M_{s,\ell}(U,V)|}{\sqrt{\Var(N_s)}}>\epsilon|\Fcal_{\ell-1}; U,V\right)^{1/2},
\end{eqnarray*}
then applying Bienaym\'e-Tchebychev inequality implies that
\begin{eqnarray*}
&&\sum_{\ell=1}^{N}\Esp\left(\frac{M_{s,\ell}^4(U,V)}{\Var(N_s)^2}|\Fcal_{\ell-1}; U,V\right)^{1/2}\times \sum_{\ell=1}^{N}\Pr\left(\frac{|M_{s,\ell}(U,V)|}{\sqrt{\Var(N_s)}}>\epsilon|\Fcal_{\ell-1}; U,V\right)^{1/2}\\
&\leq & \frac{1}{\Var(N_s)^2}\sum_{\ell=1}^{N}\Esp(M_{s,\ell}^4(U,V)|\Fcal_{\ell-1}; U,V)\times \frac{1}{\epsilon^2\Var(N_s)}\sum_{\ell=1}^{N}\Esp(M_{s,\ell}^2(U,V)|\Fcal_{\ell-1}; U,V),
\end{eqnarray*}
and by condition (C1), we get
\begin{eqnarray*}
&&\frac{1}{\Var(N_s)}\sum_{\ell=1}^{N}\Esp\left(M_{s,\ell}^2(U,V)\Ibb\left\{\frac{|M_{s,\ell}(U,V)|}{\sqrt{\Var(N_s)}}>\epsilon\right\}|\Fcal_{\ell-1}; U,V\right)\\
&\leq & \frac{1}{\epsilon^2\Var(N_s)^2}\sum_{\ell=1}^{N}\Esp(M_{s,\ell}^4(U,V)|\Fcal_{\ell-1}; U,V)\times \frac{\Var(N_s | U,V)}{\Var(N_s)}.
\end{eqnarray*}
Then, we use the following notation for expressing $M_{s,\ell}$:
\begin{eqnarray*}
M_{s,\ell}=\sum_{\alpha\in \Qcal_{s,\ell}}\{Y_s(\alpha) - \Esp(Y_s(\alpha)|\Fcal_{\ell-1}; U,V)\}=
N_{s,\ell}-\Esp(N_{s,\ell}|\Fcal_{\ell-1}; U,V).
\end{eqnarray*}
By the binomial formula we thus have
\begin{eqnarray*}
\mathbb{E}(M_{s,\ell}^4(U,V)|\Fcal_{\ell-1}; U,V)&=&
\mathbb{E}\left(\left(N_{s,\ell}-\Esp(N_{s,\ell}|\Fcal_{\ell-1}; U,V)\right)^4|\Fcal_{\ell-1}; U,V\right)\\
&=&\mathbb{E}\left(N_{s,\ell}^4|\Fcal_{\ell-1}; U,V\right)-4\mathbb{E}\left(N_{s,\ell}^3\Esp(N_{s,\ell}|\Fcal_{\ell-1}; U,V)|\Fcal_{\ell-1}; U,V\right)\\
&&\quad +6\mathbb{E}\left(N_{s,\ell}^2\Esp(N_{s,\ell}|\Fcal_{\ell-1}; U,V)^2|\Fcal_{\ell-1}; U,V\right)\\
&&\quad-4\mathbb{E}\left(N_{s,\ell}\Esp(N_{s,\ell}|\Fcal_{\ell-1}; U,V)^3|\Fcal_{\ell-1}; U,V\right)+
\Esp(N_{s,\ell}|\Fcal_{\ell-1}; U,V)^4.
\end{eqnarray*}

Using the same arguments as in the proof of Lemma \ref{Lemma:Rs}, observe that
\begin{eqnarray*}
{E}\left(N_{s,\ell}|\Fcal_{\ell-1}; U,V\right)&\leq&2\sum_{\alpha\in \Qcal_{s,\ell}}\left(\prod_{j\in\alpha^b}\phi_1(U_{k_{\ell}}, V_j)^{A^s( k_{\ell} j)}\right)\Esp\left(\left(\prod_{i\in\alpha^t, j\in \alpha^b\setminus k_{\ell}} G_{ij}^{A^s_{ij}} \right)|U,V\right)\\
&\leq&2\sum_{\alpha\in \Qcal_{s,\ell}}\phibar_s(U_{\alpha^t}, V_{\alpha^b}),
\end{eqnarray*}
and we have,
\begin{eqnarray*}
\mathbb{E}\left(N_{s,\ell}^k|\Fcal_{\ell-1}; U,V\right)=\mathbb{E}\left(N_{s,\ell}|\Fcal_{\ell-1}; U,V\right)+\sum_{t \in \Scal_k(s)}\mathbb{E}\left(N_{t,\ell}|\Fcal_{\ell-1}; U,V\right)+\mathbb{E}\left(N_{s,\ell}|\Fcal_{\ell-1}; U,V\right)^k,
\end{eqnarray*}
where $\Scal_k(s)$ denotes here the set of supermotifs of $s$ which are here formed by $k$ overlapping occurrences of $s$. 
Finally, we get
\begin{eqnarray*}
&&\frac{1}{\Var(N_s)}\sum_{\ell=1}^{N}\Esp\left(M_{s,\ell}^2(U,V)\Ibb\left\{\frac{|M_{s,\ell}(U,V)|}{\sqrt{\Var(N_s)}}>\epsilon\right\}|\Fcal_{\ell-1}; U,V\right)\\
&\leq&\frac{1}{\epsilon^2\Var(N_s)^2}\sum_{\ell=1}^{N}\Esp\left(M_{s,\ell}^4(U,V)|\Fcal_{\ell-1}; U,V\right)\times \frac{\Var(N_s | U,V)}{\Var(N_s)}\\
&\leq&\frac{2}{\epsilon^2\Var(N_s)^2}\sum_{\ell=1}^{N}|\Qcal_{s,\ell}|^4\times\left(\max_{\alpha\in\Qcal_{s,\ell}}\phibar_s(U_{\alpha^t}, V_{\alpha^b})\right)^4\times \frac{\Var(N_s | U,V)}{\Var(N_s)}.
\end{eqnarray*}
Condition (C2) holds since $\Var(N_s)^2=\Theta\left(\max\left(N^{4\pt_s+4\pb_s-4ad_+^{s}-4bd_+^{s}-2},N^{4\pt_s+4\pb_s-2ad_+^{s}-2bd_+^{s}-4}\right)\right)$ by Lemma \ref{Lemma:var:count}, 
$|\Qcal_{s,\ell}|^4=\Theta\left(\ell_t^{4\pt_s-4}\ell_b^{4\pb_s}\right)$ (see the proof of Lemma \eqref{Lemma:Cardinal}), $\phi_s^4=\Theta\left(N^{-4ad_+^s-4bd_+^s}\right)$ by \eqref{order:phis} and $\Var(N_s | U,V)/\Var(N_s)=\Theta(1)$ by Lemma \ref{Lemma:Var:deconditionning}.
\proofend

\subsection{Proof of Lemma \ref{Lemma:Cs}} \label{sec:lemma:cs} 
\proofbegin
Let show that $(\Fbar_s-\phibar_s)/\sqrt{\Var(F_s)}\to 0$ a.s. as $n\to\infty$ under the \BEDD model and condition $a+b<2/d_+^s$ ruling the graph density. 
Recall \eqref{eq:F.bar} the definition of $\Fbar_s$:
\begin{equation*}
 \Fbar_s 
= \frac{\prod_{u = 1}^{\pt_s}  \Gamma_{\dt^s_u} \prod_{v = 1}^{\pb_s}\Lambda_{\db^s_v}}{F_1^{d_+^s}},
\end{equation*}
where $\Gamma_d$ (resp $\Lambda_d$) denote the normalized empirical frequencies of the top (resp bottom) star motif with degree $d$ and $F_1$ the one of the edge.

Let begin with a Taylor expansion of order 1 of $\Fbar_s$ in parameters $(\st,\sb,\phi_1)$ denoting the top star motif, bottom star motif and edge probabilities respectively:
\begin{eqnarray*}
\Fbar_s(\Gamma,\Lambda,F_1)&=&\Fbar_s(\st,\sb,\phi_1)+\partial\Fbar_s(\st,\sb,\phi_1)\left((\Gamma,\Lambda,F_1)-(\st,\sb,\phi_1)\right)+o\left((\Gamma,\Lambda,F_1)-(\st,\sb,\phi_1)\right)\\
&=&\phibar_s+\phibar_s\partial\log(\Fbar_s(\st,\sb,\phi_1))\left((\Gamma,\Lambda,F_1)-(\st,\sb,\phi_1)\right)+o\left((\Gamma,\Lambda,F_1)-(\st,\sb,\phi_1)\right)\\
&=&\phibar_s+\phibar_s\left\{ \sum_{u = 1}^{\pt_s}  \frac{1}{\st_{\dt^s_u}}(\Gamma_{\dt^s_u}-\st_{\dt^s_u}) + \sum_{v = 1}^{\pb_s}\frac{1}{\sb_{\db^s_v}} (\Lambda_{\db^s_v}-\sb_{\db^s_v}) - \frac{\dt_+}{\phi_1}(F_1-\phi_1)\right\}\\
&&+o\left(\Gamma-\st,\Lambda-\sb,F_1-\phi_1)\right).
\end{eqnarray*}
Given the two following observations: i) the asymptotic normality of $(F_s - \phibar_s)/\sqrt{\Var(F_s)}$ holds for any motif $s$, including star motifs, under the \BEDD model and condition $a+b<2/d_+^s$ by Proposition \ref{Prop:Ls}, ii) the empirical frequencies of motifs converge to the expected ones by the law of large numbers, we get 
\begin{align*} 
&\frac{\Fbar_s-\phibar_s}{\sqrt{\Var(F_s)}}\\
&=\sum_{u = 1}^{\pt_s} \Theta\left(\frac{\phi_s}{\st_{\dt^s_u}}\sqrt{\frac{\Var(\Gamma_{\dt^s_u})}{\Var(F_s)}}\right)
+   \sum_{v = 1}^{\pb_s}\Theta\left(\frac{\phi_s}{\sb_{\db^s_v}}\sqrt{\frac{\Var(\Lambda_{\db^s_v})}{\Var(F_s)}}\right)
+  \Theta\left(\frac{\phi_s}{\phi_1}\sqrt{\frac{\Var(F_1)}{\Var(F_s)}}\right)
+o\left(1\right)\\
&= \sum_{u = 1}^{\pt_s}\Theta\left(\frac{\phi_sc_s}{\st_{\dt^s_u}c_{\st}}\sqrt{\frac{\Var(N_{\Gamma_{\dt^s_u}})}{\Var(N_s)}}\right)
+   \sum_{v = 1}^{\pb_s}\Theta\left(\frac{\phi_sc_s}{\sb_{\db^s_v}c_{\sb}}\sqrt{\frac{\Var(N_{\Lambda_{\db^s_v}})}{\Var(N_s)}}\right)
+  \Theta\left(\frac{\phi_sc_s}{\phi_1c_1}\sqrt{\frac{\Var(N_1)}{\Var(N_s)}}\right)
+o\left(1\right).
\end{align*}
Here and only here, $N_{\Gamma_d}$ (resp. $N_{\Lambda_d}$) and $c_{\st}$ (resp. $c_{\sb}$) denote, by abuse of notation, the count of top stars (resp. bottom stars) of degree $d$ and their number of positions in the graph. 
Considering only non-star motifs $s$, according to the orders of magnitude of $c_s$, $\phi_s$ and $\Var(N_s)$ given in \eqref{order:cs}, \eqref{order:phis} and Lemma \ref{Lemma:var:count} respectively, we conclude to $(\Fbar_s - \phibar_s)/\sqrt{\Var(F_s)}\rightarrow 0$ a.s. as $n\to\infty$ because $-2d(a+b)<0$ , with $d=\dt^s_u, \db^s_u$ or 1. 
\proofend

\subsection{Proof of Lemma \ref{Lemma:Var}} \label{sec:lemma:var} 

\proofbegin Let show that $\hat{\Var}(F_s)/\Var(F_s)\rightarrow 1$ a.s., as $n\to\infty$. First, observe that according to \eqref{eq:Var:Ns}, we can write:
$$\Var(N_s) =\sum_{t \in \{s\}\cup\Scal_2(s)} \Esp(N_t)-\Esp(N_s)^2=c_s\phi_s + \sum_{t \in\Scal_2(s)} c_t\phi_t- c_s^2\phi_s^2,$$
where $\Scal_2(s)$ denotes here the set of super-motifs of $s$ which are formed by two overlapping occurrences of $s$.
Then considering  $\hat{\Var}(N_s)$ its plug-in version, meaning $\Fbar_s$ replaces $\phi_s$, we get
$$\hat{\Var}(N_s)-\Var(N_s) =c_s (\Fbar_s-\phi_s) + \sum_{t \in \Scal_2(s)} c_t (\Fbar_t-\phi_t)-c_s^2( \Fbar_s^2-\phi_s^2).$$
Now we use Lemma \ref{Lemma:Cs} stating that, under the \BEDD model and condition $a+b<2/d_+^{s}$, $\Fbar_s-\phi_s=o(\sqrt{\Var(F_s)})$ for all motif $s$ and the continuous mapping theorem, to obtain that
\begin{eqnarray}\label{eq:proof:l2}
\hat{\Var}(N_s)-\Var(N_s) &=&c_so(\sqrt{\Var(F_s)}) + \sum_{t \Scal_2(s)} c_to(\sqrt{\Var(F_t)})-c_s^2o(\Var(F_s)\nonumber)\\
&=&o(\sqrt{\Var(N_s)})  + \sum_{t \in\Scal_2(s)} o(\sqrt{\Var(N_t)})-o(\Var(N_s)).
\end{eqnarray}
Let discuss now the order of $(\hat{\Var}(N_s)-\Var(N_s))/\Var(N_s)$. The first and last terms of \eqref{eq:proof:l2} divided by $\Var(N_s)$ obviously vanish. When $t\in\Scal_2(s)$, we refer to the order of magnitude of $\Var(N_s)$ given in Lemma \ref{Lemma:var:count} and its proof (see ($i$)-($ii$)-($iii$)) to get that $\sqrt{\Var(N_t)}/\Var(N_s)$ vanishes under condition $a+b<(\pt_s +\pb_s )/d_+^{s}$.
We can finally conclude to $\hat{\Var}(F_s)/\Var(F_s)\rightarrow 1$ a.s., as $n\to\infty$ under condition of Theorem \ref{Thm:AsympNorm}.

\proofend

\begin{figure}[ht]
 \begin{center}
  $\begin{array}{l|c}
    s & \quad 1 \quad \\ \hline 
    & \\ 
    & \input{Figs/1-bipartite-edge} \\ \\
    c_s & \displaystyle{m n} \\ \\ \\
    \phibar_s & \displaystyle{\phi_1} \\ 
  \end{array}$ \qquad
  $\begin{array}{l|cc}
    s & \quad 2 \quad & \quad 3 \quad \\ \hline 
    & \\ 
    & \input{Figs/2-bipartite-lambda}
    & \input{Figs/3-bipartite-V} \\ \\
    c_s 
    & \displaystyle{m \binom{n}2}
    & \displaystyle{n \binom{m}2} \\ \\
    \phibar_s 
    & \displaystyle{\gamma_2} 
    & \displaystyle{\lambda_2} 
  \end{array}$ 
  \\ \bigskip
  $\begin{array}{l|cccc}
    s & \quad 4 \quad & \quad 5 \quad & \quad 6 \quad & \quad 7 \quad \\ \hline 
    & & &  \\ 
    & \input{Figs/4-bipartite-star3up}
    & \input{Figs/5-bipartite-N}
    & \input{Figs/6-bipartite-pap} 
    & \input{Figs/7-bipartite-star3down} \\ \\
    c_s 
    & \displaystyle{n \binom{m}3}
    & \displaystyle{4 \binom{m}2 \binom{n}2}
    & \displaystyle{\binom{m}2 \binom{n}2}
    & \displaystyle{m \binom{n}3} \\ \\
    \phibar_s 
    & \displaystyle{\lambda_3} 
    & \displaystyle{\gamma_2 \lambda_2 / \phi_1} 
    & \displaystyle{\gamma_2^2 \lambda_2^2 / \phi_1^4} 
    & \displaystyle{\gamma_3} 
  \end{array}$ 
  \\ \bigskip
  $\begin{array}{l|cccccccccc}
    s & \quad 8 \quad & \quad 9 \quad & \quad 10 \quad & \quad 11 \quad & \quad 12 \quad \\ \hline 
    & & & & & \\ 
    & \input{Figs/8-bipartite-star4up} 
    & \input{Figs/9-bipartite-star3upN}
    & \input{Figs/10-bipartite-W}
    & \input{Figs/11-bipartite-VX}
    & \input{Figs/12-bipartite-XXup} \\ \\
    c_s 
    & \displaystyle{n \binom{m}{4}}
    & \displaystyle{6\binom{m}{3}\binom{n}{2}}
    & \displaystyle{6\binom{m}{3}\binom{n}{2}}
    & \displaystyle{6\binom{m}{3}\binom{n}{2}}
    & \displaystyle{\binom{m}{3}\binom{n}{2}} \\ \\
    \phibar_s 
    & \displaystyle{\lambda_4} 
    & \displaystyle{\gamma_2 \lambda_3 / \phi_1} 
    & \displaystyle{\gamma_2 \lambda_2^2 / \phi_1^2} 
    & \displaystyle{\gamma_2^2 \lambda_2 \lambda_3 / \phi_1^4} 
    & \displaystyle{\gamma_3^2 \lambda_3^2 / \phi_1^6} \\
    \multicolumn{2}{} \\ \\
    s & \quad 13 \quad & \quad 14 \quad & \quad 15 \quad & \quad 16 \quad & \quad 17 \quad \\ \hline 
    & & & & & \\ 
    & \input{Figs/13-bipartite-Istar3down}
    & \input{Figs/14-bipartite-M}
    & \input{Figs/15-bipartite-Xlambda}
    & \input{Figs/16-bipartite-XXdown}
    & \input{Figs/17-bipartite-star4down} \\ \\
    c_s 
    & \displaystyle{6\binom{m}{2}\binom{n}{3}}
    & \displaystyle{6\binom{m}{2}\binom{n}{3}}
    & \displaystyle{6\binom{m}{2}\binom{n}{3}}
    & \displaystyle{\binom{m}{2}\binom{n}{3}}
    & \displaystyle{m \binom{n}{4}} \\ \\
    \phibar_s 
    & \displaystyle{\gamma_3 \lambda_2 / \phi_1} 
    & \displaystyle{\gamma_2^2 \lambda_2 / \phi_1^2} 
    & \displaystyle{\gamma_2 \gamma_3 \lambda_2^2 / \phi_1^4} 
    & \displaystyle{\gamma_3^2 \lambda_2^3 / \phi_1^6} 
    & \displaystyle{\gamma_4} \\
  \end{array}$
  \caption{Bipartite motifs of size 2, 3, 4 and 5 as given in \cite{SCB19}. \label{fig:motifs1-5}}
 \end{center}
\end{figure}

\begin{figure}[ht]
 \begin{center}
  $\begin{array}{l|cccccc}
    s & \quad 18 \quad & \quad 19 \quad & \quad 20 \quad & \quad 21 \quad & \quad 22 \quad & \quad 23 \quad \\ \hline 
    & \\ 
    & \input{Figs/18-bipartite-star5down}
    & \input{Figs/19-bipartite}
    & \input{Figs/20-bipartite}
    & \input{Figs/21-bipartite}
    & \input{Figs/22-bipartite} 
    & \input{Figs/23-bipartite} \\ 
    \multicolumn{2}{} \\ \\
    s & \quad 24 \quad & \quad 25 \quad & \quad 26 \quad & \quad 27 \quad & \quad 28 \quad & \quad 29 \quad \\ \hline 
    & \\ 
    & \input{Figs/24-bipartite}
    & \input{Figs/25-bipartite}
    & \input{Figs/26-bipartite}
    & \input{Figs/27-bipartite}
    & \input{Figs/28-bipartite} 
    & \input{Figs/29-bipartite} \\ 
    \multicolumn{2}{} \\ \\
    s & \quad 30 \quad & \quad 31 \quad & \quad 32 \quad & \quad 33 \quad & \quad 34 \quad & \quad 35 \quad \\ \hline 
    & \\ 
    & \input{Figs/30-bipartite}
    & \input{Figs/31-bipartite}
    & \input{Figs/32-bipartite}
    & \input{Figs/33-bipartite}
    & \input{Figs/34-bipartite} 
    & \input{Figs/35-bipartite} \\ 
    \multicolumn{2}{} \\ \\
    s & \quad 36 \quad & \quad 37 \quad & \quad 38 \quad & \quad 39 \quad & \quad 40 \quad & \quad 41 \quad \\ \hline 
    & \\ 
    & \input{Figs/36-bipartite}
    & \input{Figs/37-bipartite}
    & \input{Figs/38-bipartite}
    & \input{Figs/39-bipartite}
    & \input{Figs/40-bipartite} 
    & \input{Figs/41-bipartite} \\ 
    \multicolumn{2}{} \\ \\
    s & \quad 42 \quad & \quad 43 \quad & \quad 44 \quad \\ \hline 
    & \\ 
    & \input{Figs/42-bipartite}
    & \input{Figs/43-bipartite}
    & \input{Figs/44-bipartite-star5up}
  \end{array}$
  \caption{Bipartite motifs of size 6 as given in \cite{SCB19}.} \label{fig:motifs6} 
 \end{center}
\end{figure}

\bibliographystyle{plainnat}
\bibliography{Biblio}


\end{document}

%% file: Commands.tex
\usepackage{enumitem,oldgerm}
\usepackage{stmaryrd}
\usepackage{stix}
\usepackage{amsmath}
\usepackage{amssymb}
\usepackage{amsthm}
\usepackage{xcolor}
\usepackage{xspace}
\usepackage{enumerate}
\usepackage{hyperref}
\usepackage{graphicx}
\usepackage{tikz}

\newcommand{\Bcal}{\mathcal{B}}
\newcommand{\Bias}{\mathbb{B}}
\newcommand{\Cov}{\mathbb{C}\text{ov}}

\newcommand{\Dcal}{\mathcal{D}}

\newcommand{\Esp}{\mathbb{E}}
\newcommand{\Ecal}{\mathcal{E}}

\newcommand{\Fbar}{\overline{F}}
\newcommand{\Fcal}{\mathcal{F}}

\newcommand{\Gcal}{\mathcal{G}}

\newcommand{\Ibb}{\mathbb{I}}

\newcommand{\Ncal}{\mathcal{N}}
\newcommand{\Nbb}{\mathbb{N}}
\newcommand{\Ocal}{\mathcal{O}}
\newcommand{\Pcal}{\mathcal{P}}
\renewcommand{\Pr}{\mathbb{P}}
\newcommand{\Qcal}{\mathcal{Q}}
\newcommand{\phibar}{\overline{\phi}}
\newcommand{\Rbb}{\mathbb{R}}
\newcommand{\Scal}{\mathcal{S}}

\newcommand{\Ucal}{\mathcal{U}}
\newcommand{\Vcal}{\mathcal{V}}
\newcommand{\Var}{\mathbb{V}}

\newcommand{\BEDD}{\text{B-EDD}\xspace}
\newcommand{\gt}{g} \newcommand{\gb}{h}
\newcommand{\mut}{\mu_{\gt}} \newcommand{\mub}{\mu_{\gb}}
\newcommand{\pt}{p} \newcommand{\pb}{q}
\newcommand{\dt}{d} \newcommand{\db}{e}

\newcommand{\st}{\gamma} 
\renewcommand{\sb}{\lambda}
\newcommand{\St}{\Gamma} 
\newcommand{\Sb}{\Lambda}

\renewcommand{\d}{\text{d}\xspace}

\usepackage{tikz}
\newcommand{\nodesize}{2em}
\newcommand{\edgeunit}{1.1*\nodesize}
\tikzstyle{node}=[inner sep=0]
\tikzstyle{node1}=[inner sep=0]
\tikzstyle{edge}=[-, line width=.75pt]
\tikzstyle{edgedot}=[-, dotted, line width=.75pt]
\newcommand{\topnode}{$\mdwhtcircle$}
\newcommand{\bottomnode}{$\mdwhtsquare$}
\newcommand{\topnodeter}{$\circlebottomhalfblack$}
\newcommand{\bottomnodeter}{$\squareurblack$}
\newcommand{\topnodebis}{$\mdblkcircle$}
\newcommand{\bottomnodebis}{$\mdblksquare$}

\newtheorem{lemma}{Lemma}
\newtheorem{proposition}{Proposition}
\newtheorem{remark}{Remark}
\newtheorem{theorem}{Theorem}

\newcommand{\proofbegin}{\noindent{\sl Proof.}~}
\newcommand{\proofend}{$\blacksquare$\bigskip}